\documentclass[reqno,twoside,11pt]{amsart}
\usepackage{amsmath, amsfonts, amssymb, esint, amsthm, nicefrac, bbm,
  dsfont, comment}
\usepackage{epsfig, graphicx,xcolor}
\usepackage[colorlinks=true, linkcolor=blue, citecolor=blue, urlcolor=blue]{hyperref}

\newcommand{\dx}{~\mathrm{d}x}
\newcommand{\sym}{\mathrm{sym}}

\newcommand{\Id}{\mathrm{Id}}

\newcommand{\R}{\mathbb{R}}
\def\endproof{\hspace*{\fill}\mbox{\ \rule{.1in}{.1in}}\medskip }
\newcommand*{\dbar}[1]{\bar{\bar{#1}}}
\newcommand*{\tbar}[1]{\bar{\dbar{#1}}}

\makeatletter
\renewcommand{\part}{\@startsection{part}{0}\z@
  {2\linespacing\@plus\linespacing}{\linespacing}
  {\normalfont\normalsize\bfseries\centering
 \let\@secnumfont\bfseries}}
\makeatother

\makeatletter
\renewcommand{\l@part}{%
  \@tocline{-1}{0.5em}{0pt}{}{\bfseries}}
\makeatother

\numberwithin{equation}{section}
\theoremstyle{plain}

\newtheorem*{theorem*}{Theorem}
\newtheorem{theorem}{Theorem}[section]
\newtheorem{lemma}[theorem]{Lemma}
\newtheorem{corollary}[theorem]{Corollary}
\newtheorem{proposition}[theorem]{Proposition}

\theoremstyle{definition}
\newtheorem{remark}[theorem]{Remark}
\newtheorem{example}[theorem]{Example}

\begin{document}
\title [Flexibility for (MA) implies flexibility for (II)]
{Flexibility for the Monge-Amp\`ere system implies flexibility for isometric immersions \\ into Riemannian manifolds}
\author{Dominik Inauen and Marta Lewicka}
\address{D.I.: Universit\'e de Fribourg, CH-1700, Fribourg, Switzerland}
\address{M.L.: University of Pittsburgh, Department of Mathematics, 
139 University Place, Pittsburgh, PA 15260, USA}
\email{dominik.inauen@unifr.ch, lewicka@pitt.edu}

\date{\today}
\thanks{
M.L. was partially supported by NSF grant DMS-2006439.
}

\begin{abstract} 
We prove that any flexibility result for the Monge-Amp\`ere system
implies the corresponding flexibility result for the isometric immersion system, with the same H\"older regularity. We treat the general isometric immersion system, in which both the source and the ambient metrics are arbitrary H\"older regular, in arbitrary dimension and codimension. Combining this implication with existing  results for the Monge-Amp\`ere system, we obtain new $\mathcal{C}^{1,\alpha}$ flexibility statement for isometric immersions with arbitrary source and ambient metrics. To our knowledge, these are the first such results for non-Euclidean ambient metrics.
\end{abstract} 

\maketitle
\tableofcontents

\section{Introduction}\label{sec_intro}

The purpose of this paper is to clarify the relationship between the
Monge-Amp\`ere system (\ref{MA}), through its weak formulation, given
by the von K\'arm\'an system (\ref{VK}), and the isometric immersion
system (\ref{II}), in the context of flexibility results obtained via convex
integration. We consider immersions of $(\omega\subset\R^d, g)$ into $(\mathcal{U}\subset \R^{d+k},\bar g)$ and allow both the source metric $g$ and the ambient metric $\bar g$ to have arbitrary H\"older regularity. Numerous flexibility results have been established for (\ref{VK}) and for (\ref{II}) with the Euclidean target (see Table \ref{taula1} and Table \ref{taula2}), where the corresponding analyses have been developing largely in
parallel, despite the evident similarities in their
proof techniques. In this paper, we identify structural properties of the stage constructions, which make the passage between the two systems systematic: we prove that any appropriately formulated stage for (\ref{VK}) in arbitrary dimension $d$ and codimension $k$, yields a stage for (\ref{II}), with the same defect decay parameter $S$ the and derivative growth parameter $J$. Consequently, this stage transfer preserves the H\"older flexibility exponent $\alpha$.

\medskip

\noindent The second purpose of our paper is to develop the Nash–Kuiper iteration and the pre-stage constructions, thereby turning the stage transfer into flexibility results for the system (\ref{II}) with general $g,\bar g$. We obtain uniform approximation of strictly short immersions by $\mathcal{C}^{1,\alpha}$ isometric immersions, with $\alpha$ determined by the stage parameters $J, S$ and the H\"older regularity of the two metrics $g,\bar g$. While the classical Nash-Kuiper theory provides approximation by $\mathcal{C}^{1}$ isometric immersions, its quantitative $\mathcal{C}^{1,\alpha}$ counterpart has been developed primarily for the Euclidean target $\bar g=\Id_{d+k}$. To the best of the authors' knowledge, such H\"older flexibility results have not been previously established in this generality.

\medskip

\subsection{Three PDE systems: (\ref{II}), (\ref{MA}) and (\ref{VK})} \label{sec_systems}
Let $\omega\subset\R^d$, $\mathcal U\subset\R^{d+k}$ be open and bounded sets, where $d, k\geq 1$. Given Riemannian metrics $g$ and $\bar g$, posed  on $\bar \omega$ and $\bar{ \mathcal U}$ respectively, the isometric immersion problem seeks a map $u$ satisfying: 
\begin{equation}\tag{II}\label{II}
\left\{\begin{split}  
& (\nabla u)^T(\bar g\circ u)\nabla u =g \quad \mbox{in }\; \omega,
\\ & \mbox{for } \; u:\omega\to\mathcal{U}, \quad g:\omega\to\R^{d\times d}_{\sym, >},\quad \bar g:\mathcal{U}\to\R^{(d+k)\times (d+k)}_{\sym, >}.
\end{split}\right.
\end{equation}
The central question considered here is whether, and with what regularity, every strictly short immersion can be uniformly approximated by solutions of \eqref{II}.  Here, strict shortness means:
\[g-(\nabla u)^T(\bar g\circ u)\nabla u>0,\]
in the sense of quadratic forms. For general ambient metrics, the strict shortness is an essential condition and cannot be relaxed to mere shortness, see Example \ref{counter_ex}. 

\medskip

\noindent Recall that the Monge-Amp\`ere system, introduced in \cite{lew_conv},
seeks a solution $v$ to:
\begin{equation}\tag{MA}\label{MA}
\left\{\begin{split}
& \mathfrak{Det} \,\nabla^2v \doteq \big[\langle \partial_{is}
v, \partial_{jt}v\rangle -\langle \partial_{it}v, \partial_{js}v\rangle\big]_{i,j,
  s,t:1\ldots d}= F\quad\mbox{ in }\;\omega \subset\R^d,      
\\ & \mbox{for } \; v: \omega \to \R^k, \quad F:\omega\to\R^{d^4},
\end{split}\right.
\end{equation}
where $F$ is a given table-valued field, satisfying
symmetries necessary for the resolvability of (\ref{MA}).
The system (\ref{MA}) is the multi-dimensional version of the Monge-Amp\`ere equation:
\begin{equation}\label{MAeq}
\left\{\begin{split}
& \det\nabla^2 v  \doteq \partial_{11}
v  \partial_{22} v - (\partial_{12} v)^2 = f\quad\mbox{ in }\;\omega\subset\R^2,
\\ & \mbox{for } \; v:\omega\to \R,  \quad f:\omega\to\R,
\end{split}\right.
\end{equation} 
arising from the prescribed curvature problem. Indeed, the family of Riemannian metrics on $\omega$, generated by immersions $\{u^\epsilon=(id_d, \epsilon
v)\}_{\epsilon\to 0}$, have their Riemann curvatures satisfy:
\begin{equation*}
{Riem}\big((\nabla u^\epsilon)^T\nabla u^\epsilon\big)  =
{Riem}\big(\mbox{Id}_d + \epsilon^2(\nabla v)^T\nabla v\big) 
= \epsilon^2 \mathfrak{Det} \,\nabla^2v + o(\epsilon^2).
\end{equation*}
When $d=2$, $k=1$ as in (\ref{MAeq}), the above is consistent with the
formula for the Gaussian curvature of surfaces described as the 
graph of $\{\epsilon v\}_{\epsilon\to 0}$, namely:
$$\kappa(\epsilon v) = \frac{\epsilon^2 \det\nabla^2 v}{(1+\epsilon^2|\nabla
  v|^2)^2} = \epsilon^2 \det\nabla^2v + o(\epsilon^2).$$
Thus, (\ref{MAeq}) and (\ref{MA}) appear as the leading order terms in the prescribed curvature equations. 

\medskip

\noindent The second connection between \eqref{II} and \eqref{MA} is provided by the von K\'arm\'an system:
\begin{equation}\tag{VK}\label{VK}
\left\{\begin{split}
& \frac{1}{2}(\nabla v)^T\nabla v + \sym \nabla w = A\quad\mbox{ in
}\;\omega \subset\R^d,
\\ & \mbox{for } \; v:\omega\to \R^k, \quad w:\omega\to\R^d,\quad
A:\omega\to\R^{d\times d}_\sym.
\end{split}\right.
\end{equation}
For an appropriate matrix field $A$ determined by $F$, this system yields a weak formulation for \eqref{MA}. It however also arises by expanding the isometric immersion system (\ref{II}) with $\bar g =\mathrm{Id}_{d+k}$ around the flat immersion. Indeed, setting 
$\bar u^\epsilon = ({id}_d+\epsilon^2w, \epsilon v)$ gives:
$$(\nabla \bar u^\epsilon)^T\nabla \bar u^\epsilon = \mbox{Id}_d +
\epsilon^2 \big((\nabla v)^T\nabla v + 2\, \sym\nabla w\big) + o(\epsilon^2),$$
where matching the $\epsilon^2$ terms in the above expansion and in the prescribed metrics expansion $\Id_d+2\epsilon^2 A$, is exactly \eqref{VK}.
When $d=2$, $k=1$, the above derivation is consistent with 
the left hand side of (\ref{VK}) appearing in the stretching energy of
a thin film, subject to the out-of-plane displacement $v$ and the
in-plane displacement $w$,  featured  ubiquitously in the studies of
prestrained elastic thin films \cite{lewlu, LeOchPa, lew_book}.  
As we see, the three problems (\ref{MA}), (\ref{VK}) and
(\ref{II}) and their expansions are intrinsically related.

\subsection{Main results of this paper} \label{sec_results}
Given the above mentioned connections and the historical account in section \ref{sec_histo} below, it is perhaps not so surprising that ideas have been repeatedly transferred between corresponding convex integration theories.
The first results \cite{lewpak_MA} for \eqref{MA} and \eqref{VK} were obtained by adapting constructions developed for \eqref{II}. More recently, the direction of influence has reversed. The exponent $\alpha_0=1/3$ in codimension $k=1$, was first achieved \cite{CHI} for (\ref{MA}) and the subsequent work \cite{CHI2} introduced a new corrugation scheme for \eqref{II} modeled on the defect removal mechanism of \eqref{VK}. Later developments \cite{lew_full2d2k,In2026,LL_last} continue to build on this viewpoint. 
Despite these similarities, the connections remained mainly heuristic, with no general result showing that flexibility for (\ref{MA}) directly implies flexibility for (\ref{II}). 

\medskip

\noindent In the present paper we primarily establish exactly such an implication. 
Our main results concern the stage constructions \eqref{STA-VK} and \eqref{STA-II} of the respective convex integration schemes for systems \eqref{VK} and \eqref{II}, rather than just the flexibility statements. A stage reduces the defect by adding perturbations, while controlling their derivative growth. We introduce an augmented stage formulation \eqref{STA0} for (\ref{VK}) and show that:
$$\mbox{(\ref{STA-VK0})}\Rightarrow\mbox{(\ref{STA0})}\Rightarrow\mbox{(\ref{STA-II})},\quad \mbox{(\ref{STA0})}\Rightarrow\mbox{(\ref{STA-VK})},$$
with the same parameters $J,S$ measuring the derivative growth and defect reduction. These implications are proved in Theorems \ref{prop_STA0_to_STAII}, \ref{prop_STA0_to_STAVK}, \ref{prop_STA1_to_STA0}.
Formulation \eqref{STA-VK0} slightly strengthens (\ref{STA-VK}) by encoding a structural property of the von K\'arm\'an constructions: they use estimates on the initial Hessian of the out-of-plane displacement $v$ but not its Hessian structure. Therefore, (\ref{STA-VK}) remains correct when $\nabla^2v$ is replaced by any symmetric field satisfying the same bounds, here identified as the second fundamental form relative to $\bar g$. In Examples \ref{ex_first}, and \ref{ex_last}, we show that stage constructions in \cite{lew_conv} and \cite{CHIL} indeed satisfy this formulation.   

\medskip

\noindent The mechanism for the transfer of a von K\'arm\'an stage to an isometric immersion stage is most clearly seen when the ambient metric $\bar g$ is Euclidean. We use the fields $(v,w)$ as normal and tangential perturbations to the given immersion and observe that the resulting difference in the induced metric consists of the von K\'arm\'an content of $(v,w)$, a lower order error term, and another term which is linear in $v$ and whose coefficients are given by the shape operator of the immersion. Incorporating this term in \eqref{STA0} yields the desired transfer, with a comparably short argument. 
The same strategy applies to a general $\bar g$, whose H\"older regularity requires now a regularization, introducing new error terms and making the analysis much more involved. 

\medskip

\noindent The second contribution of the paper is to complete the analysis described above by providing an appropriate \emph{pre-stage} for the general system \eqref{II}, which constructs from a strictly short immersion, a suitable referential immersion, upon which the stage can be iterated in a {\em Nash--Kuiper scheme} to yield Theorem \ref{th_final_II}. 
In keeping with the approach of this paper, we deduce the pre-stage from basic $\mathcal{C}^1$ approximation result  (\ref{STA-INI0}) for (\ref{VK}) with zero initial displacements, and the Nash–Kuiper scheme directly from  (\ref{STA-II}). The prepared referential immersion also serves as the initial iterate, avoiding a separate initial-stage construction \cite{lew_full2d2k, CHI} with a prescribed scaling relation between the approximation error and the blow-up of second derivatives. This provides a reusable framework, independent of the details of their underlying stage proofs. Consequently, we obtain that if \eqref{STA0} holds with parameters $J,S$, then for any Riemannian metric $g\in \mathcal{C}^{r,\beta}$ posed on a contractible domain $\omega\subset\R^d$ and any Riemannian metric $\bar g\in \mathcal{C}^{\bar r,\bar \beta}$ posed on $\R^{d+k},$ with $r+\beta,\bar r+\bar \beta>0$, every strictly short immersion $u\in\mathcal{C}^1(\omega, \R^{d+k})$ can be uniformly approximated by $\mathcal{C}^{1,\alpha}$ solutions of \eqref{II}, with any:
\[\alpha<\min\left\{ \frac{1}{1+2J/S},\frac{r+\beta}{2},\frac{\bar r+\bar \beta}{2}\right\}\,.\]
In particular, when $g,\bar g\in \mathcal{C}^2$, the stage constructions yielding {\em full flexibility} of $(d,k)$-dimensional Monge-Amp\`ere system imply existence of $\mathcal{C}^{1,\alpha}$ isometric immersions, for every $\alpha<1$.

\subsection{Prior developments in convex integration for systems (\ref{II}), (\ref{MA}), (\ref{VK})} \label{sec_histo}
The flexibility theory for the isometric immersion system \eqref{II} originated in the celebrated works of Nash \cite{Nash1, Nash2} and Kuiper \cite{Kuiper}, who showed that any short immersion of a $d$-dimensional Riemannian manifold into $(\R^{d+1},\Id_{d+1})$ can be uniformly approximated by $\mathcal{C}^1$ solutions of \eqref{II}. These solutions arise as limits of an iteratively constructed sequence of short immersions, where the single building block of the iteration is called a \emph{stage}. The process, nowadays referred to as convex integration, remains the foundation of all subsequent flexibility results. The corresponding Nash--Kuiper theorem for strictly short immersions into  general Riemannian manifolds is stated in  \cite[\S 1.3.18]{GR}, see also the discussion in \cite[sections 3.2-3.3]{Gromov-Nash}.

\medskip

\noindent The associated regularity question, namely determining the largest H\"older exponent $\alpha_0$ such that flexibility persists in $\mathcal{C}^{1,\alpha}$ for any $\alpha<\alpha_0$, has been a central topic of research in recent years, mostly in the setting of the Euclidean ambient metric. Starting with the work \cite{CDS} of Conti, De Lellis and Sz\'ekelyhidi which rigorously confirmed claims of Borisov \cite{Borisov1965}, and continuing with \cite{DIS1/5,CHI2,In2026}, the attainable regularity threshold in codimension $k=1$ has been improved progressively by increasingly efficient convex integration schemes; see in particular the recent scheme of Inauen \cite{In2026} which yields the best regularity exponent currently known. In higher codimension $k>1$, substantial progress has likewise been achieved through the development of more elaborate constructions that exploit perturbations in the additional normal directions and their interactions, allowing for higher regularity thresholds and, in certain $(d,k)$ regimes, for full flexibility \cite{Kallen,CaoSz2025,Zhao,lew_full2d2k,LL_last}. These results are summarized in Table \ref{taula1}. 

\medskip

\noindent Explicit convex integration constructions and numerical visualizations of $\mathcal{C}^1$ isometric embeddings have been developed by Borrelli and collaborators, for flat tori \cite{B1}, reduced spheres \cite{B2} and the hyperbolic plane \cite{B3}, giving a geometric picture of the corrugations. We also mention complementary approaches to negatively curved surfaces that include the well-posedness and stability theory \cite{V1}, and the discrete differential-geometric constructions of $\mathcal{C}^{1,1}$ surfaces  \cite{V2}.

\begin{table}[ht]
\centering
\begin{tabular}{||c|c|c||c||}
\hline
{dimension $d$} & {codimension $k$} & {exponent $\alpha_0$} &  {reference} \\ \hline
$d$ & $1$ & $0$ & \cite{Nash1, Nash2, Kuiper} \\ \hline  
$d$ & $2d_*$        & $1$           & \cite{Kallen}      \\ \hline  
$d$ & $1$        & ${1}/({1+2d_*})$           & \cite{CDS}      \\ \hline
$2$ & $1$        & ${1}/5$           & \cite{DIS1/5}      \\ \hline
$d$ & $d$        & ${1}/({2+d})$           & \cite{CaoSz2025}      \\ \hline    
$d$ & $1$        & ${1}/(1+d^2-d)$           & \cite{CHI2}  \\
\hline 
$d$  &  $d$ & $1/(1+d)$   &    \cite{Zhao}\\
\hline
$2$ & $2$        & $1$                     & \cite{lew_full2d2k} \\ \hline
$d$ & $1$        & ${1}/({2d-1})$           & \cite{In2026}  \\ \hline  
$d$ & $k$        & ${1}/({1+ 2d_*/k})$           & \cite{LL_last}     \\ \hline
\end{tabular}
\vspace{2mm}
\caption{Development of convex integration results for
(\ref{II}) with $g\in\mathcal{C}^2$ and $\bar g=\Id_{d+k}$.}
\vspace{-4mm} \label{taula1}
\end{table} 

\medskip

\noindent The present paper concerns a Nash-Kuiper theorem for \eqref{II} for $\mathcal{C}^{1,\alpha}$ isometric immersions, with general dimension $d$ and codimension $k$, and with prescribed H\"older regularities of both source $g$ and target $\bar g$ metrics.
For non-Euclidean $\bar g$, flexibility has been studied in contact geometry by Ried \cite{Ried}, building on the scheme of \cite{CDS} to obtain isometries into contact manifolds satisfying an additional differential constraint. We point out that this constraint poses an extra layer of difficulty with respect to the mere problem (\ref{II}), which the result addresses, however it does not yield the flexibility in the generality considered here.
At higher regularity, the case of non-Euclidean  $\bar g $ has also been studied in connection with the Weyl problem; see for instance the existence and regularity theory of \cite{Lu} under curvature assumptions, or the perturbative existence theorem of \cite{Li} for strictly convex surfaces in warped product spaces.  

\medskip

\noindent Convex integration for the Monge-Amp\`ere equation \eqref{MA} was initiated by Lewicka and Pakzad \cite{lewpak_MA}, who adapted the iteration process of \cite{CDS} to obtain flexibility in dimension $d=2$ and codimension $k=1$. Their H\"older exponent was subsequently improved by Cao and Sz\'ekelyhidi \cite{CS}, paralleling earlier development for (\ref{II}). A general theory in arbitrary dimension and codimension was introduced by Lewicka in \cite{lew_conv}. Several subsequent works \cite{lew_improved,lew_improved2, CHI,in_lew,in_lew2} improved the achievable regularity in dimension $d=2$,  culminating in the best so far known threshold \cite{CHI} for $k=1$, and in the proof of full flexibility \cite{in_lew2} for $k=2$. 
Recently, full flexibility has been shown \cite{CHIL} for arbitrary $d$, with $k=d_*-d+1$. These results are summarized in Table \ref{taula2}.

\begin{table}[ht]
\centering
\begin{tabular}{||c|c|c||c||}
\hline
{dimension $d$} & {codimension $k$} & {exponent $\alpha_0$} &  {reference} \\ \hline
$2$ & $1$ & $1/7$ & \cite{lewpak_MA} \\ \hline  
$2$ & $1$        & $1/5$           & \cite{CS}      \\ \hline  
$d$ & $k$        & ${1}/({1+2d_*/k})$           & \cite{lew_conv}      \\ \hline
$d$ & $2d_*$        & ${1}$           & \cite{lew_conv}  \\ \hline
$2$ & $k$        & $1/(1+4/k)$           & \cite{lew_improved}      \\ \hline  
$2$ & $1$        & ${1}/3$           & \cite{CHI}      \\ \hline    
$2$ & $k$        & $(2^k-1)/(2^{k+1}-1)$           & \cite{lew_improved2}  \\ \hline
$2$ & $4$        & $1$           & \cite{lew_improved2}  \\ \hline  
$2$ & $3$        & $1-1/{\sqrt{5}}$     & \cite{in_lew} \\ \hline
$2$ & $2$        & $1$           & \cite{in_lew2}  \\ \hline  
$d$ & $d_*-d+1$        & $1$           & \cite{CHIL}      \\ \hline  
\end{tabular}
\vspace{2mm}
\caption{Development of convex integration results for
(\ref{VK}) with  $A\in\mathcal{C}^2$.}
\vspace{-4mm} \label{taula2}
\end{table}

\subsection{Organization of the paper}

In section \ref{sec_defi}, we introduce the three stage formulations
(\ref{STA0}), (\ref{STA-II}) and (\ref{STA-VK}), whose relationship
is the first main subject of this paper. Section \ref{sec_preli} collects
several preliminary results, including
convolution estimates, decomposition into primitive metrics, and
the existence and propagation of normal frames along immersions.

\smallskip

\noindent Part 1 of the paper consists of sections \ref{seci1} - \ref{sec_NK_II3}.
In section \ref{seci1}, we prove that (\ref{STA0}) implies
(\ref{STA-II}) for general source and ambient metrics $g,\bar g$. In section \ref{sec_NK_II} we settle the Nash-Kuiper iteration scheme in Proposition \ref{thm_NK_II}, and in section \ref{sec_NK_II2} we construct the pre-stage needed to prepare a strictly short immersion that serves both as the referential immersion and the starting point for the  Nash-Kuiper iteration. In section \ref{sec_NK_II3} we complete the proofs of the flexibility result for \eqref{II}, relying only on the validity of (\ref{STA0}). This fulfills the second main objective of this paper.
We also present an example showing that, in general, the strict shortness  assumption cannot be replaced by the mere shortness.
Part 2 consists of a single section \ref{seci2}, in which we establish
the implication from (\ref{STA0}) to (\ref{STA-VK}). We also recall
the appropriate Nash-Kuiper iteration scheme in Proposition \ref{th_NK_VK}, yielding the flexibility result for the system (\ref{MA}).
Part 3 consists of sections \ref{seci3}, \ref{sec_exa}. In section \ref{seci3}, we introduce the fourth stage formulation, namely (\ref{STA-VK0}) for \eqref{VK}, which incorporates additional structural assumptions on the construction, and prove that it implies
the base formulation (\ref{STA0}). In section \ref{sec_exa} we then illustrate the verification of (\ref{STA-VK0}) through two representative examples, namely show that the minor modifications of constructions in \cite{lew_conv} and \cite{CHIL} satisfy the structural properties required in (\ref{STA-VK0}).

\section{Three stage formulations}\label{sec_defi}

In this section, we formulate the three stage-type statements that form the core of the iterative convex integration process. The first stage formulation (\ref{STA0}), introduced here for the
first time, turns out to be the strongest and implies the other
two. In section \ref{seci3}, we will show how a variant of the third formulation, namely (\ref{STA-VK0}), implies (\ref{STA0}).

\begin{equation}\tag{STA0}\label{STA0}
\left[\quad  \begin{minipage}{12.7cm}
There exist integers $J, S\geq 1$ and constants $\delta_0\in (0,1)$,
$\sigma_0>1$,  $L\in\mathbb{N}$, depending only on
$d,k, J, S$, such that the following holds.
Given $A_0\in\mathcal{C}^L(\bar\omega,\R^{d\times d}_\sym)$
and $\{b_0^{pq}\in \mathcal{C}^L(\bar\omega,\R^{k})\}_{p,q=1\ldots d}$ such that $b_0^{pq}=b_0^{qp}$ for all $p,q=1\ldots d$,
defined on the closure of an open, bounded domain $\omega\subset\R^d$ with smooth boundary, and given
positive parameters $\delta, \lambda, \sigma$,  satisfying:
\begin{equation*}
\begin{split}
& \delta\leq \delta_0,\qquad \lambda\delta^{1/2}\geq 1,
\qquad\sigma\geq \sigma_0\\ 
& \|A_0-\delta H_0\|_0\leq \frac{r_0}{2} \delta,
\qquad \|\nabla^{(m)} A_0\|_0\leq \lambda^m\delta \quad \mbox{ for }\; m=1\ldots L,
\\ & \|\nabla^{(m)} b^{pq}_0\|_0\leq \lambda^{m+1}\delta^{1/2} \quad \mbox{ for }\; m=0\ldots L, \; p,q=1\ldots d,
\end{split}
\end{equation*}
there exist $\tilde v\in\mathcal{C}^2(\bar\omega,\R^{k})$, $\tilde
w\in\mathcal{C}^2(\bar\omega,\R^{d})$ with: 
\begin{equation*}
\begin{split}
& \lambda \|\tilde v\|_0 + \|\nabla \tilde v\|_0\leq C\delta^{1/2},\qquad
\lambda \|\tilde w\|_0 + \|\nabla \tilde w\|_0\leq C\delta, \\
& \|\nabla^{2}\tilde v\|_0\leq C\lambda \sigma^J\delta^{1/2}, \qquad
\|\nabla^{2}\tilde w\|_0\leq C\lambda \sigma^J\delta, \\
&\Big\|A_0 + [\langle b_0^{pq},\tilde v\rangle]_{pq} -
\big(\frac{1}{2}(\nabla \tilde v)^T\nabla \tilde v +\sym\nabla \tilde
w\big)\Big\|_0\leq C \frac{\delta}{\sigma^S}, 
\end{split}
\end{equation*}
with constants $C$ depending only on $d, k, J, S$.
\end{minipage}\quad \right]
\end{equation}

\noindent The above can be regarded as the basic stage statement for the system (\ref{VK}), incorporating a lower-order term $[b_0^{pq}]_{pq}$, and under the assumption that the initial displacements are $v_0=0$, $w_0=0$. The
quantities $H_0\in\mathbb{R}^{d\times d}_\sym$ and $r_0\in(0,1)$
depend only on $d$, as given in Lemma \ref{lem_dec_def}.
We also point out that (\ref{STA0}) is valid in all the recent results (see Examples \ref{ex_first} and \ref{ex_last}), since they proceed by adding corrugations, with amplitudes and correction terms obtained through pointwise decomposition into primitive metrics, and techniques such as the K\"all\'en iteration and integration by parts. These techniques involve only local operations and finitely many derivatives, so the constants $C$ in (\ref{STA0}) are indeed independent of the domain $\omega$. For our subsequent purposes, it would however suffice to assume having uniform $C$ over families of domains related to a fixed  domain by diffeomorphisms whose derivatives, and those of their inverses, are uniformly bounded through the orders required in the construction.

\bigskip

\noindent Next, we formulate the stage statement for the isometric
immersion system (\ref{II}), with a given initial immersion $u$
close to a fixed reference immersion $\underline u$. As before,  $H_0\in\mathbb{R}^{d\times d}_\sym$ and $r_0\in(0,1)$
are as in Lemma \ref{lem_dec_def}. We also set $\tau_0\in (0,1)$ with the property that each enlarged domain:
\begin{equation}\label{tau0}
\omega_\tau \doteq \omega+B_\tau(0) \quad\mbox{ for all }\, 0<\tau\leq 2\tau_0,
\end{equation}
is contractible with smooth boundary.
We assume that $\underline u$ is posed on $\bar\omega_{2\tau_0}$ and define $\bar\rho\in(0,1)$:
\begin{equation}\label{rhobar}
\bar\rho \doteq \min\Big\{\rho(\underline\gamma, 2\Lambda, d, k), \frac{1}{2}\mbox{dist}(\underline u(\bar\omega_{2\tau_0}), \partial\mathcal{U})\Big\},
\end{equation}
that  depends on $\underline \gamma, 2\Lambda, d, k$ through the constant $\rho$ as specified in
Lemma \ref{lem_propa}.

\begin{equation}\tag{STA-II}\label{STA-II}
\left[\quad  \begin{minipage}{12.66cm}
Let $\omega\subset\R^d$ be open, bounded, contractible, with smooth boundary, and let $\mathcal{U}\subset\R^{d+k}$ be open and bounded. 
Given  $\bar g\in\mathcal{C}^{\bar r, \bar\beta}(\bar{\mathcal{U}}, \R^{(d+k)\times (d+k)}_{\sym, >})$ 
for some $\bar r +\bar\beta\in (0, 2]$, 
and $\underline u\in\mathcal{C}^\infty(\bar\omega_{2\tau_0}, \mathcal{U})$, and given constants $\Lambda, \underline\gamma >1$ such that:
\begin{equation}\label{imma}
\begin{split}
& \frac{1}{\Lambda}\Id_{d+k} \leq \bar g \leq \Lambda \Id_{d+k}\quad\mbox{ in } \,\bar{\mathcal{U}},
\\ & \frac{1}{\underline\gamma} \Id_d\leq (\nabla \underline
u)^T\nabla\underline u \leq \underline\gamma \Id_d\quad\mbox{ in } \,\bar\omega_{2\tau_0},
\end{split}
\end{equation}
there exist integers $J, S\geq 1$, an exponent $p>0$, and
constants $\underline\delta\in (0,1)$, $\underline\sigma>1$,
depending only on $\omega, \underline \gamma, \Lambda, d,k,J, S,p$ and on the lower bound of $\mathrm{dist}(\underline u(\bar\omega_{2\tau_0}), \partial\mathcal{U})$,
such that the following holds. Fix any $\tau\in (0,\tau_0)$. Given $g\in \mathcal{C}^{r,\beta}(\bar\omega_{2\tau}, \R^{d\times d}_\sym)$ for some $r+\beta\in (0, 2]$ and $u\in
\mathcal{C}^2(\bar\omega_{2\tau},\mathcal{U})$, and given
any parameters $\delta, \sigma$ and $\lambda\geq 1/\tau$,  satisfying the bounds (with norms on $\bar\omega_{2\tau})$:
\begin{equation*}
\begin{split}
& \delta\leq \underline\delta,\qquad \lambda^{\bar r +\bar\beta}\delta\geq 1, \qquad\sigma\geq \underline\sigma,\qquad \sigma^{p}\delta^{1/2}\leq 1,\\ 
& \|u-\underline u\|_1\leq \frac{\bar\rho}{2}, \qquad \|u\|_2\leq
\lambda\delta^{1/2},\\ 
& \|g-\delta H_0 - (\nabla u)^T(\bar g\circ u)\nabla u\|_0\leq \frac{r_0}{3}\delta,
\end{split}
\end{equation*}
there exists $\tilde u\in\mathcal{C}^2(\bar\omega_\tau,\mathcal{U})$ with the following bounds on norms taken on $\bar\omega_\tau$:
\begin{equation*}
\begin{split}
& \|\tilde u - u\|_1 \leq C_1\delta^{1/2},\qquad
\|\tilde u\|_2\leq C_2\lambda \sigma^J\delta^{1/2}, \\
& \Big\|g-\frac{\delta}{\sigma^S} H_0 - (\nabla \tilde u)^T(\bar g\circ \tilde u)\nabla
\tilde u\Big\|_0\leq \frac{r_0}{5} \frac{\delta}{\sigma^S} + \frac{\|g\|_{r,\beta}}{\lambda^{r+\beta}} + \frac{\|\bar g\|_{\bar r,\bar\beta}}{\lambda^{\bar r+\bar\beta}}.
\end{split}
\end{equation*}
Above, the constant $C_1$ depends only on $\omega, \underline \gamma, \Lambda, d,k,J, S,p$, and $C_2$ on $\omega, \mathcal{U}, \underline u, \underline \gamma, \Lambda, d,k,J, S, p, \|\bar g\|_{\bar r, \bar \beta}$ and $\bar r+\bar\beta$, but not on $\delta,\lambda,\sigma, u$ and $\tau$.
\end{minipage}\quad \right]
\end{equation}
We remark that the third term $\|\bar g\|_{\bar r,\bar\beta}/\lambda^{\bar r+\bar\beta}$ in the right hand side of the final metric error estimate in (\ref{STA-II}) can be actually adjusted to only involve the nonconstant portion of $\bar g$:
\begin{equation}\label{centered}
\Big\|g-\frac{\delta}{\sigma^S} H_0 - (\nabla \tilde u)^T(\bar g\circ \tilde u)\nabla
\tilde u\Big\|_0\leq \frac{r_0}{5} \frac{\delta}{\sigma^S} + \frac{\|g\|_{r,\beta}}{\lambda^{r+\beta}} + \inf_{\bar G \equiv \mathrm{ const}}\frac{\|\bar g - \bar G\|_{\bar r,\bar\beta}}{\lambda^{\bar r+\bar\beta}},
\end{equation}
In particular, when $\bar g = \Id_{d+k}$, then the last term above is null, and (\ref{STA-II}) becomes the formulation of the stage from \cite{LL_last, lew_full2d2k}. The proof of (\ref{centered}) will be given in Remark \ref{rem_Gbar}.

\smallskip

\noindent  We now present the third stage formulation, which concerns
the system (\ref{VK}) with general, nonzero initial displacements
$v, w$. In contrast to (\ref{STA0}) where both
these displacements were zero, this formulation allows the
convex integration procedure to directly begin from a prescribed
pair $v,w$, subject to the assumptions stated below.

\begin{equation}\tag{STA-VK}\label{STA-VK}
\left[\quad  \begin{minipage}{12.36cm}
Let $\omega\subset\R^d$ be an open, bounded domain with smooth boundary.
There exists integers $J, S\geq 1$ and constants $\underline \delta\in (0,1)$,
$\underline{\sigma}>1$, depending only on
$\omega, d,k, J, S$, such that the following holds. Fix any $\tau\in (0,\tau_0)$.
Given $A\in\mathcal{C}^{r,\beta}(\bar\omega_{2\tau},\R^{d\times
  d}_\sym)$ for some $r+\beta\in(0,2]$, and
$v\in\mathcal{C}^2(\bar\omega_{2\tau},\R^k)$, $w\in \mathcal{C}^2(\bar\omega_{2\tau},\R^d)$,
and given any parameters $\sigma$ and $\mathcal{M}\geq 1/\tau$,  satisfying the bounds (with norms taken on $\bar\omega_{2\tau})$:
\begin{equation*}
\begin{split}
& \mathcal{M}\geq \max\{\|v\|_2, \|w\|_2, 1\}, \qquad \sigma\geq \underline \sigma,\\ 
& 0< \|\mathcal{D}\|_0\leq \underline\delta \quad \mbox{ where }\;
\mathcal{D} = A -\big(\frac{1}{2}(\nabla v)^T\nabla v +\sym\nabla w\big),
\end{split}
\end{equation*}
there exist $\tilde v\in\mathcal{C}^2(\bar\omega_\tau,\R^{k})$,
$\tilde w\in\mathcal{C}^2(\bar\omega_\tau,\R^{d})$ such that (with norms on $\bar{\omega}_\tau$):
\begin{equation*}
\begin{split}
& \|\tilde v - v\|_1 \leq C\|\mathcal{D}\|_0^{1/2},\qquad
\|\tilde w - w\|_1 \leq C \|\mathcal{D}\|_0^{1/2} (1+\|\nabla v\|_0),\\
& \|\nabla^{2}\tilde v\|_0\leq C\mathcal{M} \sigma^J, \qquad
\|\nabla^{2}\tilde w\|_0\leq C\mathcal{M}\sigma^J(1+\|\nabla v\|_0),\\
&\Big\|A - \big(\frac{1}{2}(\nabla \tilde v)^T\nabla \tilde v +\sym\nabla \tilde w\big)\|_0
\leq C \Big(\frac{\|\mathcal{D}\|_0}{\sigma^S} + \frac{\|A\|_{r,\beta}}{\mathcal{M}^{r+\beta}}
\|\mathcal{D}\|_0^{(r+\beta)/2}\Big).
\end{split}
\end{equation*}
The constants $C$ above depend on $\omega, d, k, J, S$, but not on $\tau$.
\end{minipage}\quad \right]
\end{equation}

\medskip

\noindent In the following Parts 1-3 we establish the
implications between the three listed formulations  and explain how
they, in turn, yield the flexibility results for the systems (\ref{II})
and (\ref{VK}). 

\smallskip

\subsection{Notation} 

\noindent By $\mathbb{R}^{d\times d}_{\sym}$ we denote the space of symmetric
$d\times d$ matrices, and $\mathbb{R}^{d\times d}_{\sym, >}$ stands
for the set of such matrices which are positive definite. The space of
H\"older continuous vector fields 
$\mathcal{C}^{m,\alpha}(\bar\omega,\R^k)$ consists of restrictions of
all $f\in \mathcal{C}^{m,\alpha}(\mathbb{R}^d,\R^k)$ to the closure of
an open, bounded set $\omega\subset\R^d$. The
$\mathcal{C}^m(\bar\omega,\R^k)$ norm of this restriction is
denoted by $\|f\|_m$, while its H\"older norm in $\mathcal{C}^{m,
  \alpha}(\bar\omega,\R^k)$ is $\|f\|_{m,\alpha}$. By $\nabla^{(m)}f$
we denote the table-valued field of all partial derivatives of $f$ of order
$m$. If $d=k$, the symmetric part of $\nabla f$ is: $\sym\nabla f = (\nabla f + (\nabla f)^T)/2$.
By $C$ we denote a universal constant that may change from line to line.

\section{Preliminary results}\label{sec_preli}

For the proofs of the implications between the stage constructions
in section \ref{sec_defi}, we will need standard mollification and
decomposition estimates, together with results concerning the
propagation of normal frames. We collect these preliminary results in
the present section. 

\medskip

\noindent The first lemma recalls the classical convolution bounds and the 
commutator estimate \ref{stima4} from \cite{CDS}. We will deduce the triple commutator estimate \ref{stima5} from \ref{stima4}, and provide the proof of the interpolation estimate \ref{stima3.5} in  the Appendix.

\begin{lemma}\label{lem_stima} 
Let $\phi\in\mathcal{C}_c^\infty(\R^d,\mathbb{R})$ be a standard
mollifier that is nonnegative, radially symmetric, supported on the
unit ball $B(0,1)\subset\R^d$ and such that $\int_{\mathbb{R}^d} \phi \dx = 1$. Denote: 
$$\phi_l (x) \doteq \frac{1}{l^d}\phi(\frac{x}{l})\quad\mbox{ for all
}\; l\in (0,1], \;  x\in\R^d.$$
Then, for every $f,g,h\in\mathcal{C}^0(\mathbb{R}^d,\R)$, and $u\in\mathcal{C}^0(\R^d, \R^s)$, $\bar f\in\mathcal{C}^{0}(\R^s, \R)$, every
$m\geq 0$ and $r\in\{0,1\}$, $\beta\in (0,1]$, there holds:
\begin{align*}
& \|\nabla^{(m)}(f\ast\phi_l)\|_{0} \leq
\frac{C}{l^m}\|f\|_0,\tag*{(\theequation)$_1$}\vspace{1mm} \refstepcounter{equation} \label{stima1}\\
& \|f - f\ast\phi_l\|_0\leq l^{r
+\beta} \|f\|_{r,\beta},\tag*{(\theequation)$_2$} \vspace{1mm} \label{stima2}\\
& \|\nabla^{(m)}(f\ast\phi_l)\|_{0} \leq
\frac{C}{l^{(m-(r+\beta))_+}}\|f\|_{r,\beta},\tag*{(\theequation)$_3$}\vspace{1mm}  \label{stima3}\\
& \|\bar f\circ u - (\bar f\circ u)\ast\phi_l\|_0\leq \|\nabla\bar f\|_{0,\beta} \Big(l^{1+\beta} \|\nabla u\|_0^{1+\beta} + l^2\|\nabla^2u\|_0\Big),\tag*{(\theequation)$_4$} \vspace{1mm} \label{stima3.5}\\
& \|\nabla^{(m)}\big((fg)\ast\phi_l - (f\ast\phi_l)
(g\ast\phi_l)\big)\|_0\leq {C}{l^{2- m}}\|\nabla f\|_{0}
\|\nabla g\|_{0}, \tag*{(\theequation)$_5$} \label{stima4}\\
& \|\nabla^{(m)}\big((fgh)\ast\phi_l - (f\ast\phi_l)
(g\ast\phi_l)(h\ast\phi_l)\big)\|_0\\ & \qquad \leq {C}{l^{2- m}}\Big(\|f\|_0\|\nabla g\|_{0} \|\nabla h\|_{0}
+\|g\|_0\|\nabla f\|_{0} \|\nabla h\|_{0}+ \|h\|_0\|\nabla f\|_{0} \|\nabla g\|_{0} \Big), \tag*{(\theequation)$_6$} \label{stima5}
\end{align*}
with constants $C>0$ depending only on $d$ and $m$.
\end{lemma}

\medskip

\noindent The next elementary lemma provides a decomposition of symmetric
matrices into linear combinations of rank-one {matrices}, and
defines the quantities $H_0$, $r_0$, depending only on $d$:

\begin{lemma}\label{lem_dec_def}
Let $\{\eta_i\in\R^d\}_{i=1}^{d_*}$ be the $d_*= d(d+1)/2$ unit
vectors in the set $ \big\{\frac{e_p+e_q}{|e_p+e_q|}\big\}_{1\leq p\leq q\leq d}$. 
There exist linear maps $\big\{\bar a_{i}:\R^{d\times d}_\sym\to
\R\big\}_{i=1}^{d_*}$, and a constant $r_0\in(0,1)$, depending only on $d$, such that
$\R^{d\times d}_\sym = \mbox{span}\{\eta_i\otimes \eta_i\}_{i=1}^{d_*}$, with:
$$H=\sum_{i=1}^{d_*}\bar a_{i}(H) \eta_{i}\otimes \eta_{i}
\quad\mbox{ for all } H\in \R^{d\times d}_\sym,$$
and such that, denoting: $H_0 \doteq \sum_{i=1}^{d_*}\eta_{i}\otimes
\eta_{i}$, there holds:
$$|\bar a_{i} (H)-1|\leq\frac{1}{2} \quad \mbox{ for all } H\in
B(H_0,r_0)\subset \R^{d\times d}_{\sym}, \; i=1\ldots d_*.$$
\end{lemma}

\medskip

\noindent  The following result is a direct consequence of the analysis in \cite[Chapter 4, \S2, Corollary 2.5]{Hir_book} 
showing that every vector bundle over a contractible paracompact space is trivial. Since the closure of a bounded, contractible, open $\omega\subset\R^d$ with smooth boundary, is contractible \cite[Chapter 6, proof of Theorem 4.3]{Hir_book}, we obtain:

\smallskip

\begin{lemma}\label{lem_normals}
Let $u\in\mathcal{C}^{n+1}(\bar\omega, \R^{d+k})$ be an immersion defined on the closure of an open, bounded, contractible domain $\omega\subset\R^d$ with smooth boundary. Then, there exists a normal frame $E\in \mathcal{C}^n(\bar\omega, \R^{(d+k)\times k})$ to $u$, namely: 
\begin{equation}\label{norm_frame}
(\nabla u)^T E = 0,\quad E^TE=\Id_k \quad\mbox{ in } \;\bar\omega.
\end{equation}
If $u$ is smooth, then $E$ may be chosen smooth as well.
\end{lemma}

\smallskip

\noindent We recall \cite[Corollary 2.7]{LL_last} that under the additional assumption that $\omega$ is diffeomorphic to $B_1$,
the normal frame $E$ may be chosen to satisfy the bounds:
\begin{equation}\label{stim_nor}
\|\nabla^{(m)}E\|_0\leq C \|\nabla u\|_m\quad\mbox{
  for all } \; m=1\ldots n, 
\end{equation}
where $C$ depends on $\omega, d, k, n$, the modulus of continuity of $\nabla u$ on $\bar\omega$, and the parameter $\gamma$ in: 
\begin{equation*}
\frac{1}{\gamma}\Id_d\leq (\nabla u)^T\nabla u \leq \gamma\Id_d
\quad\mbox{ in } \bar\omega.
\end{equation*}
This result remains valid when $\nabla u$ is replaced by any given  frame $F\in\mathcal{C}^n(\bar\omega, \R^{(d+k)\times d})$ obeying $\frac{1}{\gamma}\Id_d\leq F^TF\leq \gamma\Id_d$.  
The  normal frame $E$ such that $F^TE=0$ and $E^TE=\Id_k$, can be then constructed with  $\|\nabla^{(m)}E\|_0\leq C\|F\|_m$, where $C$ depends only on $\omega, \gamma, d, k, n$, and the modulus of continuity of $F$.
We further remark that in codimension $k\geq 2$, the constant $C$ in (\ref{stim_nor}) cannot, in general, be chosen independently of
the immersion $u$. Specifically, in \cite[Example 2.8]{LL_last} we construct a sequence of immersions $u_n $ from a ball $\bar B_2\subset\R^2$ into $\R^4$, with uniform parameters $\gamma$, and satisfying $\|\nabla u_n\|_1\leq Cn$,
whereas every normal frame $E_n$ to $u_n$ satisfies $\|\nabla
E_n\|_0\geq Cn^2$, with positive constants $C$ independent of $n$ and
of the choice of frame. 

\bigskip

\noindent  We now state three lemmas on the propagation of normal frames
between the given tangent frames $F$ and $\tilde F$, but relative to the metrics
$\Id_{d+k}$ and $\tilde G$, respectively, without requiring the metrics' closeness. In the proof of Theorem \ref{prop_STA0_to_STAII}, we will apply these results with $F=\nabla\underline u$, $\tilde F=\nabla u_0$,
$\tilde G=\bar g_0\circ u_0$ with $g_0$ being a smoothing of $\bar g$. We formulate the lemmas for general
frames and metrics, independently of this application. The proofs will be given in the Appendix.

\begin{lemma}\label{lem_det}
Let $F\in\mathcal{C}(\bar\omega,\R^{(d+k)\times d})$ and $G\in\mathcal{C}(\bar\omega,\R^{(d+k)\times (d+k)}_{\sym, >})$ be defined on the closure of an
open, bounded $\omega\subset\R^d$, and assume that for some $\gamma,\Lambda\geq 1$ there holds:
\begin{equation}\label{assu_imm}
\frac{1}{\gamma}\Id_d\leq F^TGF\leq \gamma\Id_d, \qquad \frac{1}{\Lambda}\Id_{d+k}\leq G\leq \Lambda\Id_{d+k}
\quad\mbox{ in } \bar\omega,
\end{equation}
Then we have:
$$\big(\frac{d}{\gamma\Lambda}\big)^{1/2}\leq |F|\leq (d\gamma\Lambda)^{1/2} \quad\mbox{ and }\quad 
\frac{1}{\gamma^d}\leq\det(F^TG F)\leq \gamma^d\quad 
\mbox{ in } \, \bar\omega.$$
\end{lemma}

\begin{lemma}\label{lem_Tu}
Given $F\in\mathcal{C}^n(\bar\omega,\R^{(d+k)\times d})$, $G\in\mathcal{C}^n(\bar\omega,\R^{(d+k)\times (d+k)}_{\sym, >})$ defined on the closure of an
open, bounded $\omega\subset\R^d$, that satisfy (\ref{assu_imm}), define
$T=F\big(F^TGF\big)^{-1}\in\mathcal{C}^{n}(\bar\omega,\R^{(d+k)\times d})$. Then:
$$\|T\|_0\leq C \quad \mbox{ and }\quad \|\nabla^{(m)}T\|_0\leq C
\hspace{-1mm} \sum_{\sum sp_s=m} \; \prod_{s=1}^m\big(\|\nabla^{(s)}F\|_0 + \|\nabla^{(s)}G\|_0\big)^{p_s},$$
for all $ m=1\ldots n$, where $C$ depends only on $\gamma, \Lambda, d, k, n$.
\end{lemma}

\begin{lemma}\label{lem_propa}
Let $F\in\mathcal{C}^n(\bar\omega,\R^{(d+k)\times d})$ and $\tilde G\in\mathcal{C}^n(\bar\omega,\R^{(d+k)\times (d+k)}_{\sym, >})$ be two fields defined on the closure of an
open, bounded $\omega\subset\R^d$, that satisfy: 
\begin{equation}\label{assu_imm2}
\frac{1}{\gamma}\Id_d\leq F^TF\leq \gamma\Id_d, \qquad \frac{1}{\Lambda}\Id_{d+k}\leq \tilde G\leq \Lambda\Id_{d+k}
\quad\mbox{ in } \bar\omega.
\end{equation}
Let $E\in\mathcal{C}^{n}(\bar\omega, \R^{(d+k)\times k})$
be a normal frame to $F$, namely:
\begin{equation*}
F^TE=0,\qquad E^TE=\Id_k \quad\mbox{ in } \;\bar\omega.
\end{equation*}
There exists a constant $\rho\in (0,1)$ depending on
$\gamma, \Lambda, d,k$, such that if $\tilde F\in\mathcal{C}^{n}(\bar\omega,
\R^{(d+k)\times d})$ satisfies $\|\tilde F- F\|_0\leq \rho$, then:
\begin{equation}\label{out1}
\frac{1}{4\gamma\Lambda}\Id_d \leq \tilde F^T\tilde G\tilde F\leq 4\gamma\Lambda\Id_d  \quad\mbox{ in } \bar\omega,
\end{equation}
and moreover $\tilde F$ has a normal frame 
$\tilde E\in\mathcal{C}^{n}(\bar\omega, \R^{(d+k)\times k})$, namely:
$$\tilde F^T\tilde G\tilde E =0 ,\qquad \tilde E^T\tilde G\tilde E=\Id_k \quad\mbox{ in } \;\bar\omega,$$ 
defined so that for all $m=1\ldots n$ and $i=1\ldots k$ there holds:
\begin{equation}\label{out2}
\|\nabla^{(m)}\tilde E^i\|_0 \leq C \hspace{-3mm}\sum_{\sum_{s=1}^m sp_s=m} \;
\prod_{s=1}^m \Big(\|\nabla^{(s)} \tilde G\|_0 + 
\|\nabla^{(s)} \tilde F\|_0+\|\nabla^{(s)} F\|_0 + \sum_{j\leq i}\|\nabla^{(s)} E^j\|_0\Big)^{p_s},
\end{equation}
with constants $C$ depending only on $\gamma, \Lambda, d, k, n$.
\end{lemma}

\part{(\ref{STA0}) $\Rightarrow$ (\ref{STA-II}) and flexibility results for (\ref{II}) }

\section{Stage (\ref{STA0}) implies stage (\ref{STA-II})}\label{seci1}

The purpose of this section is to show that the stage (\ref{STA0}) for the Monge-Ampère system (\ref{MA}) directly implies the stage (\ref{STA-II}) for the isometric immersion system (\ref{II}), for any $g, \bar g$. The proof is divided into eight steps, with its key calculation being (\ref{shakl}). There, we observe that updating an immersion $u$ by the normal and tangential displacements $v,w$ obtained from (\ref{STA0}), yields an immersion $\tilde u$ whose metric tensor differs from that of $u$ by three contributions:
\begin{itemize}
\item[(i)] the von K\'arm\'an content of $v,w$, 
\item[(ii)] a term linear in $v$ determined by the second fundamental form of $u$ relative to the ambient metric $\bar g$, 
\item[(iii)] additional error terms collected in $\mathcal{R}$.
\end{itemize}
The aforementioned second fundamental form, identified in step 3, provides then the coefficients $b_0^{pq}$ in (\ref{STA0}). The remaining errors are estimated in step 7 and shown to be of smaller order when the initial defect is sufficiently small.
Along the way, we utilize Lemma \ref{lem_normals} to construct a normal frame for the referential immersion $\underline u$, and Lemma \ref{lem_propa} to propagate this frame to the nearby immersion $u$ while preserving the required derivative bounds.

\begin{theorem}\label{prop_STA0_to_STAII}
Assume that
(\ref{STA0}) holds, with some $J, S\geq 1$. Then (\ref{STA-II}) holds
as well, with the same exponents $J,S$ and $p=S+1$.
\end{theorem}  

\begin{proof}
{\bf 1. (Mollification and definition of new parameters)} Let $\omega, \mathcal{U}, \underline u,\underline \gamma, \Lambda, \tau, g, \bar g, u$ be as in the
statement of (\ref{STA-II}) and let $\delta_0, \sigma_0, J,S,L$ be as
in (\ref{STA0}). Fix positive parameters $\delta, \sigma, \lambda$
satisfying the assumptions of (\ref{STA-II}), where we set:
\begin{equation}\label{redef}
\underline \delta = \min\Big\{2\delta_0,
\frac{\bar\rho^2}{4}\Big\}, \qquad \underline
\sigma = \max\Big\{\sigma_0, 7, \frac{6\tbar C}{r_0}\Big\}
\end{equation}
with respect to a constant $\tbar C$ that will be specified below, depending $\underline \gamma, \Lambda, d,k,J,S$. In particular, recalling (\ref{rhobar}) and since $\delta_0,\sigma_0$ in (\ref{STA0}) depend only on $d,k,J,S$, we note the dependence of $\underline\delta$ and $\underline\sigma$ as stated in the formulation of (\ref{STA-II}).
In the course of the proof, we will also utilize other constants $\bar C, \dbar C>1$, defined throughout as needed and depending on $\omega, \mathcal{U},
\underline u, \underline \gamma, \Lambda$, $d,k, J,S, L, \|\bar g\|_{\bar r, \bar\beta}, \bar r+\bar\beta$.
Define the mollified quantities:
$$g_0 = g\ast\phi_{l}, \qquad \bar g_0 = \bar g\ast\phi_{l},\qquad u_0=u\ast\phi_{l}, \quad\mbox{ where } \; l=\frac{1}{\bar C \lambda}. $$
Above, $\phi_l$ is the standard $d$-dimensional (for $g_0, u_0$) or $d+k$-dimensional (for $\bar g_0$) mollifier as in Lemma \ref{lem_stima}. In particular, taking $\bar C$ sufficiently large, we may request that:
\begin{equation}\label{imm_2l}
\frac{1}{2\Lambda}\Id_{d+k}\leq \bar g_0(y)\leq 2\Lambda\,\Id_{d+k}\quad\mbox{ for  $y\in\mathcal{U}$ with $\mathrm{dist}(y, \partial\mathcal{U})>l$}.
\end{equation}
Also, since $\lambda\geq 1/\tau$ we may guarantee that $2l\leq \tau$, so that all the $\omega$-dependent convolutions use only data on $\bar\omega_{2\tau}$.
By \ref{stima1}, \ref{stima2} and recalling that $\|u\|_2\leq\lambda\delta^{1/2}$, we get:
\begin{equation}\label{uno}
\begin{split}  
& \|g_0-g\|_0\leq \frac{\|g\|_{r,\beta}}{\lambda^{r+\beta}},\qquad
\|\bar g_0-\bar g\|_0\leq \frac{\|\bar g\|_{\bar r,\bar\beta}}{\lambda^{\bar r+\bar\beta}}, \qquad 
\|u_0-u\|_1\leq \frac{\|u\|_2}{\bar C\lambda}\leq \delta^{1/2},
\\ & \|\nabla^{(m+2)}u_0\|_0\leq C \bar C^m\lambda^m\|\nabla^{(2)}u\|_0
\leq C\bar C^m\lambda^{m+1}\delta^{1/2}\quad\mbox{for
  all } m=0\ldots L. 
\end{split}
\end{equation}  
with $C$ that depends on $d, k, L$.
By Lemma \ref{lem_normals}, there exists a normal frame to $u$:
$$\underline E\in\mathcal{C}^\infty(\bar\omega_{2\tau_0}, \R^{(d+k)\times k})\quad\mbox{ satisfying: }\quad (\nabla \underline u)^T \underline E= 0,\quad
{\underline E}^T {\underline E}=\Id_{k}\quad\mbox{ in }\; \bar\omega_{2\tau_0}.$$
We fix the frame $\underline E$ once for all, independently of $\tau$ and of the stage input data $u,\delta,\lambda,\sigma$. By (\ref{uno}) and the assumption, we also get:
$$\|\nabla u_0-\nabla\underline u\|_0\leq \|\nabla u_0-\nabla u\|_0+
\|\nabla u -\nabla\underline u\|_0\leq \delta^{1/2} + \frac{\bar\rho}{2}\leq \bar\rho,$$
since $\delta\leq \underline\delta\leq (\bar\rho/2)^2$ in (\ref{redef}). Recall that $\bar\rho$ is not greater than the threshold $\rho$ in Lemma \ref{lem_propa} for $\underline\gamma, 2\Lambda$. We may hence apply Lemma \ref{lem_propa} to the quantities:
$$F=\nabla\underline u,\quad E=\underline E, \quad \tilde F=\nabla u_0, \quad\tilde G = \bar g_0\circ u_0.$$
The first assertion (\ref{out1}) yields that:
\begin{equation}\label{obv}
\frac{1}{8\underline \gamma\Lambda}\Id_d\leq (\nabla u_0)^T(\bar g_0\circ u_0)\nabla u_0\leq 8\underline\gamma\Lambda\,\Id_d \quad \mbox{ in }\,\bar\omega_\tau.
\end{equation}
The second assertion yields existence of the normal frame
$E_0\in\mathcal{C}^\infty(\bar\omega_\tau, \R^{(d+k)\times k})$, namely: $$(\nabla u_0)^T (\bar g_0\circ u_0)E_0= 0, \qquad
E_0^T(\bar g_0\circ u_0)E_0=\Id_{k}\quad \mbox{ in }\;\bar\omega_\tau,$$ 
with bounds that we will develop in the next step.

\smallskip

{\bf 2. (Estimates on the normal frame to $\nabla u_0$)} We claim that,  with $C$ depending
on $\omega, \mathcal{U}, \underline u, \underline\gamma, \Lambda, d, k, L$ and $\|\bar g\|_{\bar r, \bar\beta}$, there holds: \begin{equation}\label{due}
\|\nabla^{(m)}E_0\|_0 \leq C \bar C^m\lambda^m\delta^{1/2} \quad\mbox{for all } \, m=1\ldots L.
\end{equation}
To prove (\ref{due}), we need to estimate all terms in the following consequence of (\ref{out2}):
\begin{equation}\label{out22}
\begin{split}
\|\nabla^{(m)} E_0\|_0 \leq C \hspace{-3mm}\sum_{\sum_{s=1}^m sp_s=m} \;
\prod_{s=1}^m \Big(&\|\nabla^{(s)} (\bar g_0\circ u_0)\|_0 + 
\|\nabla^{(s+1)}\underline u\|_0+\|\nabla^{(s+1)} u_0\|_0  \\ & + \|\nabla^{(s)} \underline E\|_0\Big)^{p_s} \quad\mbox{ for all }\; m=1\ldots L,
\end{split}
\end{equation}
with constant $C$ depending on $\underline\gamma, \Lambda, d, k, L$.
The last three terms in the parentheses above are all bounded by $C\bar C^s\lambda^s\delta^{1/2}$ in view of (\ref{uno}) and since  $\lambda^s\delta^{1/2}\geq \lambda\delta^{1/2}\geq (\lambda^{\bar r +\bar\beta}\delta)^{1/2}\geq 1$. We also recall that $\underline E$ depends on $\underline u$ only, through the application of Lemma \ref{lem_normals}.
We now estimate the term $\|\nabla^{(s)} (\bar g_0\circ u_0)\|_0$. First, \ref{stima3} implies that, with $C$ depending on $\mathcal{U}, d,k,L$:
$$\|\nabla^{(m)}\bar g_0\|_0\leq C \|\bar g\|_{\bar r, \bar\beta} (\bar C\lambda)^{(m-(\bar r+\bar\beta))_+} \quad\mbox{ for } m=1\ldots L+2.$$
We claim that the above may be rewritten as:
\begin{equation}\label{bargg}
\|\nabla^{(m)}\bar g_0\|_0\leq C\|\bar g\|_{\bar r, \bar \beta}(\bar C\lambda)^m\delta^{1/2}  \quad\mbox{ for } m=1\ldots L+2.   
\end{equation}
Indeed, the estimate is obvious for $m\geq 2$. For $m=1$ and $\bar r=0$ we have: $ (\bar C\lambda)^{(m-(\bar r+\bar\beta))_+} = (\bar C\lambda)^{1-(\bar r+\bar\beta)}\leq \bar C\lambda\delta^{1/2}$, while with $\bar r=1$ the same term trivially equals: $1\leq \bar C\lambda\delta^{1/2}$.

\medskip

\noindent Second, from (\ref{uno}) and the assumptions, with $C$ depending on $\omega, \underline\gamma, d,k,L$ (since $\bar \rho<1$):
\begin{equation}\label{use}
\begin{split}
& \|\nabla u_0\|_0\leq \|\nabla \underline u\|_0 + \frac{\bar\rho}{2} + \delta^{1/2}\leq C \\ & \|\nabla^{(m+1)}u_0\|_0\leq C(\bar C\lambda)^m\delta^{1/2} \quad \mbox{ for } m=1\ldots L+1.
\end{split}
\end{equation}
In particular, we directly see that, with $C$ depending on $\omega, \mathcal{U}, \underline\gamma, d,k, L$:
\begin{equation}\label{nowe1}
\|\nabla (\bar g_0\circ u_0)\|_0\leq \|\nabla \bar g_0\|_0\|\nabla u_0\|_0 \leq C \|\bar g\|_{\bar r, \bar\beta} (\bar C\lambda)\delta^{1/2},
\end{equation}
On the other hand, for $m=2\ldots L+2$ we apply the Fa\`a di Bruno formula, again with the final constant $C$ depending on $\omega, \mathcal{U}, \underline\gamma, d,k, L$:
\begin{equation*}
\begin{split}
& \|\nabla^{(m)}(\bar g_0\circ u_0)\|_0 \leq C \sum_{\sum sp_s =m}    
\|\nabla^{(p_1+\ldots p_m)}\bar g_0\|_0\prod_{s=1}^m\|\nabla^{(s)}u_0\|_0^{p_s} \\ & \leq 
C \|\bar g\|_{\bar r, \bar\beta} \sum_{\sum sp_s =m} 
(\bar C\lambda)^{((p_1+\ldots p_m)-(\bar r+\bar\beta))_+} 
\prod_{s=2}^m\big((\bar C\lambda)^{s-1}\delta^{1/2}\big)^{p_s}
\\ & \leq C \|\bar g\|_{\bar r, \bar\beta} \Big((\bar C\lambda)^{m-(\bar r +\bar\beta)} + \sum_{\sum sp_s =m, \; p_1<m} \hspace{-2mm}
(\bar C\lambda)^{((p_1+\ldots p_m)-(\bar r+\bar\beta))_+} 
(\bar C\lambda)^{m-(p_1+\ldots p_m)} \delta^{1/2}\Big),
\end{split}    
\end{equation*}
The first quantity in the right hand side corresponds to the case $p_1=m$, while the second to $p_1<m$. Considering the
cases $p_1+\ldots p_m < \bar r+\bar \beta$ and $p_1+\ldots p_m \geq \bar r+\bar \beta$ separately, we get:
\begin{equation*}
\begin{split}
\|\nabla^{(m)}(\bar g_0\circ u_0)\|_0 & \leq C \|\bar g\|_{\bar r, \bar\beta} \Big((\bar C\lambda)^{m-(\bar r +\bar\beta)} + 
(\bar C\lambda)^{m-1}\delta^{1/2}\Big),
\end{split}    
\end{equation*}
because, recalling that $\bar r+\bar\beta\leq 2$, we see that $p_1+\ldots p_m< \bar r+\bar \beta$ implies $p_1+\ldots p_m=1$ (only possible when $\bar r +\bar\beta>1$).
In particular, the above yields that:
\begin{equation}\label{nowe}
\begin{split}
\|\nabla^{(m)}(\bar g_0\circ u_0)\|_0 & \leq C \|\bar g\|_{\bar r, \bar\beta} (\bar C\lambda)^{m} \Big(\frac{1}{\bar C^{\bar r+\bar\beta}} + \frac{1}{\bar C}\Big) \delta \quad\mbox{ for all }\, m=2\ldots L+2,
\end{split}    
\end{equation}
where the final bound is due to the assumption $\lambda^{\bar r+\bar\beta}\delta\geq 1$. We remark in passing, that the necessity of the above estimate up to $L+2$ rather than $L$ is justified by its application to bound the term $IV$ in step 4 below. Recalling (\ref{out22}), the bounds (\ref{nowe}) and (\ref{nowe1})  complete the proof of (\ref{due}), with constants $C$ as indicated.

\smallskip

{\bf 3. (Identification of $b_0$ as the second fundamental form of $u_0$)} For each $p,q=1\ldots d$, define a smooth field $b_0^{pq}\in\mathcal{C}^\infty(\bar\omega_\tau, \R^k)$ by setting:
\begin{equation}\label{sec_fun} 
\begin{split}
b_0^{pq}= & E_0^T(\bar g_0\circ u_0)\partial_{pq}u_0 + \frac{1}{2}
E_0^T\big((D\bar g_0)_{\mid u_0} \partial_pu_0\big)\partial_q u_0
+ \frac{1}{2}
E_0^T\big((D\bar g_0)_{\mid u_0} \partial_qu_0\big)\partial_p u_0
\\ & - \frac{1}{2} \Big[\Big\langle 
\big((D\bar g_0)_{\mid u_0} E_0e_i\big)\partial_p u_0, \partial_qu_0
\Big\rangle\Big]_{i=1\ldots k}.
\end{split}
\end{equation}
Clearly, $b_0^{pq}=b_0^{qp}$ for all $p,q=1\ldots d$, as required in the assumption of (\ref{STA0}).
Writing in coordinates and using the Einstein summation convention, we rewrite the above as:
\begin{equation*}
\begin{split}
(b_0^{pq})^i & = \big\langle E_0e_i, (\bar g_0\circ u_0)\partial_{pq}u_0\big\rangle \\ & \quad \, + \frac{1}{2} (E_0e_i)^\gamma\Big(
\partial_\alpha(\bar g_0)_{\beta\gamma} + \partial_\beta(\bar g_0)_{\alpha\gamma} - \partial_\gamma(\bar g_0)_{\alpha\beta}\Big)
_{\mid u_0} (\partial_pu_0)^\alpha(\partial_q u_0)^\beta 
\\ & = \Big\langle E_0e_i, (\bar g_0\circ u_0)\Big(\partial_{pq}u_0 + \big[(\Gamma_{\alpha\beta}^j)_{\mid u_0}(\partial_pu_0)^\alpha(\partial_q u_0)^\beta \big]_{j=1\ldots d+k}
\Big)\Big\rangle \\ & = \Big(E_0^T(\bar g_0\circ u_0)\big(\partial_{pq}u_0 + (\Gamma\circ u_0) [\partial_pu_0, \partial_q u_0]\big)\Big)^i,
\end{split}
\end{equation*}
where $\Gamma = \{\Gamma^j_{\alpha\beta}\}_{j,\alpha, \beta=1\ldots d+k}$ denote the Christoffel symbols of the metric $\bar g_0$. The above formula identifies $b_0$ in (\ref{sec_fun}) as the second 
fundamental form of the immersion $u_0$ relative to the ambient metric $\bar g_0$.
We now proceed to estimating the derivatives of $b_0^{pq}$. To this end, we first note the bound, with constant $C$ depending on $\omega, \mathcal{U}, \underline\gamma, d, k, L$:
\begin{equation}\label{barggg}
\|\nabla^{(m)}(D\bar g_0\circ u_0)\|_0\leq C\|\bar g\|_{\bar r, \bar\beta} (\bar C\lambda)^{m+1}\delta^{1/2} \quad\mbox{ for }\; m=0\ldots L.    
\end{equation}
Indeed, for $m=0$ the above follows from (\ref{bargg}), while for $m\geq 1$ we may use Fa\`a di Bruno's formula with the same reasoning as in step 2, in:
\begin{equation*}
\begin{split}
& \|\nabla^{(m)}(D\bar g_0\circ u_0)\|_0 \leq C \sum_{\sum sp_s =m}    
\|\nabla^{(1+p_1+\ldots p_m)}\bar g_0\|_0\prod_{s=1}^m\|\nabla^{(s)}u_0\|_0^{p_s} \\ & \leq 
C \|\bar g\|_{\bar r, \bar\beta} \sum_{\sum sp_s =m} 
(\bar C\lambda)^{((1+p_1+\ldots p_m)-(\bar r+\bar\beta))_+} 
\prod_{s=2}^m\big((\bar C\lambda)^{s-1}\delta^{1/2}\big)^{p_s}
\leq C \|\bar g\|_{\bar r, \bar\beta} (\bar C\lambda)^{m+1}\delta^{1/2}.
\end{split}    
\end{equation*}
Above, terms in which only the first derivatives of $u_0$ occur, must have $p_1=m$ and then $(\bar C\lambda)^{m+1-(\bar r +\bar\beta)}\leq (\bar C\lambda)^{m+1}\delta^{1/2}$. In terms in which some derivatives of order at least $2$ occur, the factor $\delta^{1/2}$ is supplied by the product, and the total frequency power is at most $m+1$.

\medskip

\noindent By (\ref{imm_2l}), (\ref{uno}), (\ref{due}), (\ref{nowe}), (\ref{use}), (\ref{nowe1}) and (\ref{barggg}), we deduce:
\begin{equation}\label{bd_B}
\begin{split}
\|\nabla^{(m)}b_0^{pq}\|_0 & \leq C \hspace{-1mm}\sum_{r+s+t=m} \|\nabla^{(r)} E_0\|_0\|\nabla^{(s)}(\bar g_0\circ u_0)\|_0\|\nabla^{(t+2)}u_0\|_0
\\ & \quad + C \hspace{-3mm}\sum_{r+s+t+v=m} \|\nabla^{(r)} E_0\|_0\|\nabla^{(s)}(D\bar g_0\circ u_0)\|_0\|\nabla^{(t+1)}u_0\|_0\|\nabla^{(v+1)}u_0\|_0
\\ & \leq C \hspace{-1mm}\sum_{r+s+t=m} (\bar C\lambda)^r(\bar C\lambda)^s (\bar C\lambda)^{t+1}\delta^{1/2} 
\\ & \quad + C \hspace{-3mm}\sum_{r+s+t+v=m} (\bar C\lambda)^r(\bar C\lambda)^{s+1}\delta^{1/2} (\bar C\lambda)^{t}(\bar C\lambda)^{v}
\\ & \leq C (\bar C\lambda)^{m+1}\delta^{1/2}\leq \dbar C\lambda^{m+1}\delta^{1/2} \quad\mbox{ for all }\; m=0\ldots L, \quad p,q=1\ldots d,
\end{split}
\end{equation}
with a cumulative constant $\dbar C$ that depends on $\omega, \mathcal{U}, \underline u, \underline\gamma,\Lambda,
d, k, J,S, L, \|\bar g\|_{\bar r, \bar\beta}, \bar r+\bar\beta$.

\smallskip

{\bf 4. (Identification of $A_0$) }
We now define the smooth field $A_0\in\mathcal{C}^\infty(\bar\omega_\tau, \R^{d\times d}_\sym)$ by setting:
$$A_0=\frac{1}{2}\Big(g_0-\frac{\delta}{\sigma^S}H_0 - (\nabla
u_0)^T(\bar g_0\circ u_0)\nabla u_0\Big).$$
In order to validate the $A_0$-related assumption bounds in (\ref{STA0}), we first write:
\begin{equation}\label{decompo}
\begin{split}
& (\nabla u_0)^T(\bar g_0\circ u_0)\nabla u_0 - \big((\nabla
u)^T(\bar g\circ u)\nabla u\big)\ast\phi_{l} = I+ II + III,
\\ & \mbox{where: }\; 
I = (\nabla u_0)^T\big((\bar g\circ u)\ast\phi_l\ast\phi_l\big)\nabla u_0 - \Big((\nabla u)^T \big((\bar g\circ u)\ast\phi_l\big)\nabla u\Big)\ast\phi_l, 
\\ & \qquad\quad \; II = (\nabla u_0)^T\Big(\bar g_0\circ u_0 - \big(\bar g\circ u)\ast\phi_l\ast\phi_l\Big) \nabla u_0,
\\ & \qquad\quad  III \;= \Big((\nabla u)^T\big((\bar g\circ u)\ast\phi_l - \bar g\circ u\big)\nabla u\Big)\ast\phi_l.
\end{split}
\end{equation}
To bound the derivatives of $I$, we use the triple commutator estimate \ref{stima5} and \ref{stima3} in:
\begin{equation}\label{ola1}
\begin{split}
\|\nabla^{(m)} I\|_0 & \leq C (\bar C\lambda)^{m-2} \Big(\|\nabla u\|_0\|\nabla\big((\bar g\circ u)\ast\phi_l\big)\|_0\|\nabla^2u\|_0 + \|(\bar g\circ u)\ast\phi_l\|_0 \|\nabla^2 u\|_0^2\Big) \\ & \leq 
C \|\bar g\|_{\bar r, \bar\beta}(\bar C\lambda)^{m-2} \Big(
(\bar C \lambda)^{(1-(\bar r+\bar\beta))_+}\lambda\delta^{1/2} +\lambda^2\delta\Big) \\ & \leq C \|\bar g\|_{\bar r, \bar\beta}(\bar C\lambda)^{m} \delta \Big(\frac{
(\bar C \lambda)^{(1-(\bar r+\bar\beta))_+}\lambda }{(\bar C\lambda)^2\delta^{1/2}} + \frac{1}{\bar C^2}\Big) \\ & \leq 
\frac{ C \|\bar g\|_{\bar r, \bar\beta}}{\bar C^{\bar r+\bar\beta}}(\bar C\lambda)^{m}  \delta \quad\mbox{ for }\;m=0\ldots L,
\end{split}    
\end{equation}
with $C$ depending on $\omega, \mathcal{U}, \underline\gamma, d,k,L$.
Indeed, both terms in parentheses above are majorized by $1/\bar C^{\bar r +\bar\beta}$. This is clear for the second term and also for the first term in case $\bar r+\bar\beta\geq 1$, because $\lambda \delta^{1/2}\geq 1$ and $\bar C^2 \geq \bar C^{\bar r+\bar\beta}$.
When $\bar r+\bar\beta< 1$, the first term becomes, as claimed:
$$ \frac{(\bar C \lambda)^{(1-(\bar r+\bar\beta))_+}\lambda }{(\bar C\lambda)^2\delta^{1/2}} =
\frac{\lambda }{(\bar C\lambda)^{1+(\bar r+\bar\beta)}\delta^{1/2}} 
\leq  \frac{1}{\bar C^{1+ (\bar r + \bar \beta)} (\lambda^{\bar r+ \bar\beta}\delta)^{1/2}} \leq \frac{1}{\bar C^{1+ (\bar r + \bar \beta)}}. $$
Towards bounding the derivatives of $II$ in (\ref{decompo}), we write:
\begin{equation}\label{split456}
\begin{split}
& \bar g_0\circ u_0 - (\bar g\circ u)\ast\phi_l\ast\phi_l = IV + V +VI
\\ & \mbox{where: } IV = \bar g_0\circ u_0 - (\bar g_0\circ u_0)\ast\phi_l, \quad \\ & \qquad\quad\; V = (\bar g_0\circ u_0 - \bar g\circ u)\ast\phi_l, \\ &
\qquad \quad \; VI = \big(\bar g\circ u - (\bar g\circ u)\ast\phi_l\big)\ast\phi_l. 
\end{split}
\end{equation}
By \ref{stima2} and (\ref{nowe}), we get, with $C$ depending on $\omega, \mathcal{U}, \underline\gamma,  d,k,L$, for all $m=0\ldots L$:
$$\|\nabla^{(m)}IV\|_0\leq \frac{1}{(\bar C\lambda)^2}\|\nabla^{(m+2)}(\bar g_0\circ u_0)\|_0\leq C\|\bar g\|_{\bar r, \bar\beta} (\bar C\lambda)^{m}\Big(\frac{1}{\bar C^{\bar r+\bar\beta}} + \frac{1}{\bar C}\Big)\delta.$$
Similarly, by \ref{stima1}, \ref{stima2} we arrive at:
\begin{equation*}
\begin{split}
\|\nabla^{(m)}V\|_0 & \leq C (\bar C\lambda)^m\|\bar g_0\circ u_0 - \bar g\circ u\|_0 \leq C (\bar C\lambda)^m\Big(\|\bar g_0- \bar g\|_0 + \|\bar g\circ u_0 - \bar g\circ u\|_0\Big)
\\ & \leq C (\bar C\lambda)^m\Big( \frac{\|\bar g\|_{\bar r,\bar\beta}}{(\bar C\lambda)^{\bar r+\bar\beta}} + \|\bar g\|_{\bar r,\bar\beta}\|u_0-u\|_0^{\bar r+\bar\beta} + \|\bar g\|_{\bar r,\bar\beta}\|u_0-u\|_0\Big),
\end{split}
\end{equation*}
where the second term in the right hand side parentheses corresponds to the case $\bar r+\bar\beta<1$ and the third term to the case $\bar r+\bar\beta\geq 1$. Consequently, we use that $\|u\|_2\leq \lambda\delta^{1/2}$ to get that for all $m=0\ldots L$ and with $C$ depending on $\omega, \mathcal{U}, d,k,L$:
\begin{equation*}
\begin{split}
\|\nabla^{(m)}V\|_0 & \leq C \|\bar g\|_{\bar r,\bar\beta} (\bar C\lambda)^m\Big( \frac{1}{(\bar C\lambda)^{\bar r+\bar\beta}} + \big(\frac{\|u\|_2}{(\bar C\lambda)^2}\big)^{\bar r+\bar\beta} + \frac{\|u\|_2}{(\bar C\lambda)^2}\Big)
\\ & \leq C \|\bar g\|_{\bar r,\bar\beta} (\bar C\lambda)^m\Big(
\frac{\delta}{\bar C^{\bar r +\bar\beta}} + \big(\frac{\delta^{1/2}}{\bar C^2\lambda}\big)^{\bar r+\bar\beta} + \frac{\delta^{1/2}}{\bar C^2\lambda}\Big)
\\ & \leq C \|\bar g\|_{\bar r,\bar\beta} (\bar C\lambda)^m\Big(
\frac{\delta}{\bar C^{\bar r +\bar\beta}} + \frac{\delta}{\bar C^{2(\bar r+\bar\beta)}} + \frac{\delta}{\bar C^2}\Big) 
\leq \frac{C\|\bar g\|_{\bar r,\bar\beta} }{\bar C^{\bar r+\bar\beta}}(\bar C\lambda)^m\delta.
\end{split}
\end{equation*}
Finally, we get with $C$ depending on $\omega, \mathcal{U}, \underline\gamma,  d,k,L$:
\begin{equation}
\|\nabla^{(m)}VI\|_0 \leq C (\bar C\lambda)^m
\|\bar g\circ u - (\bar g\circ u)\ast\phi_l\|_0 
\leq \frac{C\|\bar g\|_{\bar r,\bar\beta} }{\bar C^{\bar r+\bar\beta}}(\bar C\lambda)^m\delta,\label{ola6}
\end{equation}
for all $m=0\ldots L$. Indeed, when $\bar r=0$ then by \ref{stima2}:
$$\|(\bar g\circ u)\ast\phi_l - \bar g\circ u\|_0\leq \frac{\|\bar g\circ u\|_{\bar r+\bar\beta}}{(\bar C\lambda)^{\bar r+\bar \beta}} \leq 
\frac{1}{(\bar C\lambda)^{\bar r+\bar \beta}}\big(\|\bar g\|_0 + \|\bar g\|_{\bar r+\bar\beta} \|\nabla u\|_0^{\bar r+\bar\beta}\big) \leq \frac{C\|\bar g\|_{\bar r, \bar\beta}}{(\bar C\lambda)^{\bar r +\bar\beta}},$$
whereas for $\bar r=1$, we use the estimate \ref{stima3.5} and \eqref{uno} to get:
\begin{equation*}
\|(\bar g\circ u)\ast\phi_l - \bar g\circ u\|_0  \leq \|\nabla\bar g\|_{0,\beta} \Big(\frac{\|\nabla u\|_0^{\bar r+\bar\beta}}{(\bar C\lambda)^{\bar r+\bar \beta}} + \frac{\|\nabla^2 u\|_0}{(\bar C\lambda)^{2}}\Big)
\leq C\|\bar g\|_{\bar r, \bar\beta}\Big(\frac{\delta}{\bar C^{\bar r +\bar\beta}} +\frac{\delta^{1/2}}{\bar C^2\lambda}\Big),
\end{equation*}
concluding \eqref{ola6} in both cases.
Consequently, by the formula in (\ref{split456}):
\begin{equation*}
\begin{split}
\|\nabla^{(m)}\big(\bar g_0\circ u_0 - (\bar g\circ u)\ast\phi_l\ast\phi_l\big)\|_0 \leq \frac{C\|\bar g\|_{\bar r,\bar\beta} }{\min\{\bar C^{\bar r+\bar\beta},\bar C\}}(\bar C\lambda)^m\delta \quad\mbox{ for all } m=0\ldots L.
\end{split}
\end{equation*}
We are ready to estimate derivatives of $II$ in (\ref{decompo}), with $C$ depending on $\omega, \mathcal{U}, \underline\gamma,d,k,L$:
\begin{equation}\label{ola2}
\begin{split}
\|\nabla^{(m)} II\|_0 & \leq C\sum_{r+s+t=m}
\|\nabla^{(r+1)}u_0\|_0\|\nabla^{(s)}\big(\bar g_0\circ u_0 - (\bar g\circ u)\ast\phi_l\ast\phi_l\big)\|_0\|\nabla^{(t+1)}u_0\|_0
\\ & \leq C\sum_{r+s+t=m} (\bar C\lambda)^r \frac{C\|\bar g\|_{\bar r,\bar\beta} }{\min\{\bar C^{\bar r+\bar\beta},\bar C\}}(\bar C\lambda)^s\delta  (\bar C\lambda)^t
\\ & \leq \frac{C\|\bar g\|_{\bar r,\bar\beta} }{\min\{\bar C^{\bar r+\bar\beta},\bar C\}}(\bar C\lambda)^m\delta \quad\mbox{ for }\;m=0\ldots L.
\end{split}    
\end{equation}
For the derivatives of $III$ in (\ref{decompo}), arguing as above using \ref{stima1}, \ref{stima2} and (\ref{uno}) yields:
\begin{equation}\label{ola3}
\begin{split}
\|\nabla^{(m)} III\|_0 & \leq C (\bar C\lambda)^m \|(\bar g\circ u)\ast\phi_l - \bar g\circ u\|_0\|\nabla u\|_0^2\\ &  \leq \frac{C\|\bar g\|_{\bar r, \bar\beta}}{\bar C^{\bar r +\bar\beta}}(\bar C\lambda)^m\delta
\quad\mbox{ for }\;m=0\ldots L,
\end{split}    
\end{equation}
with $C$ depending on $\omega, \mathcal{U}, \underline\gamma,d,k,L$.

\smallskip

\noindent At this point, we deduce from (\ref{ola1}), (\ref{ola2}), (\ref{ola3}), through the decomposition (\ref{decompo}), that:
\begin{equation}\label{fabio}
\begin{split}
& \|\nabla^{(m)}\big((\nabla u_0)^T(\bar g_0\circ u_0)\nabla u_0 - ((\nabla u)^T(\bar g\circ u)\nabla u)\ast\phi_{l} 
\big)\|_0 \\ & \leq \frac{C\|\bar g\|_{\bar r,\bar\beta} }{\min\{\bar C^{\bar r+\bar\beta},\bar C\}}(\bar C\lambda)^m\delta \quad\mbox{ for }\;m=0\ldots L,
\end{split}
\end{equation}
with $C$ depending on $\omega, \mathcal{U}, \underline\gamma, d,k,L$. Taking $\bar C$ sufficiently large, in function of $C$, $\|\bar g\|_{\bar r, \bar\beta}$, $\bar r+\bar\beta$, the last assumption in (\ref{STA-II}) and (\ref{fabio}) yield, provided that $\sigma^S\geq\sigma\geq 7$:
\begin{equation}\label{bdsA0}
\begin{split}
& \big\|A_0 - \frac{1}{2}\big(\delta-\frac{\delta}{\sigma^S}\big)H_0\big\|_0 \leq \frac{1}{2}\big\|\big(g-\delta
H_0 - (\nabla u)^T(\bar g\circ u)\nabla u\big)\ast\phi_{l}\big\|_0 
\\ & \qquad\qquad \qquad\qquad \qquad\quad +\frac{1}{2} \big\|(\nabla u_0)^T(\bar g_0\circ u_0)\nabla u_0 - ((\nabla
u)^T(\bar g\circ u)\nabla u)\ast\phi_{l}\big\|_0 \\ & \qquad\qquad ~\, \leq \frac{r_0}{6}\delta + \frac{C\|\bar g\|_{\bar r,\bar\beta} }{\min\{\bar C^{\bar r+\bar\beta},\bar C\}} \delta \leq \frac{r_0}{5}\delta \leq \frac{r_0}{4}\big(\delta-\frac{\delta}{\sigma^S}\big), 
\\ & \|\nabla^{(m)}A_0\|_0\leq \frac{1}{2}\big\|\nabla^{(m)}\big(g-\delta
H_0 - (\nabla u)^T(\bar g\circ u)\nabla u\big)\ast\phi_{l}\big\|_0
\\ & \qquad \qquad \qquad + \frac{1}{2} \big\|\nabla^{(m)}\big((\nabla u_0)^T(\bar g_0\circ u_0)\nabla u_0 - ((\nabla
u)^T(\bar g\circ u)\nabla u)\ast\phi_{l}\big)\big\|_0 \\ &  \qquad\qquad ~\, \leq C\bar C^m\lambda^m\frac{r_0}{3}\delta + \frac{C\|\bar g\|_{\bar r,\bar\beta} }{\min\{\bar C^{\bar r+\bar\beta},\bar C\}} (\bar C\lambda)^m\delta 
\leq \dbar C \lambda^m\delta \quad\mbox{ for }\;m=1\ldots L,
\end{split}
\end{equation}  
with a cumulative constant $\dbar C$ that depends on $\omega, \mathcal{U}, \underline u, \underline\gamma,\Lambda,
d, k, J,S, L, \|\bar g\|_{\bar r, \bar\beta}, \bar r+\bar\beta$.

\smallskip

{\bf 5. (Redefinition of frequency parameters)}
Define:
\begin{equation}\label{redef2}
\delta_{VK}= \frac{1}{2}\big(\delta-\frac{\delta}{\sigma^S}\big), \qquad \lambda_{VK} = 4\dbar C \lambda, \qquad
\sigma_{VK}=\tbar C\sigma,
\end{equation}
where $\tbar C$ is another constant that depends on $\underline\gamma, \Lambda, d, k, J, S$, chosen in the next step.
By (\ref{bdsA0}), (\ref{bd_B}), in combination with (\ref{redef}) and if only $\sigma^S\geq \sigma\geq 2$, we get:
\begin{equation*}
\begin{split}  
& \|A_0 - \delta_{VK}H_0\|_0\leq \frac{r_0}{2}\delta_{VK},
\\ & \|\nabla^{(m)}A_0\|_0\leq \dbar C\lambda^m\delta\leq 4\dbar C\lambda^m\delta_{VK} \leq \lambda_{VK}^m\delta_{VK}\quad \mbox{ for }\; m=1\ldots L,
\\ & \|\nabla^{(m)}b_0^{pq}\|_0\leq \dbar C \lambda^{m+1}\delta^{1/2}
\leq 4\dbar C\lambda^{m+1}\delta_{VK}^{1/2}\leq \lambda_{VK}^{m+1}\delta_{VK}^{1/2}
\quad\mbox{for } \;m=0\ldots L, \; p,q=1\ldots d, \\
& \delta_{VK}\leq\delta_0, \quad \sigma_{VK}\geq \sigma_0,
\\ & \lambda_{VK}\delta_{VK}^{1/2}= \lambda\delta^{1/2} 4\dbar C \big(\frac{1-1/\sigma^S}{2}\big)^{1/2}\geq 4 \big(\frac{1-1/\sigma^S}{2}\big)^{1/2}\geq 1. 
\end{split}
\end{equation*}  
This validates the assumption bounds in (\ref{STA0}).
Therefore, the stage statement (\ref{STA0}) applied on $\bar\omega_{\tau}$ yields existence of $\tilde v\in \mathcal{C}^2(\bar\omega_\tau,\R^k)$, $\tilde
w\in\mathcal{C}^2(\bar\omega_\tau,\R^d)$, with:
\begin{equation}\label{quattro}
\begin{split}
& \lambda_{VK} \|\tilde v\|_0 + \|\nabla \tilde v\|_0\leq C\delta_{VK}^{1/2},\qquad
\lambda_{VK} \|\tilde w\|_0 + \|\nabla \tilde w\|_0\leq C\delta_{VK}, \\
& \|\nabla^{2}\tilde v\|_0\leq C\lambda_{VK} \sigma_{VK}^J\delta_{VK}^{1/2}, \qquad
\|\nabla^{2}\tilde w\|_0\leq C\lambda_{VK} \sigma_{VK}^J\delta_{VK}, \\
&\Big\|A_0 + [\langle b_0^{pq},\tilde v\rangle]_{pq} -
\big(\frac{1}{2}(\nabla \tilde v)^T\nabla \tilde v +\sym\nabla \tilde
w\big)\Big\|_0 \leq C \frac{\delta_{VK}}{\sigma_{VK}^S},
\end{split}
\end{equation}
where the constant $C$ depends on $d,k,J,S$.

\smallskip

{\bf 6. (The updated immersion $\tilde u$ and its metric error)} 
We define the map $\tilde u\in\mathcal{C}^2(\bar\omega_\tau, \mathcal{U})$, with respect to the constructed normal frame $E_0$ and the tangent frame $T_0$:
$$\tilde u = u_0+E_0\tilde v + T_0\tilde w \qquad\mbox{where }\;
T_0=(\nabla u_0)\big((\nabla u_0)^T(\bar g_0\circ u_0)\nabla u_0\big)^{-1}.$$
We will see in step 8 that $\tilde u$ is an immersion. In order to justify the image $\tilde u(\bar \omega_\tau)$ 
being a subset of $\mathcal{U}$, we apply Lemma \ref{lem_Tu} to $F=\nabla u_0$ and $G=\bar g_0\circ u_0$, and recall (\ref{nowe1}) and (\ref{nowe}), to get:
\begin{equation}\label{Tzero}
  \|T_0\|_0\leq C, \qquad \|\nabla^{(m)} T_0\|_0\leq C (1+\|\bar g\|_{\bar r, \bar \beta})^m (\bar C\lambda)^m \delta^{1/2} \quad\mbox{ for all }\; m=1\ldots L,
\end{equation}  
where $C$ depends on $\omega, \mathcal{U},  \underline\gamma, \Lambda, d, k, L$. Consequently, by (\ref{uno}) and (\ref{quattro}):
\begin{equation*}
\begin{split}
\|\tilde u - u\|_0& \leq \|u_0-u\|_0+\|E_0\tilde v\|_0+\|T_0\tilde w\|_0\ \\ & \leq \delta^{1/2} + \frac{C\delta_{VK}^{1/2}}{\lambda_{VK}}+\frac{C\delta_{VK}}{\lambda_{VK}}\leq \delta^{1/2}\Big(1+ \frac{C}{\dbar C}\Big),
\end{split}
\end{equation*}
where $C$ depends on $\underline\gamma, \Lambda, d, k, J, S$. Enlarging $\dbar C$ if needed, so that the right hand side above is bounded by $2\delta^{1/2}$, which implies the claim through (\ref{redef}), since:
\[\|\tilde u-\underline u\|_0\leq \|\tilde u- u\|_0+\|u-\underline u\|_0\leq 2\delta^{1/2}+\frac{\bar\rho}{2} \leq \frac{3}{2}\bar\rho \leq \frac{3}{4}\mathrm{dist}\left(\underline u(\bar \omega_{2\tau_0}),\partial \mathcal U\right).\]
By a direct calculation, we obtain:
\begin{equation}\label{dprep}
\begin{split}
& (\nabla \tilde u)^T(\bar g_0\circ \tilde u)\nabla \tilde u - (\nabla u_0)^T(\bar g_0\circ u_0)\nabla u_0
\\ & = 2\,\sym\big((\nabla u_0)^T(\bar g_0\circ u_0)\nabla (E_0\tilde v)\big) + 2\,\sym\big((\nabla
u_0)^T(\bar g_0\circ u_0)\nabla (T_0\tilde w)\big) \\ & \quad + \nabla (E_0\tilde v)^T (\bar g_0\circ u_0) \nabla (E_0\tilde v)
+ 2\,\sym\big( \nabla (E_0\tilde v)^T (\bar g_0\circ u_0) \nabla (T_0\tilde w)\big) \\ & \quad + 
\nabla (T_0\tilde w)^T (\bar g_0\circ u_0)\nabla (T_0\tilde w)
+ (\nabla \tilde u)^T\big(\bar g_0\circ \tilde u - \bar g_0\circ u_0\big)\nabla\tilde u.
\end{split}
\end{equation}
Since $(\nabla u_0)^T(\bar g_0\circ u_0)E_0 = 0$, the first term in the right hand side above is a matrix field whose each $p,q$ component, $p,q=1\ldots d$, equals:
$$-2 \Big\langle E_0^T(\bar g_0\circ u_0)\partial_{pq}u_0, \tilde v \Big\rangle  - 
\Big\langle E_0^T\big((D\bar g_0)_{\mid u_0} \partial_pu_0\big)\partial_q u_0, \tilde v \Big\rangle
- \Big\langle
E_0^T\big((D\bar g_0)_{\mid u_0} \partial_qu_0\big)\partial_p u_0, \tilde v \Big\rangle.$$
On the other hand, the last term in the right hand side of (\ref{dprep}) contains the following term:
$$\Big\langle \Big[\Big\langle 
\big((D\bar g_0)_{\mid u_0} E_0e_i\big)\partial_p u_0, \partial_qu_0
\Big\rangle\Big]_{i=1\ldots k}, \tilde v \Big\rangle$$
By definition (\ref{sec_fun}), the sum of the two displayed expressions above, equals:
$$- 2  [\langle b_0^{pq},\tilde v\rangle]_{pq}.$$
Further, the second term in the right hand side of (\ref{dprep})
contains the term $2\, \sym\nabla \tilde w$, as $(\nabla u_0)^T(\bar g_0\circ u_0)T_0 =
\Id_d$, whereas the third term contains the quadratic expression $(\nabla \tilde v)^T\nabla\tilde v$, since $E_0^T(\bar g_0\circ u_0)E_0=\Id_k$. Consequently, (\ref{dprep}) becomes:
\begin{equation}\label{shakl}
\begin{split}
& (\nabla \tilde u)^T(\bar g_0\circ \tilde u)\nabla \tilde u - (\nabla u_0)^T(\bar g_0\circ u_0)\nabla u_0
= -2  [\langle b_0^{pq},\tilde v\rangle]_{pq} + \big(\nabla \tilde
v)^T\nabla\tilde v + 2\,\sym\nabla\tilde w + \mathcal{R} \vspace{1mm}
\\ & \mbox{where: } 
\mathcal{R} = 2\,\sym\big[\langle \partial_p u_0, (\bar g_0\circ u_0)(\partial_q T_0)\tilde w\rangle\big]_{pq} + 2\,\sym\big[\langle (\partial_p E_0)\tilde v, (\bar g_0\circ u_0)E_0\partial_q\tilde v\rangle\big]_{pq}
\\ & \qquad\qquad\quad \; + \big[\langle(\partial_pE_0)\tilde v, (\bar g_0\circ u_0)(\partial_q E_0)\tilde v\rangle\big]_{pq} + 2\,\sym \big(\nabla(E_0\tilde v)^T(\bar g_0\circ u_0) \nabla (T_0\tilde w)\big) 
\\ & \qquad\qquad\quad \; + \nabla(T_0\tilde w)^T(\bar g_0\circ u_0) \nabla (T_0\tilde w) \\ & \qquad\qquad\quad \; +
(\nabla u_0)^T \big(\bar g_0\circ\tilde u- \bar g_0\circ u_0- (D\bar g_0)_{\mid u_0}(\tilde u- u_0)\big)\nabla u_0  
\\ & \qquad\qquad\quad \; + (\nabla u_0)^T (D\bar g_0)_{\mid u_0} (T_0\tilde w)\nabla u_0 
\\ & \qquad\qquad\quad \; + 2\,\sym \big((\nabla u_0)^T (\bar g_0\circ \tilde u - \bar g_0\circ u_0)\nabla (E_0\tilde v +T_0\tilde w)\big)
\\ & \qquad\qquad\quad \; + \nabla (E_0\tilde v +T_0\tilde w)^T 
(\bar g_0\circ \tilde u - \bar g_0\circ u_0)\nabla (E_0\tilde v +T_0\tilde w).
\end{split}
\end{equation}
We may now write:
\begin{equation}\label{cinque}
\begin{split}
& g-\frac{\delta}{\sigma^S}H_0-(\nabla \tilde u)^T(\bar g\circ \tilde u)\nabla \tilde u \\ &
= 2A_0+(g-g_0) - \Big((\nabla \tilde u)^T(\bar g_0\circ \tilde u)\nabla \tilde u - (\nabla u_0)^T(\bar g_0\circ u_0)\nabla u_0\Big) \\ & \quad \;\; + (\nabla \tilde u)^T(\bar g_0\circ \tilde u - \bar g\circ \tilde u)\nabla \tilde u
\\ & = 2\Big(A_0 + [\langle b_0^{pq},\tilde v\rangle]_{pq}-
\big(\frac{1}{2}(\nabla \tilde v)^T\nabla \tilde v +\sym\nabla \tilde
w\big)\Big) + (g-g_0) - \mathcal{R} \\ & \quad \;\; + (\nabla \tilde u)^T(\bar g_0\circ \tilde u - \bar g\circ \tilde u)\nabla \tilde u.
\end{split}
\end{equation}
To estimate the four terms in the right hand side above, we recall
(\ref{quattro}) and (\ref{uno}), so that it only remains to bound the
error term $\mathcal{R}$ and the auxiliary term $(\nabla \tilde u)^T(\bar g_0\circ \tilde u - \bar g\circ \tilde u)\nabla \tilde u$. 

\smallskip

{\bf 7. (Estimating the remainder terms)} 
We estimate the first five terms in $\mathcal{R}$, by using (\ref{use}), (\ref{Tzero}), (\ref{quattro}) and (\ref{due}), with $C$ that depends on $\omega, \mathcal{U}, \underline u, \underline\gamma, \Lambda, d,k, J, S, L, \|\bar g\|_{\bar r,\bar\beta}$ and, presently, another constant $\tbar C$ that depends only on $\underline\gamma, \Lambda, d,k,J,S$, namely:
\begin{equation*}
\begin{split}
&\tbar C \Big(C\bar C \lambda\delta^{1/2}\frac{\delta_{VK}}{\lambda_{VK}} +C\bar C \lambda\delta^{1/2}\frac{\delta_{VK}}{\lambda_{VK}} 
+ (C\bar C\lambda \delta^{1/2})^2 \frac{\delta_{VK}}{\lambda_{VK}^2}
\\ & \quad + \big(C\bar C\lambda\delta^{1/2}\frac{\delta_{VK}^{1/2}}{\lambda_{VK}} + \delta_{VK}^{1/2}\Big)
\big(C\bar C\lambda\delta^{1/2}\frac{\delta_{VK}}{\lambda_{VK}} + \delta_{VK}\big) + \big(C\bar C\lambda\delta^{1/2}\frac{\delta_{VK}}{\lambda_{VK}} + \delta_{VK}\big)^2\Big) \leq \tbar C\delta^{3/2},
\end{split}  
\end{equation*}  
valid provided that $\dbar C$ is chosen sufficiently large in the definition of $\lambda_{VK}$.
We now estimate the remaining four terms in $\mathcal{R}$, with the same dependence of constants $C$, $\tbar C$ and always taking $\dbar C$ sufficiently large. The sixth term is bounded by (\ref{use}), \ref{stima3} and (\ref{quattro}):
\begin{equation*}
\begin{split}
& \|(\nabla u_0)^T \big(\bar g_0\circ\tilde u- \bar g_0\circ u_0- (D\bar g_0)_{\mid u_0}(\tilde u- u_0)\big)\nabla u_0 \|_0\leq \frac{1}{2}\|\nabla u_0\|_0^2\|\nabla^2\bar g_0\|_0\|\tilde u- u_0\|^2_0 \\ & \leq
C\|\nabla^2\bar g_0\|_0\big(\|\tilde v\|_0+\|\tilde w\|_0\big)^2
\leq C(\bar C\lambda)^{2-(\bar r+\bar\beta)}\Big(\frac{\delta_{VK}^{1/2}}{\lambda_{VK}} + \frac{\delta_{VK}}{\lambda_{VK}} \Big)^2
\\ & \leq C\frac{(\bar C \lambda)^2}{\lambda^{\bar r+\bar\beta}}\frac{\delta}{\lambda_{VK}^2}
\leq C\bar C^2 \lambda^2\frac{\delta^2}{\lambda_{VK}^2}\leq \delta^2,
\end{split}  
\end{equation*}  
in view of the second order Taylor expansion.
For the seventh term, we observe that $(D\bar g_0)_{\mid u_0} \nabla u_0 = \nabla (\bar g_0\circ u_0)$ and therefore by (\ref{use}), (\ref{nowe1}), (\ref{quattro}), (\ref{obv}) and the definition of the field $T_0$, we get:
\begin{equation*}
\begin{split}
\|(\nabla u_0)^T(D \bar g_0)_{\mid u_0}(T_0\tilde w)\nabla u_0\|_0 & \leq
\|\nabla u_0\|^2\|\nabla(\bar g_0\circ u_0)\|_0\big\|\big((\nabla u_0)^T(\bar g_0\circ u_0)\nabla u_0\big)^{-1}\big\|_0\|\tilde w\|_0 \\ & \leq C
\bar C \lambda\delta^{1/2}\frac{\delta_{VK}}{\lambda_{VK}} 
\leq \delta^{3/2},
\end{split}  
\end{equation*}  
For the eighth and ninth terms, we first use \ref{stima3}, (\ref{quattro}), (\ref{Tzero}) to get:
\begin{equation*}
\begin{split}
\|\bar g_0\circ \tilde u - \bar g_0\circ u_0\|_0 & \leq \|\nabla \bar g_0\|_0 \|E_0\tilde v + T_0\tilde w\|_0\leq C(\bar C\lambda)^{(1-(\bar r+\bar\beta))_+}
\Big(\frac{\delta_{VK}^{1/2}}{\lambda_{VK}} + \frac{\delta_{VK}}{\lambda_{VK}}\Big)
\\ & \leq C\bar C \lambda^{(1-(\bar r+\bar\beta))_+}\frac{\delta^{1/2}}{\lambda_{VK}}
= C\bar C \lambda^{(1-(\bar r+\bar\beta))_+-1}\frac{\delta^{1/2}}{\dbar C}\leq \frac{C\bar C}{\dbar C}\delta.
\end{split}    
\end{equation*}
Indeed, when $\bar r =0$ then $\lambda^{(1-(\bar r+\bar\beta))_+-1}= \lambda^{-(\bar r+\bar\beta)}\leq\delta$, while when $\bar r=1$ then the same term equals $\lambda^{-1}\leq \delta^{1/2}$. We also have:
$$\|\nabla(E_0\tilde v +T_0\tilde w)\|_0\leq C\Big(\delta_{VK}^{1/2} + \bar C\lambda\delta^{1/2}\frac{\delta_{VK}^{1/2}}{\lambda_{VK}} + \delta_{VK} + \bar C\lambda\delta^{1/2}\frac{\delta_{VK}}{\lambda_{VK}}\Big) \leq C \delta^{1/2}.$$
The above two displayed formulas yield the estimate on the eighth and ninth terms in $\|\mathcal{R}\|_0$:
\begin{equation*}
\begin{split}
& \|(\nabla u_0)^T (\bar g_0\circ \tilde u - \bar g_0\circ u_0)\nabla (E_0\tilde v +T_0\tilde w)\big\|_0 \\ & \qquad\qquad + \|\nabla (E_0\tilde v +T_0\tilde w)^T 
(\bar g_0\circ \tilde u - \bar g_0\circ u_0)\nabla (E_0\tilde v +T_0\tilde w)\|_0 \\ & \leq   C\|\bar g_0\circ \tilde u - \bar g_0\circ u_0\|_0\Big(\|\nabla(E_0\tilde v +T_0\tilde w)\|_0 + \|\nabla(E_0\tilde v +T_0\tilde w)\|_0^2\Big) 
\\& \leq \frac{C\bar C}{\dbar C}\delta^{3/2}\leq \delta^{3/2}.
\end{split}    
\end{equation*}
In conclusion, we obtain the bound:
\begin{equation}\label{RSbound}
\|\mathcal{R}\|_0 \leq \tbar C\delta^{3/2} \leq \tbar C \frac{\delta}{\sigma^{S+1}}.   
\end{equation}
We now consider the final term in the right hand side of (\ref{cinque}). By (\ref{use}) and (\ref{uno}):
\begin{equation}\label{addi}
\|(\nabla \tilde u)^T(\bar g_0\circ \tilde u - \bar g\circ \tilde u)\nabla \tilde u\|_0\leq \|\nabla \tilde u\|_0^2\|\bar g_0-\bar g\|_0\leq C\frac{\|\bar g\|_{\bar r,\bar\beta}}{(\bar C\lambda)^{\bar r+\bar \beta}}\leq
\frac{\|\bar g\|_{\bar r,\bar\beta}}{\lambda^{\bar r+\bar \beta}},
\end{equation}
provided that $\bar C$ is large (here, in dependence on $\underline\gamma, \Lambda, d,k,J,S, \bar r+\bar\beta$). 
In conclusion, (\ref{cinque}) yields in view of (\ref{STA0}) and (\ref{uno}), (\ref{RSbound}), (\ref{addi}), with $C$ depending on
$d, k, J, S$:
\begin{equation}\label{addi2}
\big\|g-\frac{\delta}{\sigma^S}H_0-(\nabla \tilde u)^T(\bar g\circ\tilde u)\nabla \tilde u\big\|_0 
\leq C\frac{\delta_{VK}}{\sigma_{VK}^S} + \frac{\|g\|_{r,\beta}}{\lambda^{r+\beta}} + \tbar C \frac{\delta}{\sigma^{S+1}} + \frac{\|\bar g\|_{\bar r,\bar\beta}}{\lambda^{\bar r+\bar \beta}}.
\end{equation}
Recalling (\ref{redef}), (\ref{redef2}) and choosing first $\tbar C$ and then $\underline \sigma$ sufficiently large, we arrive at:
$$ C\frac{\delta_{VK}}{\sigma_{VK}^S} + \tbar C \frac{\delta}{\sigma^{S+1}} \leq \frac{C}{\tbar C}\frac{\delta}{\sigma^S} + \frac{\tbar C}{\underline\sigma} \frac{\delta}{\sigma^{S}} \leq \frac{r_0}{5}\frac{\delta}{\sigma^S},$$
which yields the claimed bound:
\begin{equation*}
\big\|g-\frac{\delta}{\sigma^S}H_0-(\nabla \tilde u)^T(\bar g\circ \tilde u)\nabla \tilde u\big\|_0 
\leq \frac{r_0}{5}\frac{\delta}{\sigma^S} + \frac{\|g\|_{r,\beta}}{\lambda^{r+\beta}} +\frac{\|\bar g\|_{\bar r,\bar\beta}}{\lambda^{\bar r+\bar \beta}}.
\end{equation*}

\smallskip

{\bf 8. (Final bounds on $\tilde u$)} By (\ref{imm_2l}), (\ref{uno}), (\ref{due}), 
(\ref{redef2}), (\ref{quattro}), (\ref{Tzero}), we get:
\begin{equation*}
\begin{split}
\|\tilde u - u\|_1& \leq \|u_0-u\|_1+\|\tilde u - u_0\|_1\leq
\delta^{1/2} + \|E_0\tilde v\|_1+ \|T_0\tilde w\|_1 \\ & \leq
\delta^{1/2} + \|E_0\|_0\|\tilde v\|_1 + \|\nabla E_0\|_0\|\tilde v\|_0 +
\|T_0\|_0\|\tilde w\|_1 + \|\nabla T_0\|_0\|\tilde w\|_0
\\ & \leq \delta^{1/2} + C\delta_{VK}^{1/2} + \dbar C \lambda\delta^{1/2}\frac{\delta_{VK}^{1/2}}{\lambda_{VK}}
+ C\delta_{VK} + \dbar C \lambda\delta^{1/2}\frac{\delta_{VK}}{\lambda_{VK}}\leq C_1\delta^{1/2},
\end{split}
\end{equation*} 
where $C$ depends only on $\underline\gamma, \Lambda, d,k,J,S$. Thus, taking $\dbar C$ in (\ref{redef2}) sufficiently large, we see that 
$C_1$ depends on $\underline\gamma, \Lambda, d, k, J, S$. It remains to estimate the second derivatives:
\begin{equation*}
\begin{split}
\|\nabla^{(2)}\tilde u\|_0& \leq \|\nabla^{(2)} u_0\|_0 + \|\nabla^{(2)}(\tilde u-u_0)\|_0\leq
C\lambda\delta^{1/2} + \|\nabla^{(2)}(E_0\tilde v)\|_0+ \|\nabla^{(2)}(T_0\tilde w)\|_0 \\ & \leq
C\big(\lambda\delta^{1/2} + \lambda_{VK}\sigma_{VK}^J\delta_{VK}^{1/2}
+ \lambda\delta^{1/2}\delta_{VK}^{1/2} + \lambda^2\delta^{1/2}\frac{\delta_{VK}^{1/2}}{\lambda_{VK}}
\\ & \qquad \; + \lambda_{VK}\sigma_{VK}^J\delta_{VK} + \lambda\delta^{1/2}\delta_{VK} 
+ \lambda^2\delta^{1/2}\frac{\delta_{VK}}{\lambda_{VK}} \big) \leq C_2\lambda\sigma^J\delta^{1/2},
\end{split}
\end{equation*}    
again by (\ref{uno}), (\ref{due}), (\ref{quattro}), (\ref{Tzero}), 
with $C_2$ absorbing $\bar C, \dbar C$ and hence depending on the entire list of quantities $\omega, \mathcal{U}, \underline u, \underline\gamma, \Lambda, d, k, J, S, \|\bar g\|_{\bar r, \bar\beta}, \bar r+\bar\beta$.  This ends the proof.
\end{proof}

\medskip

\noindent

\begin{remark}\label{rem_Gbar}
We now indicate how the proof of Theorem \ref{prop_STA0_to_STAII} above implies the refined estimate (\ref{centered}). Fix a constant matrix $\bar G\in\R^{(d+k)\times (d+k)}_{\sym}$ and observe that, writing $\bar g_0 - \bar g=(\bar g-\bar G)\ast\phi_l - (\bar g-\bar G)$, the middle bound in (\ref{uno}) can be replaced by:
$$\|\bar g_0-\bar g\|_0\leq \frac{\|\bar g-\bar G\|_{\bar r, \bar\beta}}{\lambda^{\bar r+\bar\beta}}.$$
The bound on $\|\bar g_0-\bar g\|_0$ is used in exactly two places in our proof: in the estimate of the derivatives of term $V$ in step 4, and in the formula (\ref{addi}). 
Hence, in order to validate (\ref{centered}), the construction of $\tilde u$ through the fields $\bar g_0, E_0, b_0^{pq}, A_0, T_0$ and the choice of all parameters in it, remain unchanged and independent of $\bar G$. In the bounds on $\|\nabla^{(m)}V\|_0$ we use (\ref{uno}) as stated, whereas the right hand side of (\ref{addi}) gets replaced by ${\|\bar g-\bar G\|_{\bar r, \bar\beta}}/ \lambda^{\bar r+\bar\beta}$. Consequently, in the final defect estimate (\ref{addi2}), only the last term is updated, and we get its altered form:
$$\big\|g-\frac{\delta}{\sigma^S}H_0-(\nabla \tilde u)^T(\bar g\circ\tilde u)\nabla \tilde u\big\|_0 
\leq C\frac{\delta_{VK}}{\sigma_{VK}^S} + \frac{\|g\|_{r,\beta}}{\lambda^{r+\beta}} + \tbar C \frac{\delta}{\sigma^{S+1}} + \frac{\|\bar g-\bar G\|_{\bar r,\bar\beta}}{\lambda^{\bar r+\bar \beta}}.$$
Adjusting $\bar C, \dbar C, \tbar C$ as before, and minimizing over all constant $\bar G$, there follows (\ref{centered}). Obviously, in the discussed formula, the infimum  may be replaced by taking the minimum.
\end{remark}

\section{Nash-Kuiper iteration scheme based on (\ref{STA-II})}\label{sec_NK_II}

The result we prove in this section is referred to as the {\em Nash-Kuiper iteration scheme} for system (\ref{II}). It provides an inductive application of the stage (\ref{STA-II}), starting from a  prepared referential immersion $\underline u$ whose metric defect is sufficiently small in terms of its immersability constant $\underline \gamma$. Such an immersion is later obtained through the {\em pre-stage} in Proposition \ref{lem_prestage}, followed by the refinement in Corollary \ref{lem_initial_immersion}. Successive applications of (\ref{STA-II}) then produce a sequence of short $\mathcal{C}^2$ immersions converging in $\mathcal{C}^{1,\alpha}$, with $\alpha$ in the range (\ref{range}).
The proof follows the outline of \cite[Theorem 1.3]{lew_full2d2k}, with all constants specified through explicit conditions.

\begin{proposition}\label{thm_NK_II}
Let $\omega\subset\R^d$ be open, bounded, contractible, with smooth boundary, and let $\mathcal{U}\subset\R^{d+k}$ be open and bounded. 
Fix $\tau_0\in (0,1)$ satisfying  (\ref{tau0}). Let 
$g\in\mathcal{C}^{r+\beta}(\bar\omega_{2\tau_0}, \R^{d\times d}_{\sym, >})$ and $\bar g\in\mathcal{C}^{\bar r+\bar \beta}(\bar{\mathcal{U}}, \R^{(d+k)\times (d+k)}_{\sym, >})$ with some $r+\beta, \bar r+\bar\beta \in (0,2]$, and let $\Lambda\geq 1$ with:
$$\frac{1}{\Lambda}\Id_{d+k}\leq \bar g\leq \Lambda\,\Id_{d+k} \quad\mbox{ in }\;\bar{\mathcal{U}}.$$
Fix $\underline\gamma\geq 1$, the integers $J,S\geq 1$, an exponent $p>0$, a threshold parameter $d_\mathcal{U}>0$, and a regularity exponent $\alpha$ within the range:
\begin{equation}\label{range}
0<\alpha<\min\Big\{\frac{1}{1+ 2J/S}, \frac{r+\beta}{2}, \frac{\bar r+\bar \beta}{2}\Big\},
\end{equation}
Then, there exists $\underline\delta\in (0,1)$, depending on $\omega, \underline\gamma, \Lambda,d,k,J,S,p,\alpha, r+\beta, \bar r+\bar\beta, d_\mathcal{U}$ such that the following holds. Given $\underline u \in\mathcal{C}^\infty(\bar\omega_{2\tau_0}, \mathcal{U})$ and $\delta_0\in (0,\underline\delta)$ which satisfy the following conditions:
\begin{itemize}
\item[(i)] $\displaystyle \frac{1}{\underline\gamma}\Id_d\leq (\nabla\underline u)^T\nabla \underline u\leq \underline\gamma\,\Id_d$ in $\bar\omega_{2\tau_0}$, \vspace{1mm}
\item[(ii)] $\mathrm{dist}(\underline u(\bar\omega_{2\tau_0}),\partial\mathcal{U})\geq d_\mathcal{U}$, \vspace{1mm}
\item[(iii)] (\ref{STA-II}) holds with metrics $g, \bar g$, exponents $J,S,p$, the base domain $\omega$, and $\underline u$ as the reference immersion, \vspace{1mm} 
\item[(iv)] $\displaystyle \big\|g-\delta_0H_0 - (\nabla\underline u)^T(\bar g\circ \underline u)\nabla \underline u\big\|_0\leq\frac{r_0}{3}\delta_0$ with norm taken on $\bar\omega_{2\tau_0}$,
\end{itemize}
there exists an immersion $\bar u\in\mathcal{C}^{1,\alpha}(\bar\omega, \mathcal{U})$ such that:
\begin{equation}\label{hasab}
\|\bar u - \underline u\|_1\leq 2C_1\delta_0^{1/2}\quad\mbox{ and }\quad
(\nabla \bar u)^T(\bar g\circ\bar u)\nabla \bar u = g \quad\mbox{in }\; \bar\omega,
\end{equation}
where $C_1$ is the constant in (\ref{STA-II}), depending only on $\omega, \underline\gamma, \Lambda, d,k,J,S,p$.
\end{proposition}
\begin{proof}
{\bf 1. (Setting the iteration parameters).} Let $g,\bar g, u$ be  as in the statement, together with $\underline\gamma, \Lambda, J,S,p,d_\mathcal{U}$ and the  exponent $\alpha$ in (\ref{range}). 
All norms of $g, \underline u$ are assumed to be taken on $\bar\omega_{2\tau_0}$.
Fix any $\theta$ in the following range (for example, we may choose the midpoint):
\begin{equation}\label{th_t}
\alpha<\theta<\min\Big\{\frac{1}{1+2J/S}, \frac{r+\beta}{2}, \frac{\bar r+\bar \beta}{2} \Big\}.
\end{equation}
In particular, we note that:
\begin{equation}\label{ko1}
\frac{J}{S}<\frac{1}{2\theta} - \frac{1}{2}.
\end{equation}
Further, choose $b>1$ sufficiently close to $1$, in order to ensure that:
\begin{equation}\label{ko1.5}
\begin{split}
& \frac{p(b-1)}{S}<\frac{1}{2}, \qquad b<\min\Big\{\frac{r+\beta}{2\theta}, \frac{\bar r + \bar\beta}{2\theta}\Big\},
\\ & \frac{1-\alpha}{2}-\alpha b \Big(\frac{1}{2\theta}-\frac{1}{2}\Big)>0.
\end{split}
\end{equation}
The last above condition is viable, since at $b=1$ its left hand side equals $\frac{1}{2}(1-{\alpha}/{\theta})$, which is positive by (\ref{th_t}). Consequently, the same inequality holds also for $b$ sufficiently close to $1$.

\smallskip

\noindent Recall the stage formulation (\ref{STA-II}) and denote by $\underline\delta_{II}>0$ and $\underline\sigma_{II}\geq 1$ the threshold parameters in there, depending on $\omega, \underline\gamma, \Lambda, d,k,J,S,p, d_\mathcal{U}$, with the property that any $u, \delta, \lambda, \sigma$ satisfying the indicated conditions yield existence of an immersion $\tilde u$ with indicated bounds, involving constants $C_1, C_2\geq 1$. We set, similarly to (\ref{rhobar}):
$$\rho_* = \min\big\{\rho(\underline\gamma, 2\Lambda, d,k), \frac{d_\mathcal{U}}{2}\big\},$$
and define the current threshold $\underline\delta\in (0,1)$, via conditions:
\begin{equation}\label{ko1.6}
\underline\delta\leq \underline\delta_{II}, \qquad \underline\delta^{-(b-1)/S}\geq \underline\sigma_{II}, \qquad \underline \delta^{-(\ln b)/2}\geq {2},\qquad 
2C_1\underline\delta^{1/2}\leq\frac{\rho_*}{2}.
\end{equation}
It is clear that $\underline\delta$ depends on $\omega, \underline\gamma, \Lambda,d,k,J,S,p,\alpha, r+\beta, \bar r+\bar\beta, d_\mathcal{U}$, because of $b$
depending on $J,S,p,\alpha, r+\beta, \bar r+\bar\beta$. 
Let now $\underline u\in\mathcal{C}^\infty(\bar\omega_{2\tau_0},\mathcal{U})$ and $\delta_0\in (0,\underline \delta)$ have properties (i)-(iv). Define:
\begin{equation}\label{ko1.7}
\begin{split}
a=\frac{1}{\delta_0} \quad & \mbox{ and } \quad \lambda_0\geq \max\Big\{a^{\frac{1}{2\theta}}, \frac{2}{\tau_0}\Big\} \\ & \mbox{ such that: }\quad 
\|\underline u\|_2\leq \frac{\lambda_0}{a^{1/2}}, \qquad 
\frac{r_0}{5} + a^b\Big(\frac{\|g\|_{r,\beta}}{\lambda_0^{r+\beta}} 
+ \frac{\|\bar g\|_{\bar r,\bar\beta}}{\lambda_0^{\bar r+\bar \beta}}\Big)\leq \frac{r_0}{3}.
\end{split}
\end{equation}
We finally set  for all $n\geq 0$, the main iteration parameters:
\begin{equation}\label{ko2}
\delta_n =\frac{1}{a^{b^n}},\quad \lambda_n = \lambda_0C_2^n a^{\frac{b^n-1}{2
    \theta}},\quad \sigma_{n+1} =
\Big(\frac{\delta_n}{\delta_{n+1}}\Big)^{1/S} \quad\mbox{so that }\,
\delta_{n+1} = \frac{\delta_n}{\sigma_{n+1}^S}.
\end{equation}
Since $\lambda_n\geq 2^n\lambda_0\geq 2^{n+1}/\tau_0$, the assumption on the relation of $\lambda$ and $\tau$ listed in (\ref{STA-II}) holds.

\smallskip

{\bf 2. (Validating the parameters relations for induction)} We will obtain $\bar u$ as the limit of 
$\{u_n\in\mathcal{C}^{2}(\bar\omega_{\tau_0/2^n},\mathcal{U})\}_{n=0}^\infty$, converging in
$\mathcal{C}^{1,\alpha}(\bar\omega)$. These immersions will be constructed iteratively,
starting from $u_0 = \underline u\in \mathcal{C}^\infty(\bar\omega_{\tau_0},\mathcal{U})$ and successively applying (\ref{STA-II}) with:
$$u=u_{n},\qquad \delta =\delta_n, \qquad\lambda=\lambda_n, \qquad\sigma=\sigma_{n+1},\qquad \tau=\frac{\tau_0}{2^{n+1}},$$
to get $\tilde u = u_{n+1}$. In this step, we validate that for all $n\geq 0$ there holds:
\begin{align*}
& \sigma_{n+1}\geq\underline\sigma_{II},\qquad \sigma_{n+1}^p\delta_n^{1/2}\leq 1,
\qquad \delta_{n+1}\leq\underline\delta_{II},\qquad \lambda_{n+1}^{\bar r +\bar\beta}\delta_{n+1}\geq 1,
\tag*{(\theequation)$_1$}\refstepcounter{equation} \label{ko61} \\
& \frac{r_0}{5}\frac{\delta_n}{\sigma_{n+1}^S} 
+ \frac{\|g\|_{r,\beta}}{\lambda_n^{r+\beta}} 
+ \frac{\|\bar g\|_{\bar r,\bar \beta}}{\lambda_n^{\bar r+\bar \beta}}
\leq \frac{r_0}{3}\delta_{n+1},
\tag*{(\theequation)$_2$} \label{ko62}\\
&  C_2\lambda_n\delta_n^{1/2}\sigma_{n+1}^J\leq \delta_{n+1}^{1/2}\lambda_{n+1},
\tag*{(\theequation)$_3$} \label{ko63}
\end{align*}
In order to show the four bounds in \ref{ko61}, we use the definition
(\ref{ko2}). The first of the claimed bounds follows in virtue of the second condition in (\ref{ko1.6}):
$$\sigma_{n+1}=a^{(b^{n+1}-b^n)/S} = a^{b^n(b-1)/S}\geq \underline \delta^{-b^n(b-1)/S}\geq \underline\delta^{-(b-1)/S}\geq \underline\sigma_{II}.$$
For the second bound, we use the first condition in (\ref{ko1.5}):
$$\sigma_{n+1}^p\delta_n^{1/2}= a^{b^n((b-1)p/S - 1/2)}\leq 1.$$
The third bound is clear since $\delta_{n+1}\leq \delta_0\leq
\underline\delta\leq\underline\delta_{II}$, recalling the first condition in (\ref{ko1.6}). The fourth and final bound in \ref{ko61} is due to having $C_2\geq 1$ and $\lambda_0\geq a^{\frac{1}{2\theta}}$:
$$ \lambda_{n+1}^{\bar r +\bar\beta}\delta_{n+1} = \big(\lambda_0 C_2^{n+1}\big)^{\bar r +\bar\beta}
a^{b^{n+1}(\frac{\bar r +\bar\beta}{2\theta} -1) - \frac{\bar r+\bar\beta}{2\theta}}\geq \Big(\frac{\lambda_0}{a^{\frac{1}{2\theta}}}\Big)^{\bar r +\bar\beta} a^{b(\frac{\bar r +\bar\beta}{2\theta} -1)} 
\geq \Big(\frac{\lambda_0}{a^{\frac{1}{2\theta}}}\Big)^{\bar r +\bar\beta}\geq 1$$
This ends the proof of \ref{ko61}.
To show \ref{ko62}, we observe that, by definition:
$$\frac{r_0}{5}\frac{\delta_n}{\sigma_{n+1}^S} =\frac{r_0}{5}\delta_{n+1},$$ 
and thus it suffices to justify the following calculation and apply the last condition in (\ref{ko1.7}):
\begin{equation*} 
\delta_{n+1}\lambda_n^{r+\beta} = \big(\lambda_0 C_2^n\big)^{r+\beta} a^{
\frac{b^n-1}{2\theta}(r+\beta)-b^{n+1} } \geq \lambda_0^{r+\beta} a^{b^n(\frac{r+\beta}{2\theta}-b)
  -\frac{r+\beta}{2\theta}} \geq\frac{\lambda_0^{r+\beta}}{a^b}.
\end{equation*}
Indeed, the above is a consequence of having $C_2\geq 1$ and the second condition in (\ref{ko1.5}). The parallel bound $\delta_{n+1}\lambda_n^{\bar r+\bar\beta} \geq \lambda_0^{\bar r+\bar \beta}/{a^b}$ is obtained exactly the same, leading to:
$$\Big(\frac{\|g\|_{r,\beta}}{\lambda_n^{r+\beta}} + \frac{\|\bar g\|_{\bar r,\bar\beta}}{\lambda_n^{\bar r+\bar\beta}} \Big)\frac{1}{\delta_{n+1}}\leq 
\Big(\frac{\|g\|_{r,\beta}}{\lambda_0^{r+\beta}} + \frac{\|\bar g\|_{\bar r,\bar\beta}}{\lambda_0^{\bar r+\bar\beta}} \Big) a^b\leq \frac{r_0}{3}-\frac{r_0}{5},$$
which is precisely \ref{ko62}. Finally, \ref{ko63} follows by (\ref{ko1}) in:
\begin{equation*} 
\frac{\lambda_{n+1}}{\lambda_n}\big(\frac{\delta_{n+1}}{\delta_n}\big)^{1/2}\frac{1}{\sigma_{n+1}^J}
= C_2 a^{b^n\frac{b-1}{2\theta} - b^n(b-1)/2-b^n(b-1)J/S} =
C_2 a^{b^n(b-1)(\frac{1}{2\theta}-\frac{1}{2}-\frac{J}{S}) }\geq C_2.
\end{equation*}
This ends the validation of the parameter progression bounds \ref{ko61}-\ref{ko63}. 

\smallskip

{\bf 3. (The remaining bound for the iteration)} We thus see that as long as:
\begin{equation}\label{ko3}
\|u_n - \underline u\|_1\leq \frac{\bar \rho}{2} \quad \mbox{ on }\, \bar\omega_{\tau_0/2^n},
\end{equation}
with $\bar \rho$ given in (\ref{rhobar}) in function of $\underline u$, one can indeed define $u_{n+1}$ by an  application of (\ref{STA-II}). 
Noting the bound $b^i\geq 1+(\ln b)i$ valid for all $i\geq 0$ and $b\geq 1$, we then get: 
\begin{equation}\label{ko7} 
\begin{split}
\|u_{n+1}-\underline u\|_1 & \leq C_1\sum_{i=0}^{n}\delta_i^{1/2}= C_1\sum_{i=0}^{n} \frac{1}{a^{b^i/2}}\leq C_1\sum_{i=0}^{n}
\Big(\frac{1}{a}\Big)^{1/2+(\ln b)i/2} \\ & \leq \frac{C_1}{a^{1/2}}\sum_{i=0}^\infty
\Big(\frac{1}{a}\Big)^{(\ln b)i/2} = \frac{C_1}{a^{1/2}} \frac{a^{(\ln b)/2}}{a^{(\ln b)/2}-1} \\ & \leq \frac{2C_1}{a^{1/2}} = 2C_1\delta_0^{1/2}
\leq \frac{\rho_*}{2}\leq \frac{\bar\rho}{2} \quad \mbox{ on }\, \bar\omega_{\tau_0/2^{n+1}},
\end{split}
\end{equation}
if only $a^{(ln b)/2}\geq 2$, which holds in view of the third condition in (\ref{ko1.6}) since $1/a=\delta_0\leq \underline \delta$. In the final two estimates in (\ref{ko7}) we also used the last condition in (\ref{ko1.6}) and $\rho_*\leq \bar\rho$. Consequently, the sequence
$\{u_n\in\mathcal{C}^2(\bar\omega_{\tau_0/2^n},\mathcal{U})\}_{n=0}^\infty$ is well
defined and it satisfies, for all $n\geq 0$:
\begin{equation}\label{ko8}
\begin{split}
& \|u_{n+1}-u_n\|_1\leq C_1\delta_{n}^{1/2} \quad \mbox{ on }\, \bar\omega_{\tau_0/2^{n+1}},
\\ & \|u_{n}\|_2\leq \delta_{n}^{1/2}\lambda_{n},\qquad  \big\|g-\delta_n H_0-(\nabla u_n)^T(\bar g\circ u_n)\nabla u_n\big\|_0\leq \frac{r_0}{3}\delta_{n}
\quad \mbox{ on }\, \bar\omega_{\tau_0/2^n}.
\end{split}
\end{equation}

\smallskip

{\bf 4. (Convergence in $\mathcal{C}^{1,\alpha}$)} We are now ready to conclude the proof, arguing that
$\{u_n\}_{n=0}^\infty$ is Cauchy in 
$\mathcal{C}^{1,\alpha}(\bar\omega)$. To this end, we use the
interpolation inequality in H\"older spaces:
$$\|u_{n+1}-u_n\|_{1,\alpha}\leq C \|u_{n+1}-u_n\|_2^\alpha\|u_{n+1}-u_n\|_1^{1-\alpha}.$$
Since $\{\delta_n^{1/2}\lambda_n\}_{n=0}^\infty$ is an increasing sequence,
the above and (\ref{ko8}) imply for all $n\geq 0$:
\begin{equation}\label{ko9}
\begin{split}
& \|u_{n+1}-u_n\|_{1,\alpha} \leq C(\delta_{n+1}^{1/2}
\lambda_{n+1})^{\alpha}\delta_n^{(1-\alpha)/2} \\ & = C \big(\lambda_0C_2^{n+1}\big)^\alpha \big(a^{
  \frac{b^{n+1}-1}{2\theta}- b^{n+1}/2}\big)^\alpha a^{-b^n(1-\alpha)/2}=
C \big(\lambda_0 C_2^{n+1}\big)^\alpha a^{-\frac{\alpha}{2\theta}} a^{-b^n q}, 
\end{split}
\end{equation}
where $q$ is the following constant, independent of $n$ and positive by the last condition in (\ref{ko1.5}):
$$q= \frac{1-\alpha}{2} - \alpha b\big(\frac{1}{2\theta}-\frac{1}{2}\big).$$
We now argue that the series $\sum_{n=1}^\infty\|u_{n+1}-u_n\|_{1,\alpha}$ is Cauchy, as its increments in (\ref{ko9}) are bounded by the quantities whose series converges. For this last property, we use the ratio test:
$$ \frac{\big(\lambda_0 C_2^{n+2}\big)^\alpha  a^{-b^{n+1} q}}{\big(\lambda_0 C_2^{n+1}\big)^\alpha  a^{-b^n q}} = C_2^\alpha\Big(\frac{1}{a}\Big)^{b^n(b-1)q} \to 0 \quad \mbox{ as }\, n\to\infty.$$
In conclusion, $\bar
u\in\mathcal{C}^{1,\alpha}(\bar\omega,\mathcal{U})$ may be defined as the limit: 
$$u_n\to \bar u \qquad\mbox{as }\; n\to\infty \quad \mbox{in }\;
\mathcal{C}^{1,\alpha}(\bar\omega,\R^{d+k}).$$
Indeed, property (\ref{ko3}) whose validity for all $n\geq 0$ we proved above, implies that $\|\bar u - \underline u\|_0\leq \bar \rho/2$, so that the definition of $\bar\rho$ in (\ref{rhobar}) ensures $\bar u(\bar\omega)\subset\mathcal U$. Similarly, passing to the limit in (\ref{ko7}), yields the first claimed condition in (\ref{hasab}):
$$\|\bar u - \underline u\|_1\leq 2C_1\delta_0^{1/2}.$$
For the second condition, we use the last statement in (\ref{ko8}) to get:
$$\big\|g-(\nabla u_n)^T(\bar g\circ u_n)\nabla u_n\big\|_0 \leq  \Big(\frac{r_0}{3}+|H_0|\Big)\delta_n
\to 0 \quad\mbox{as }\; n\to\infty.$$
Since the left hand side above clearly converges to $\|g-(\nabla \bar u)^T(\bar g\circ \bar u)\nabla \bar u\|_0$ as $n\to\infty$, there follows the second claim in (\ref{hasab}). This ends the proof.
\end{proof}

\section{Pre-stage for (\ref{II}) and constructing the family of referential immersions} \label{sec_NK_II2}

In this section we establish one of the building blocks for the convex integration algorithm, namely the {\em pre-stage}. This construction reduces the defect of a given strictly short immersion $u$ to an arbitrarily small size, while keeping the modified immersion $\underline u$ in a prescribed $\mathcal{C}^0$ neighbourhood of $u$ and retaining a uniform bound on its immersability constant, albeit possibly increasing its second derivatives.
For $\bar g=\Id_{d+k}$ (confer Example \ref{counter_ex}), the content of Proposition \ref{lem_prestage} below follows from the proof of the celebrated Nash-Kuiper theorem (see \cite{Nash2, Kuiper} and also the note in \cite[\S 1.3.18]{GR}). Here, consistently with the approach of the present paper, we deduce this result for arbitrary $g, \bar g$, from the basic h-principle for the system (\ref{MA}), relative to zero reference displacements, proved in \cite[Theorem 6.2]{lew_conv}:\\
\begin{equation}\tag{STA-INI0}\label{STA-INI0}
\left[\quad  \begin{minipage}{12.19cm}
Let $\omega\subset\R^d$ be an open, bounded domain. 
Given $A\in\mathcal{C}(\bar\omega,\R^{d\times d}_{\sym, >})$
and given parameters $\epsilon, \delta>0$, there exist $v\in\mathcal{C}^1(\bar\omega,\R^{k})$, $
w\in\mathcal{C}^1(\bar\omega,\R^{d})$ with: 
\begin{equation*}
\begin{split}
& \|v\|_0 + \|w\|_0\leq \epsilon, \\
& \|\nabla v\|_0\leq C\|A\|_0^{1/2}, \qquad \|\nabla w\|_0\leq  C\|A\|_0,
\\ &\Big\|A - \big(\frac{1}{2}(\nabla v)^T\nabla v +\sym\nabla w\big)\Big\|_0\leq \delta,
\end{split}
\end{equation*}
where constants $C$ depend only on $\omega, d, k$.
\end{minipage}\quad \right]
\end{equation}

\medskip

\noindent We note in passing that, by approximation, one can have both displacement fields $v$ and $w$ above smooth, which is what we use in the proofs below. Two points are now in order. First, the immersion $\underline u$ constructed below is {\em almost} the referential immersion from which the successive applications of the stage (\ref{STA-II}) proceed in the Nash-Kuiper iteration of Proposition \ref{thm_NK_II}. The actual referential immersion will be obtained in Corollary \ref{lem_initial_immersion}, by additionally centering the defect at a multiple of $H_0$ while preserving the uniformity properties in (\ref{coro_asser}).

\medskip

\noindent The second point concerns the method of proof. As in Theorem \ref{prop_STA0_to_STAII}, we seek to construct $\underline u$ using $v$ and $w$, obtained in (\ref{STA-INI0}), as normal and tangential corrections. The resulting change in the metric tensor consists of the von Kármán content, a contribution involving the appropriate curvature form, and additional error terms collected in $\mathcal{R}$. We would like to treat the latter two contributions as lower-order errors. However, (\ref{shakl2}) includes in $\mathcal{R}$ the term:
\begin{equation}\label{bad_term}
\nabla(Tw)^T(\bar g\circ u)\nabla(Tw),
\end{equation}
which is quadratic in $\nabla w$ and hence not necessarily small, because the original defect need not be small. This caused no difficulty in the proof of Theorem \ref{prop_STA0_to_STAII}: there, the original defect was of order $\delta\ll 1$ and (\ref{bad_term})  was hence of order $\delta^2$.

\medskip

\noindent The obstruction arises when using the Kuiper-corrugation version of (\ref{STA-INI0}) above, but it is not present in the Nash-spiral construction used in \cite[Theorem 4.1]{lew_full2d2k}, which however necessitates codimension $k\geq2$. In the present setting we wish to include $k=1$; the remedy is to divide the defect into $N$ equal portions and correct them successively. Each portion is then of order $\delta/N$, so its corresponding quadratic error is of order $(\delta/N)^2$, and summing over the $N$ corrections gives $\delta^2/N$, which can be made arbitrarily small by choosing $N$ sufficiently large.

\smallskip
\noindent Another approach has been used in  \cite{DIS1/5} and also in subsequent papers \cite{CHI2, In2026}, where $k=1$ and $\bar g= \mathrm{Id}_{d+1}$. More precisely, \cite[Section 8.1]{DIS1/5} applied the classical Nash decomposition into finitely many primitive metrics, approximately cancelled by successive Conti--De Lellis--Sz\'ekelyhidi corrugations. These corrugations have profiles $\Gamma^t,\Gamma^n$ that satisfy the identity $2\partial_\xi\Gamma^t +(\partial_\xi \Gamma^t)^2+(\partial_\xi\Gamma^n)^2=s^2$, see \cite[Proposition 3.5]{DIS1/5}, allowing for incorporating the leading quadratic tangential contribution (\ref{bad_term}) into the primitive metrics, leaving an error of the desired lower order at each step. Consequently, after finitely many perturbations, the defect can be made arbitrarily small by choosing the frequencies large enough.

\medskip

\begin{proposition}\label{lem_prestage}
Satisfaction of (\ref{STA-INI0}) on an open, bounded, contractible $\omega\subset\R^d$ with smooth boundary implies the following. 
Let $g\in \mathcal{C}(\bar\omega, \R^{d\times d}_{\sym, >})$ and $\bar g\in\mathcal{C}(\bar{\mathcal{U}},\R^{(d+k)\times(d+k)}_{\sym, >})$ be two Riemannian metrics, posed on $\bar\omega$ and $\bar{\mathcal{U}}$, respectively, where $\mathcal{U}\subset\R^{d+k}$ is an open and bounded set. Let $u\in\mathcal{C}^1(\bar\omega, \mathcal{U})$ be an immersion, satisfying:
$$(\nabla u)^T(\bar g\circ u)\nabla u<g \quad\mbox{ in } \; \bar\omega.$$ 
Then, for every $\epsilon, \delta>0$ there exists $\underbar
u\in\mathcal{C}^\infty(\bar\omega,\mathcal{U})$, such that:
\begin{equation*}
\begin{split}
& \frac{1}{\underline \gamma} \Id_d\leq (\nabla \underbar u)^T\nabla
\underbar u \leq \underline\gamma\Id_d \quad\mbox{ and } \quad
(\nabla \underline u)^T(\bar g\circ\underline u)\nabla \underline u<g \quad\mbox{ in } \,\bar\omega,\\
& \|u-\underbar u\|_0\leq \epsilon \quad\mbox{ and }\quad
\|g - (\nabla \underline u)^T(\bar g\circ \underline u)\nabla \underline u\|_0\leq \delta,
\end{split}
\end{equation*}
where  $\underline \gamma \geq 1$
depends only on (the ellipticity and immersability parameters of) $g, \bar g, u$. 
\end{proposition}
\begin{proof}
{\bf 1.} Let $\omega, \mathcal{U}, g, \bar g, u$ and $\epsilon, \delta$ be as in the statement, and denote by $\gamma, \Lambda\geq 1$ the related ellipticity and immersability constants, guaranteeing that: 
\begin{equation}\label{imm3}
\begin{split}
& \frac{1}{\gamma}\Id_d\leq (\nabla u)^T\nabla u\leq \gamma\Id_d \quad 
\mbox{and}\quad \frac{1}{\gamma}\Id_d\leq g \leq \gamma\Id_d \qquad\mbox{ in } \; \bar\omega,
\\ & \frac{1}{\Lambda}\Id_{d+k}\leq \bar g \leq \Lambda\Id_{d+k} \qquad\mbox{ in } \;\bar{\mathcal{U}}.
\end{split}
\end{equation}
By density of $\mathcal{C}^\infty$ in $\mathcal{C}^1$, we find a smooth field $u_0\in\mathcal{C}^\infty(\bar\omega,\mathcal{U})$ such that:
\begin{equation}\label{imm4}
\begin{split}
&\|u_0-u\|_0\leq \frac{\epsilon}{2}, \\
& \mathcal{D}_0\doteq g -(\nabla u_0)^T(\bar g\circ u_0)\nabla u_0 > 0 \quad \mbox{and}\quad \frac{1}{2\gamma}\Id_d\leq (\nabla u_0)^T\nabla u_0\leq 2\gamma\Id_d \quad\mbox{ in } \; \bar\omega.
\end{split}
\end{equation}
In particular, the last assertion above implies that:
\begin{equation}\label{imm4.5}
\frac{1}{2\gamma\Lambda}\Id_d\leq (\nabla u_{0})^T(\bar g\circ u_{0})\nabla u_0 \leq 2\gamma\Lambda\,\Id_d \quad\mbox{ in }\,\bar\omega.
\end{equation}
Fix parameters $\theta\in (0,1)$ and $N\in\mathbb{N}$, so that:
$$\theta\|\mathcal{D}_0\|_0\leq \frac{\delta}{3}, \qquad N \geq C\frac{\|\mathcal{D}_0\|_0^2}{\delta},\qquad C\|\mathcal{D}_0\|_0^2\Id_d\leq N\theta\mathcal{D}_0 \quad\mbox{in }\;\omega.$$
where $C$ is a constant that depends only on $\omega, d,k,\gamma, \Lambda$ that we specify below. We will inductively construct immersions 
$\{u_i\in\mathcal{C}^\infty(\bar\omega, \mathcal{U})\}_{i=1}^N$ via pairs of displacements $\{(v_i, w_i)\in\mathcal{C}^\infty(\bar\omega, \R^k\times \R^d)\}_{i=1}^N$ obtained by an application of (\ref{STA-INI0}) with the following parameters:
$$A=\frac{1-\theta}{2N}\mathcal{D}_0\in\mathcal{C}(\bar\omega, \R^{d\times d}_{\sym, >}), \qquad \epsilon_i>0,\qquad \delta_i>0,$$
where $\epsilon_i, \delta_i$ are chosen sufficiently small, in function of all the quantities obtained thus far.
Then, with a constant $C$ depending only on $\omega, d, k$, there holds:
\begin{equation}\label{imm5}
\begin{split}
&\|v_i\|_0 + \|w_i\|_0\leq \epsilon_i \\
& \|\nabla v_i\|_0\leq C\Big(\frac{\|\mathcal{D}_0\|_0}{N}\Big)^{1/2}, \quad\|\nabla w_i\|_0\leq C \frac{\|\mathcal{D}_0\|_0}{N},
\\ & \Big\| A - \big(\frac{1}{2}(\nabla v_i)^T\nabla v_i + \sym\nabla w_i\big)\Big\|_0\leq \delta_i.
\end{split}
\end{equation}
Along the induction process, we will also prove a relaxed version of the bound (\ref{imm4.5}), namely:
\begin{equation}\label{imm6}
\frac{1}{3\gamma\Lambda}\Id_d\leq (\nabla u_{i})^T(\bar g\circ u_{i})\nabla u_i \leq 3\Lambda\gamma\,\Id_d \quad\mbox{ in }\,\bar\omega \quad\mbox{ for all }\; i=0\ldots N.
\end{equation}

\smallskip

{\bf 2.} Assume that an immersion $u_i\in\mathcal{C}^\infty(\bar\omega, \mathcal{U})$ satisfying (\ref{imm6}) has been constructed, for some $i=0\ldots N-1$. By Lemma \ref{lem_normals}, we get a smooth normal frame $E_{i,0}\in\mathcal{C}^\infty(\bar\omega, \R^{(d+k)\times k})$, with:
$$(\nabla u_i)^TE_{i,0}=0, \qquad E_{i,0}^TE_{i,0}=\Id_k \qquad \mbox{ in } \;\bar\omega.$$
We now fix a smooth $G_i\in\mathcal{C}^\infty(\bar \omega, \R^{(d+k)\times (d+k)}_{\sym, >})$ such that:
\begin{equation}\label{imm7}
\|G_i-\bar g\circ u_i\|_0\leq \eta_{i+1},
\end{equation}
with $\eta_{i+1}$ small enough. In particular, we request that, in view of (\ref{imm6}):
\begin{equation}\label{imm7.4}
\frac{1}{4\gamma\Lambda}\Id_d\leq (\nabla u_{i})^TG_i\nabla u_i \leq 4\Lambda\gamma\,\Id_d \quad\mbox{ in }\,\bar\omega, 
\end{equation}
allowing to define the tangent frame $T_i=\nabla u_i\big((\nabla u_i)^TG_i\nabla u_i\big)^{-1}\in\mathcal{C}^\infty(\bar\omega, \R^{(d+k)\times d})$.
We observe in passing that the same condition as in (\ref{imm7.4}) also holds for $\big((\nabla u_i)^TG_i\nabla u_i\big)^{-1}$, which yields $\|\big((\nabla u_i)^TG_i\nabla u_i\big)^{-1}\|_0\leq 4\Lambda\gamma$, while (\ref{imm6}) and (\ref{imm3}) imply $\|\nabla u_i\|_0\leq 2\Lambda\gamma^{1/2}$, and hence:
\begin{equation}\label{imm7.45}
\|T_i\|_0\leq \|\nabla u_i\|_0\|\big((\nabla u_i)^TG_i\nabla u_i\big)^{-1}\|_0
\leq 8\Lambda^2\gamma^{3/2}.
\end{equation}
Further, we apply Lemma \ref{lem_propa} to $F=\tilde F=\nabla u_i$, $E=E_{i,0}$ and $\tilde G= G_i$, and obtain a normal frame $E_{i}\in\mathcal{C}^\infty(\bar\omega, \R^{(d+k)\times k})$ such that the following identities hold on $\bar\omega$:
\begin{equation}\label{imm7.5}
(\nabla u_i)^TG_iE_i=0,\qquad E_i^TG_iE_i=\Id_k,\qquad E_i^TG_iT_i=0, \qquad (\nabla u_i)^TG_iT_i=\Id_d.
\end{equation}
We now set the ensuing immersion $u_{i+1}$, whose image is contained in $\mathcal{U}$ provided $\epsilon_{i+1}$ is small:
$$u_{i+1}= u_i+E_iv_{i+1} + T_iw_{i+1}\in\mathcal{C}^\infty(\bar\omega, \mathcal{U}).$$
We record the following bounds, arising in view of (\ref{imm5}) and (\ref{imm7.45}), and also in view of the second condition in (\ref{imm7.5}) which yields $\|E_i\|_0^2\leq 2\Lambda$, provided that $\eta_{i+1}$ in (\ref{imm7}) has been chosen small enough to assure that $(2\Lambda)^{-1}\Id_{d+k}\leq G_i\leq 2\Lambda\Id_{d+k}$:
\begin{equation}\label{imm7.6}
\begin{split}
& \|u_{i+1}-u_i\|_0\leq C\epsilon_{i+1},
\\ & \|\nabla u_{i+1}\|_0\leq \|\nabla u_i\|_0 + \epsilon_{i+1}(\|\nabla E_i\|_0+\|\nabla T_i\|_0) + \|E_i\|_0\|\nabla v_{i+1}\|_0+ \|T_i\|_0\|\nabla w_{i+1}\|_0
\\ & \qquad \qquad \,
\leq \epsilon_{i+1}(\|\nabla E_i\|_0+\|\nabla T_i\|_0) + 
C(1+\|\mathcal{D}_0\|_0),
\end{split}
\end{equation}
Above, the constant $C$ depends only on $\omega, d, k, \gamma,\Lambda$.
A direct calculation, as the one leading to formula (\ref{shakl}) in the proof of Theorem \ref{prop_STA0_to_STAII}, implies:
\begin{equation}\label{shakl2}
\begin{split}
& (\nabla u_{i+1})^T(\bar g\circ u_{i+1})\nabla u_{i+1} - (\nabla u_i)^T(\bar g\circ u_i)\nabla u_i
\\ & \qquad\qquad\qquad\qquad = -2  [\langle b_i^{pq}, v_{i+1}\rangle]_{pq} + \Big(\big(\nabla v_{i+1})^T\nabla v_{i+1} + 2\,\sym\nabla w_{i+1} \Big)+ \mathcal{R}_{i+1}, \vspace{1mm}
\\ & \mbox{where: } 
\mathcal{R}_{i+1} = 2\,\sym\big[\langle \partial_p u_i, G_i(\partial_q T_i) w_{i+1}\rangle\big]_{pq} + 2\,\sym\big[\langle (\partial_p E_i) v_{i+1}, G_iE_i\partial_q v_{i+1}\rangle\big]_{pq}
\\ & \qquad\qquad\qquad \; + \big[\langle(\partial_pE_i) v_{i+1}, G_i(\partial_q E_i)v_{i+1}\rangle\big]_{pq} + 2\,\sym \Big(\nabla(E_iv_{i+1})^TG_i\nabla (T_i w_{i+1})\Big) 
\\ & \qquad\qquad\qquad \; + \nabla(T_i w_{i+1})^TG_i \nabla (T_i w_{i+1}) \\ & \qquad\qquad\qquad \; - (\nabla u_i)^T (\bar g\circ u_i - G_i)\nabla u_i +
(\nabla u_{i+1})^T \big(\bar g\circ u_{i+1}- G_i\big)\nabla u_{i+1},
\end{split}
\end{equation}
and where $\{b_i^{pq}\in\mathcal{C}^\infty(\bar\omega, \R^k)\}_{p,q=1\ldots d}$, with $b_i^{pq}=b_i^{qp}$, are defined similarly to the first three terms in the second fundamental form formula (\ref{sec_fun}) in step 3 of the proof of Theorem \ref{prop_STA0_to_STAII}:
$$b_i^{pq} = E_i^TG_i\partial_{pq}u_i + \frac{1}{2}E_i^T(\partial_qG_i)\partial_pu_i + \frac{1}{2}E_i^T(\partial_pG_i)\partial_qu_i$$

\smallskip

{\bf 3.} We consider the magnitudes of each of the three leading terms in the expansion (\ref{shakl2}):
$$-2  [\langle b_i^{pq}, v_{i+1}\rangle]_{pq}, \qquad \big(\nabla v_{i+1})^T\nabla v_{i+1} + 2\,\sym\nabla w_{i+1}, \qquad \mathcal{R}_{i+1}. $$
It is clear that the components of the first term $\langle b_i^{pq}, v_{i+1}\rangle$ are as small as needed, provided that $\epsilon_{i+1}$ is chosen small enough, because of the first bound in (\ref{imm5}). The second term will be dealt with using the last assertion in (\ref{imm5}). We now estimate norms of the seven terms in $\mathcal{R}_{i+1}$. 
The first four terms are as small as needed provided that $\epsilon_{i+1}$ is small enough, since each of them contains $v_{i+1}$ or $w_{i+1}$, where we recall (\ref{imm5}). This is clear in the first three terms, whereas for the fourth term we use the identities in (\ref{imm7.5}) to rewrite it as:
\begin{equation*}
\begin{split}
& \sym \Big(\nabla (E_iv_{i+1})^T(G_i\nabla (T_i w_{i+1}))\Big) = 
\sym \big[\big\langle (\partial_pE_i)v_{i+1}, G_i(\partial_q T_i) w_{i+1}\big\rangle\big]_{pq} \\ & \qquad\qquad + \sym \big[\big\langle (\partial_pE_i)v_{i+1}, G_iT_i\partial_q w_{i+1}\big\rangle\big]_{pq} + \sym \big[\big\langle E_i\partial_pv_{i+1}, G_i(\partial_q T_i) w_{i+1}\big\rangle\big]_{pq} 
\end{split}    
\end{equation*}
Regarding the fifth term in $\mathcal{R}_{i+1}$, we may likewise rewrite it as:
\begin{equation*}
\begin{split}
& \nabla (T_iw_{i+1})^TG_i\nabla (T_iw_{i+1}) = (\nabla w_{i+1})^TT_i^TG_iT_i\nabla w_{i+1} \\ & \qquad\qquad + 2\,\sym\big[\big\langle T_i(\partial_pw_{i+1}), G_i(\partial_qT_i)w_{i+1}\big\rangle\big]_{pq}
+ \big[\big\langle (\partial_pT_i) w_{i+1}, G_i(\partial_qT_i)w_{i+1}\big\rangle\big]_{pq}
\end{split}    
\end{equation*}
The first term above is potentially "large", so we do not estimate it:
$$ (\nabla w_{i+1})^TT_i^TG_iT_i\nabla w_{i+1} = (\nabla w_{i+1})^T \big((\nabla u_i)^TG_i\nabla u_i\big)^{-1} \nabla w_{i+1},$$
and only observe that the final two terms are again controlled by the smallness of $\epsilon_{i+1}$. 

\smallskip

\noindent The sixth term in $\mathcal{R}_{i+1}$ is estimated by (\ref{imm7}):
$$ \| (\nabla u_i)^T (\bar g\circ u_i - G_i)\nabla u_i \|_0\leq \|\nabla u_i\|_0^2 \, \eta_{i+1},$$
while for the final, seventh term, we have, by (\ref{imm7.6}) and (\ref{imm7}):
\begin{equation*}
\begin{split}
& \|(\nabla u_{i+1})^T \big(\bar g\circ u_{i+1}- G_i\big)\nabla u_{i+1}\|_0\leq
\|\nabla u_{i+1}\|_0^2 \|\bar g\circ u_{i+1}- G_i\|_0\\ & \leq 
\|\nabla u_{i+1}\|_0^2 \Big(\|\bar g\circ u_{i+1}- \bar g\circ u_{i}\|_0 + \|\bar g\circ u_{i}- G_i\|_0\Big)
\\ & \leq \Big(\epsilon_{i+1}(\|\nabla E_i\|_0+\|\nabla T_i\|_0) + 
C(1+\|\mathcal{D}_0\|_0)\Big)^2
\Big(\omega_{\bar g}(C\epsilon_{i+1}) + \eta_{i+1}\Big),
\end{split}    
\end{equation*}
where $C$ depends only on $\omega, d,k,\gamma,\Lambda$ and $\omega_{\bar g}$ denotes the modulus of continuity of $\bar g$ on $\bar{\mathcal{U}}$.

\smallskip

\noindent In conclusion, we see that upon fixing $\delta_{i+1}$, then taking $\eta_{i+1}$ in \eqref{imm7} sufficiently small, in function of all the so far defined quantities, and further taking $\epsilon_{i+1}$ small, we obtain:
\begin{equation*}
\begin{split}
& \|\mathcal{R}_{i+1} - (\nabla w_{i+1})^T \big((\nabla u_i)^TG_i\nabla u_i\big)^{-1} \nabla w_{i+1}\|_0\leq \frac{\delta_{i+1}}{2},
\\ &  2\big\| [\langle b_i^{pq}, v_{i+1}\rangle]_{pq}\big\|_0\leq \frac{\delta_{i+1}}{2},
\\ & \Big\|2A - \Big((\nabla v_{i+1})^T\nabla v_{i+1}+2\,\sym\nabla w_{i+1}\Big)\Big\|_0\leq 2\delta_{i+1},
\end{split}    
\end{equation*}
We may hence write, in view of (\ref{shakl2}) and recalling that $2A=\frac{1-\theta}{N}\mathcal{D}_0$:
\begin{equation}\label{imm8}
\begin{split}
& (\nabla u_{i+1})^T(\bar g\circ u_{i+1})\nabla u_{i+1} - (\nabla u_i)^T(\bar g\circ u_i)\nabla u_i
\\ & \qquad\qquad\qquad\qquad = \frac{1-\theta}{N}\mathcal{D}_0 +(\nabla w_{i+1})^T \big((\nabla u_i)^TG_i\nabla u_i\big)^{-1} \nabla w_{i+1} + \bar{\mathcal{R}}_{i+1}
\\ & \mbox{where: } \|\bar{\mathcal{R}}_{i+1}\|_0\leq 3\delta_{i+1}.
\end{split}    
\end{equation}

\smallskip

{\bf 4.} We now prove the immersability condition (\ref{imm6}). From (\ref{imm8}) it follows that:
\begin{equation}\label{imm9}
\begin{split}
 (\nabla u_{i+1})^T(\bar g\circ u_{i+1})\nabla u_{i+1}  = \, &
(\nabla u_0)^T(\bar g\circ u_0)\nabla u_0 + (i+1)
\frac{1-\theta}{N}\mathcal{D}_0 \\ & + \sum_{j=0}^{i}(\nabla w_{j+1})^T \big((\nabla u_j)^TG_j\nabla u_j\big)^{-1} \nabla w_{j+1} + \sum_{j=0}^i\bar{\mathcal{R}}_{j+1}.
\end{split}    
\end{equation}
Since each term $(\nabla w_{j+1})^T \big((\nabla u_j)^TG_j\nabla u_j\big)^{-1} \nabla w_{j+1}$ is nonnegative definite by (\ref{imm7.4}), and $\mathcal{D}_0$ is positive definite as well, we obtain by (\ref{imm4.5}) provided that $\sum_{j=1}^N{\delta}_{j}$ is sufficiently small:
$$ (\nabla u_{i+1})^T(\bar g\circ u_{i+1})\nabla u_{i+1} \geq  (\nabla u_0)^T(\bar g\circ u_0)\nabla u_0 + \sum_{j=0}^i\bar{\mathcal{R}}_{j+1}\geq 
\frac{1}{3\gamma\Lambda}\Id_d \quad\mbox{in }\, \bar\omega.$$
For the upper bound, we first write:
\begin{equation}
\label{imm10}
(\nabla u_0)^T(\bar g\circ u_0)\nabla u_0 + (i+1)
\frac{1-\theta}{N}\mathcal{D}_0 = g-  \mathcal{D}_0 
+ (i+1) \frac{1-\theta}{N}\mathcal{D}_0 \leq g-\theta\mathcal{D}_0,
\end{equation}
as $i+1\leq N$, and therefore:
$$ (\nabla u_{i+1})^T(\bar g\circ u_{i+1})\nabla u_{i+1} \leq g -\theta\mathcal{D}_0 + \sum_{j=0}^{i}(\nabla w_{j+1})^T \big((\nabla u_j)^TG_j\nabla u_j\big)^{-1}\nabla w_{j+1}+ \sum_{j=0}^i\bar{\mathcal{R}}_{j+1}.$$
Second, by (\ref{imm7.4}) and (\ref{imm5}), we get, with $C$ depending on $\omega, d,k,\gamma, \Lambda$:
\begin{equation}\label{imm11}
\begin{split}
\sum_{j=0}^i\|(\nabla w_{j+1})^T \big((\nabla u_j)^TG_j\nabla u_j\big)^{-1}\nabla w_{j+1}\|_0 & \leq \sum_{j=0}^i\|\nabla w_{j+1}\|_0^2 \big\|\big((\nabla u_j)^TG_j\nabla u_j\big)^{-1}\big\|_0 \\ & \leq CN \Big(\frac{\|\mathcal{D}_0\|_0}{N}\Big)^2 = C \frac{\|\mathcal{D}_0\|_0^2}{N}.
\end{split}    
\end{equation}
The last quantity is again small, provided that $N$ has been chosen large enough and $\sum_{j=1}^N\delta_j$ is small enough. In conclusion:
\begin{equation}\label{imm_fin}
(\nabla u_{i+1})^T(\bar g\circ u_{i+1})\nabla u_{i+1} \leq g-\frac{\theta}{2}\mathcal{D}_0 \leq \gamma\,\Id_d\leq 3\gamma\Lambda\,\Id_d\quad\mbox{ in }\bar\omega.
\end{equation}
This ends the proof of the inductive estimate (\ref{imm6}).

\smallskip

{\bf 5.} We declare $\underline u = u_N\in\mathcal{C}^\infty(\bar\omega, \mathcal{U})$ and prove the claimed assertions. By (\ref{imm6}) we obtain:
$$\frac{1}{3\gamma\Lambda^2}\Id_d\leq(\nabla\underline u)^T\nabla \underline u \leq 3\gamma\Lambda^2\Id_d \quad\mbox{ in }\; \bar\omega.$$
The second assertion follows by the first bound displayed in (\ref{imm_fin}):
$$ (\nabla \underline u)^T(\bar g\circ \underline u)\nabla \underline u
= (\nabla u_{N})^T(\bar g\circ u_{N})\nabla u_{N} \leq g - \frac{\theta}{2}\mathcal{D}_0<g.$$
The third assertion is due to (\ref{imm4}) and (\ref{imm7.6}), where again the $N$-tuple of $\epsilon_i$ is taken small enough.
For the final, fourth assertion, we use (\ref{imm9}) and (\ref{imm10}) to get:
\begin{equation*}
(\nabla\underline u)^T(\bar g\circ\underline u)\nabla \underline u = g-\theta\mathcal{D}_0 + \sum_{j=0}^{N-1}(\nabla w_{j+1})^T \big((\nabla u_j)^TG_j\nabla u_j\big)^{-1} \nabla w_{j+1} + \sum_{j=0}^{N-1}\bar{\mathcal{R}}_{j+1},  \end{equation*}
which by (\ref{imm11}) and the assumed smallness conditions, implies:
$$\big\| g - (\nabla\underline u)^T(\bar g\circ\underline u)\nabla \underline u\big\|_0 \leq \theta\|\mathcal{D}_0\|_0 + C\frac{\|\mathcal{D}_0\|_0^2}{N}
+ \frac{\delta}{3}\leq \delta.$$
The proof is done.
\end{proof}

\section{Flexibility result for (\ref{II}) and an example} \label{sec_NK_II3}

In this section, we prepare to iterate (\ref{STA-II}) via the Nash-Kuiper scheme of Proposition \ref{thm_NK_II} by first applying Proposition \ref{lem_prestage} to construct a family of referential immersions $\underline u$ whose metric defects are centered at $\delta H_0$. In the proof of the general flexibility result in Theorem \ref{th_final_II}, we then select one such immersion corresponding to a sufficiently small parameter $\underline\delta$. We specialize this result in Corollary \ref{specII}, using the exponents for which (\ref{STA-II}) is currently known to hold. We conclude with Example \ref{counter_ex} in dimension $d=2$ and codimension $k=1$, showing that the strict shortness assumption on $u$ cannot, in general, be weakened to mere shortness.

\medskip

\noindent Recall definition of the matrix $H_0\in\R^{d\times d}_\sym$, given in Lemma \ref{lem_dec_def}. We then have:

\begin{corollary}\label{lem_initial_immersion} 
Under the same assumptions as in Proposition \ref{lem_prestage}, there exists $\underline\gamma\geq 1$ and $\underline\delta_*>0$ 
depending only on the data $g,\bar g, u$ and the dimension $d$,
such that the following holds. For every $\epsilon\in (0,1)$ and $\delta_0\in (0, \underline\delta_*)$ there exists an immersion $\underline u \in\mathcal{C}^\infty(\bar\omega, \mathcal{U})$, satisfying:
\begin{equation}\label{coro_asser}
\begin{split}
& \frac{1}{\underline \gamma} \Id_d\leq (\nabla \underbar u)^T\nabla
\underbar u \leq \underline\gamma\Id_d \quad\mbox{ in } \,\bar\omega,
\\ & \|u-\underline u\|_0\leq\epsilon \quad\mbox{ and }\quad \mathrm{dist}(\underline u(\bar\omega), \partial\mathcal{U})\geq \frac{3}{4} \mathrm{dist}(u(\bar\omega), \partial\mathcal{U}), 
\\ &  \big\|g - \delta_0 H_0 - (\nabla \underline u)^T(\bar g\circ \underline u)\nabla \underline u\big\|_0\leq \frac{r_0}{3}\delta_0.
\end{split}
\end{equation}
\end{corollary}
\begin{proof}
We choose $\underline\delta_*\in (0,1)$ sufficiently small to ensure that:
\begin{equation}\label{grain}
g-\underline\delta_* H_0\geq \frac{1}{2}g \quad \mbox{ and } \quad
g-(\nabla u)^T(\bar g\circ u)\nabla u >
\underline\delta_* H_0\quad\mbox{ in }\; \bar\omega.
\end{equation}
Fix $\epsilon<\mathrm{dist}(u(\bar\omega), \partial\mathcal{U})/4$
and $\delta_0<\underline\delta_*$, and apply Proposition
\ref{lem_prestage} with the metrics:
$g- \delta_0 H_0$ on $\bar\omega$ and $\bar g$ on
$\bar{\mathcal{U}}$, the immersion $u$,  and parameters $\epsilon, r_0\delta_0/3$.
We obtain $ \underline u\in \mathcal{C}^\infty(\bar\omega,\mathcal U)$
such that:
$$ \|u - \underline u\|_0\leq\epsilon, \qquad 
\| g - \delta_0 H_0 - (\nabla\underline u)^T
(\bar g\circ\underline u)\nabla\underline u \|_0
\leq\frac{r_0}{3}\delta_0,$$
where the distance assertion is also valid in view of the first
condition above. Finally, we observe that the immersability constant $\underline\gamma$
of $\underline u$ depends only on the ellipticity constant of $g-\delta_0 H_0$ (which is independent of $\delta_0$ by the first condition in (\ref{grain})), the ellipticity constant of $\bar g$, and the immersability constant of $u$. The proof is done.
\end{proof}

\bigskip

\noindent Since (\ref{STA-INI0}) holds by \cite[Theorem 6.2]{lew_conv}, we thus obtain the following flexibility result for (\ref{II}), with the regularity range (\ref{range}) determined by the integers $J,S$ in (\ref{STA0}) and the regularity of the Riemannian metrics $g$ and $\bar g$:

\begin{theorem}\label{th_final_II}
Assume that (\ref{STA0}) holds with some integers $J,S$.
Let $g\in\mathcal{C}^{r,\beta}(\bar\omega, \R^{d\times d}_{\sym,>})$ and $\bar g\in\mathcal{C}^{\bar r, \bar\beta}(\bar{\mathcal{U}},\R^{(d+k)\times(d+k)}_{\sym, >})$ be two Riemannian metrics, with some $r+\beta, \bar r+\bar \beta>0$, where $\omega\subset\R^d$ is open, bounded, contractible, with smooth boundary, and $\mathcal{U}\subset\R^{d+k}$ is open and bounded.
Let $u\in\mathcal{C}^1(\bar\omega, \mathcal{U})$ be an immersion, satisfying:
\begin{equation}\label{imm_str_short}
(\nabla u)^T(\bar g\circ u)\nabla u< g \quad \mbox{in }\; \bar\omega.
\end{equation}
Then, for every $\epsilon>0$ and for every $\alpha$ in
the range:
$$ 0<\alpha<\min\Big\{\frac{1}{1+2J/S}, \frac{r+\beta}{2}, \frac{\bar r+\bar\beta}{2}\Big\}$$
there exists $\bar u\in\mathcal{C}^{1,\alpha}(\bar\omega,\mathcal{U})$, such that:
\begin{equation*}
\|\bar u - u\|_0\leq \epsilon \quad \mbox{and} \quad (\nabla
\bar u)^T(\bar g\circ \bar u)\nabla \bar u =g \quad \mbox{in }\;\bar\omega. 
\end{equation*}
\end{theorem}
\begin{proof}
Let $g,\bar g, u$ be as in the statement and fix $\epsilon\in (0,1)$ and the H\"older regularity exponent $\alpha$ in the indicated range. 
Extending $g,u$ to a neighborhood of $\bar\omega$, and decreasing $\tau_0$ in (\ref{tau0}) if necessary, we may assume that $u$ remains a strictly short immersion into $\mathcal{U}$ on $\bar\omega_{2\tau_0}$. 
We fix those extensions, together with $\tau_0$, before choosing $\delta_0$.
Denote by $\Lambda, d_\mathcal{U}$ the quantities in:
\begin{equation*}
\begin{split}
& \frac{1}{\Lambda}\Id_{d+k}\leq \bar g\leq \Lambda\,\Id_{d+k}\quad\mbox{ in }\, \bar{\mathcal{U}},    
\\ & d_\mathcal{U} = \frac{3}{4}\mathrm{dist}(u(\bar\omega_{2\tau_0}),\partial\mathcal{U}),
\end{split}    
\end{equation*}
and denote by $\underline\delta_*>0$, $\underline\gamma\geq 1$ constants obtained through Corollary \ref{lem_initial_immersion} applied on $\omega_{2\tau_0}$, in  view of the validity of (\ref{STA-INI0}). By (\ref{STA0}) and Theorem \ref{prop_STA0_to_STAII} applied on $\omega$, the stage (\ref{STA-II}) holds, with $J,S$ as in (\ref{STA0}), with $p=S+1$, and with  $C_1$ depending only on $\omega, \underline\gamma, \Lambda, d,k,J,S$.

\medskip

\noindent We apply Proposition \ref{thm_NK_II} on $\omega$, with the parameters set above and with a fixed exponent $\alpha$ in the indicated range, to obtain the threshold parameter $\underline\delta>0$ for the Nash-Kuiper iteration. We observe in passing that if either metric $g, \bar g$ has regularity greater than $2$, then we may use $2$ as its regularity; this does not change the admissible range in (\ref{range}), since trivially $1/(1+2J/S)<1$.
We subsequently apply Corollary \ref{lem_initial_immersion} on $\omega_{2\tau_0}$, with parameters:
$$\frac{\epsilon}{2} \qquad \mbox{ and }\qquad \delta_0<\min\{\underline\delta_*, \underline\delta\} \quad\mbox{satisfying} \quad 2C_1\delta^{1/2}_0\leq \frac{\epsilon}{2},$$
to find a referential immersion $\underline u\in\mathcal{C}^\infty(\bar\omega_{2\tau_0}, \mathcal{U})$, with the following properties:
\begin{equation*}
\begin{split}
& \frac{1}{\underline\gamma}\Id_{d}\leq (\nabla\underline u)^T\nabla\underline u\leq \underline\gamma\,\Id_{d}\quad\mbox{ in }\, \bar\omega_{2\tau_0},    
\\ & \|u-\underline u\|_0\leq \frac {\epsilon}{2} \quad \mbox{ on }\, \bar\omega_{2\tau_0}, \qquad \mathrm{dist}(\underline u(\bar\omega_{2\tau_0}),\partial\mathcal{U})\geq d_\mathcal{U},
\\ & \big\|g-\delta_0H_0 - (\nabla \underline u)^T(\bar g\circ \underline u)\nabla \underline u\big\|_0\leq \frac{r_0}{3}\delta_0 \quad \mbox{ on }\, \bar\omega_{2\tau_0}.
\end{split}    
\end{equation*}
These are exactly the assumptions (i), (ii), (iv) in Proposition \ref{thm_NK_II}. Consequently, there exists $\bar u\in\mathcal{C}^{1,\alpha}(\bar\omega, \mathcal{U})$ with $(\nabla \bar u)^T(\bar g\circ \bar u)\nabla \bar u =g$ and with:
$$\|\bar u - u\|_0\leq \|\underline u - u\|_0 + \|\bar u - \underline u\|_0
\leq \frac{\epsilon}{2} + 2C_1\delta_0^{1/2}\leq \epsilon. $$
The proof is done.
\end{proof}

\medskip

\noindent In conclusion, the versions of (\ref{STA0}) proved in \cite{CHIL, in_lew2}, the analysis done directly for the system (\ref{II}) in codimension $k=1$ in \cite{CHI2, In2026}, as well as the forthcoming result \cite{CIL_new} imply, in virtue of Theorem \ref{th_final_II}, the current state of the art regularity  for flexibility of the system (\ref{II}):

\begin{corollary}\label{specII}
Let $g\in\mathcal{C}^{r,\beta}(\bar\omega, \R^{d\times d}_{\sym,>})$ and $\bar g\in\mathcal{C}^{\bar r, \bar\beta}(\bar{\mathcal{U}},\R^{(d+k)\times(d+k)}_{\sym, >})$ be two Riemannian metrics, with some $r+\beta, \bar r+\bar \beta>0$, where $\omega\subset\R^d$ is open, bounded, contractible, with smooth boundary, and $\mathcal{U}\subset\R^{d+k}$ is open and bounded.
Let $u\in\mathcal{C}^1(\bar\omega, \mathcal{U})$ be an immersion, satisfying:
\begin{equation}\label{short_again}
(\nabla u)^T(\bar g\circ u)\nabla u< g \quad \mbox{in }\; \bar\omega.
\end{equation}
Then, for every $\epsilon>0$ and for every $\alpha$ in the range:
$$ 0<\alpha<\min\Big\{\frac{1}{2(d-s)+1}, \frac{r+\beta}{2}, \frac{\bar r+\bar\beta}{2}\Big\}
\quad\mbox{ when }\quad k\in\big[(s-1)_*+1,s_*+1\big),$$
where $s=1\ldots d$, there exists $\bar u\in\mathcal{C}^{1,\alpha}(\bar\omega,\mathcal{U})$, such that:
\begin{equation*}
\|\bar u - u\|_0\leq \epsilon \quad \mbox{and} \quad (\nabla
\bar u)^T(\bar g\circ \bar u)\nabla \bar u =g \quad \mbox{in }\;\bar\omega. 
\end{equation*}
\end{corollary}

\bigskip

\noindent We close this section by examining whether the strict shortness condition (\ref{short_again}) is essential, or whether it can be relaxed to the mere shortness, namely:
\begin{equation}\label{imm_short}
(\nabla u)^T(\bar g\circ u)\nabla u \leq g \quad \mbox{in }\; \bar\omega.
\end{equation}
It is certainly so when either  $g$ or $\bar g$ is Euclidean. Indeed, when $\bar g=\Id_{d+k}$ then any short immersion $u$ satisfying (\ref{imm_short}) may be approximated by $u_\epsilon = (1-\epsilon) u$ which is strictly short:
$$g - (\nabla u_\epsilon)^T\nabla u_\epsilon = g- (1-\epsilon)^2 (\nabla u)^T\nabla u = (1-\epsilon)^2(g-(\nabla u)^T\nabla u) + (2\epsilon - \epsilon^2)g >0,$$
for small $\epsilon>0$. In case $g=\Id_d$ and $\omega= B_1$, we define $u_\epsilon(x) = u((1-\epsilon)x)$ and similarly get:
$$\Id_d - (\nabla u_\epsilon)^T(\bar g\circ u_\epsilon)\nabla u_\epsilon = \Id_d- (1-\epsilon)^2 \big((\nabla u)^T(\bar g\circ u)\nabla u\big)_{\mid (1-\epsilon )x} >0.$$
The same construction works on any domain $\omega$ that is star-shaped about a point. However, in the general case of arbitrary Riemannian metrics $g$ and $\bar g$, neither the approximation argument as above nor the actual assertion of Theorem \ref{th_final_II} or Corollary \ref{specII}  are valid under the weakened assumption \eqref{imm_short}, because of the following example:

\smallskip

\begin{example}\label{counter_ex}\footnote{Suggested by ChatGPT.}
Let $d=2, k=1$, $\omega=B_1\subset\R^2$, $\mathcal{U}=B_2\times (-1,1)\subset\R^3$ and declare, writing in polar $(r,\theta)$ or in cylindrical $(r,\theta,z)$ coordinates, respectively:
\begin{equation*}
\begin{split}
& g(r,\theta)= \big(1+ \chi(r)(2+\cos\theta)\big)dr^2 + a(r)d\theta^2,
\\ & \bar g(r,\theta, z) = dr^2 + a(r) d\theta^2 + dz^2,
\\ & u(r,\theta) = (r,\theta,0).
\end{split}
\end{equation*}
Above, $\chi(r)$ is a smooth, non-decreasing cut-off function, interpolating between $0$ on the interval $[0, 1/4]$ and $1$ on $[1/2, 1]$, and $a(r)$ is a smooth non-decreasing function, interpolating between $r^2$ on $[0,1/4]$ and $1$ on $[1/2, 2]$. In particular, $g$ and $\bar g$ are smooth and there holds:
\begin{equation}\label{fori}
\begin{split}
& g= \Id_2, \quad \bar g=\Id_3 \quad\mbox{ for } r\leq 1/4, 
\\ & g= (3+\cos\theta)dr^2 + d\theta^2, \quad \bar g = dr^2+d\theta^2+dz^2 
\quad\mbox{ for } r\geq 1/2. 
\end{split}
\end{equation}
We directly check that $u$ satisfies the non-strict shortness condition (\ref{imm_short}):
$$(\nabla u)_{\mid(r,\theta)}^T \bar g_{\mid(r,\theta,0)}(\nabla u)_{\mid(r,\theta)} = dr^2 + a(r) d\theta^2 
\leq \big(1+\chi(r)(2+\cos\theta)\big) dr^2+ a(r)d\theta^2= g_{\mid (r,\theta)}.$$
Assume that $u$ has a nearby exact isometric immersion $\tilde u$, written in cylindrical coordinates as $(\tilde u^r, \tilde u^\theta, \tilde u^z)$, so that the image $\tilde u(\mathcal{A})$
of the annulus $\mathcal{A}=\{r\in (1/2+1/5, 1-1/5)\}\subset B_1$ lies in the region with the radial coordinate $\tilde u^r\geq 1/2$. By the second line of (\ref{fori}), it follows that:
$$(\nabla \tilde u)_{\mid(r,\theta)}^T (\nabla \tilde u)_{\mid(r,\theta)} = \big((\nabla \tilde u)^T (\bar g\circ \tilde u)\nabla \tilde u\big)_{\mid(r,\theta)} = g_{\mid(r,\theta)}= (3+\cos\theta) dr^2+ d\theta^2 \quad\mbox{in } \mathcal{A}. $$
Equivalently, we may write on $\mathcal{A}$:
\begin{equation}\label{coor}
|\partial_r\tilde u|^2 = 3+\cos\theta, \qquad \langle \partial_r\tilde u, \partial_\theta\tilde u \rangle = 0, \qquad |\partial_\theta\tilde u|^2 =1.       
\end{equation}
The third equality implies $\partial_\theta\tilde u^\theta \leq |\partial_\theta \tilde u^\theta |\leq 1$. Since $\tilde u$ uniformly approximates the identity map $u$, the closed curve $\theta\mapsto \tilde u(r,\theta)$ winds around the origin once, so that  $\tilde u^\theta (r,2\pi) - \tilde u^\theta(r,0)=2\pi$, implying $\partial_\theta\tilde u^\theta =1$. Consequently, $\tilde u^r$ and $\tilde u^z$ are independent from $\theta$, because then:
$$ \partial_\theta\tilde u = (0,1,0)\quad\mbox{in } \mathcal{A}. $$
The second equality in (\ref{coor}) now yields $\partial_r\tilde u^\theta=0$, so in view of the first equality we obtain:
$$ |\partial_r\tilde u^r|^2 + |\partial_r\tilde u^z|^2 = 3+ \cos\theta \quad\mbox{in } \mathcal{A}. $$
This contradicts the independence of $\tilde u^r$ and $\tilde u^z$ from $\theta$ and hence contradicts existence of $\tilde u$. 
\end{example}

\part{(\ref{STA0}) $\Rightarrow$ (\ref{STA-VK}) and flexibility results for (\ref{VK}) }

\section{Stage (\ref{STA0}) implies stage (\ref{STA-VK})}\label{seci2}

In this section, we show that the basic stage (\ref{STA0}), formulated for the Monge-Amp\`ere system (\ref{MA}) in its weak form (\ref{VK}), relative to the  zero initial displacements but including a lower-order term $b_0$, directly implies the stage (\ref{STA-VK}) for arbitrary initial  $v,w$.
The key calculation is (\ref{shakl3}). It shows that correcting $v$ by $\tilde v$ and $w$ by $-(\nabla v)^T\tilde v+\tilde w$, where $\tilde v,\tilde w$ are obtained by applying (\ref{STA0}), updates the von K\`arm\`an content by precisely the content associated to $\tilde v,\tilde w$, plus a lower-order term in $\tilde v$ whose coefficients are $b_0^{pq}=\partial_{pq}v$.

\medskip

\noindent We then recall two results from \cite{lew_conv}, namely the pre-stage and the Nash-Kuiper iteration scheme, which yield the convex integration result for (\ref{MA}) in the regularity range (\ref{ran_VK}).

\begin{theorem}\label{prop_STA0_to_STAVK}
Assume that (\ref{STA0}) holds, with some $J, S\geq 1$. Then (\ref{STA-VK}) holds
as well, with the same exponents $J,S$.
\end{theorem}  
\begin{proof}
{\bf 1.} Let $\omega, \tau, A, v, w$ be as in the statement on (\ref{STA-VK})
and let $\delta_0, \sigma_0, J,S,L$ be as in (\ref{STA0}).
Fix also positive parameters $\mathcal{M}, \sigma$, which together
with $\mathcal{D} = A-\big(\frac{1}{2}(\nabla v)^T\nabla v + \sym\nabla w\big)$
satisfy the assumptions of (\ref{STA-VK}), where we define:
\begin{equation}\label{redef22}
\underline \delta = \frac{\delta_0}{{\bar C}^2 +\dbar C}, \qquad \underline
\sigma = \sigma_0,
\end{equation}
with respect to the constants $\bar C, \dbar C$ depending only on $d,k, J, S, L$ and given below. Set:
$$\delta=({\bar C}^2 +\dbar C)\|\mathcal{D}\|_0, \qquad \lambda =
\frac{\mathcal{M}}{\|\mathcal{D}\|_0^{1/2}},$$ 
and define the following regularized quantities, which on  $\bar\omega_{\tau}$
use only the data on $\bar\omega_{2\tau}$:
$$\bar A = A\ast\phi_{1/\lambda}, \qquad v_0=v\ast\phi_{1/\lambda},
\qquad w_0=w\ast\phi_{1/\lambda},$$
in virtue of $1/\lambda\leq \underline\delta^{1/2}\tau\leq\tau$. By Lemma \ref{lem_stima}, we further obtain:
\begin{equation}\label{uno2}
\begin{split}  
& \|\bar A-A\|_0\leq \frac{\|A\|_{r,\beta}}{\lambda^{r+\beta}} =
\frac{\|A\|_{r,\beta}}{{\mathcal{M}}^{r+\beta}} \|\mathcal{D}\|_0^{(r+\beta)/2},
\\ & \|v_0-v\|_1 + \|w_0-w\|_1\leq \frac{\|v\|_2 + \|w\|_2}{\lambda}\leq
2\frac{\mathcal{M}}{\lambda} = 2\|\mathcal{D}\|_0^{1/2},
\\ & \|\nabla^{(m+2)}v_0\|_0+ \|\nabla^{(m+2)}w_0\|_0\leq C
\mathcal{M} \lambda^m = \bar C\lambda^{m+1} \|\mathcal{D}\|_0^{1/2} \quad\mbox{for
  all } m=0\ldots L. 
\end{split}
\end{equation}  
with $\bar C$ that depends on $d, k, L$. We now define:
$$\quad A_0 = \bar A -
\big(\frac{1}{2}(\nabla v_0)^T\nabla v_0 + \sym\nabla w_0\big) +\delta H_0, \qquad 
b_0^{pq}=\partial_{pq}v_0 \quad\mbox{ for all } p,q=1\ldots d.$$
Using Lemma \ref{lem_stima} again, we arrive at:
\begin{equation*}
\begin{split}
 \|A_0-\delta H_0\|_0 & \leq \|\mathcal{D}\ast\phi_{1/\lambda}\|_0+ 
\|(\nabla v_0)^T\nabla v_0 - ((\nabla v)^T\nabla v)\ast\phi_{1/\lambda}\|_0\\ & \leq \|\mathcal{D}\|_0 + \lambda^{-2}\mathcal{M}^2\leq \frac{1}{2}r_0\dbar C\|\mathcal{D}\|_0
\end{split}
\end{equation*}
where $\dbar C$ depends on $d,k,L$ is chosen appropriately, and additionally ensures that:
\begin{equation*}
\begin{split}
\|\nabla^{(m)}A_0\|_0&\leq \|\nabla^{(m)}(\mathcal{D}\ast\phi_{1/\lambda} )\|_0
+ \|\nabla^{(m)}\big((\nabla v_0)^T\nabla v_0 - ((\nabla v)^T\nabla v)\ast\phi_{1/\lambda}\big)\|_0
\\ & \leq C\big(\lambda^m \|\mathcal{D}\|_0+
\lambda^{m-2}\mathcal{M}^2\big) \leq \dbar C\lambda^m \|\mathcal{D}\|_0 \mbox{ for all }\; m=1\ldots L.
\end{split}
\end{equation*}  
Recalling (\ref{redef22}), we conclude:
\begin{equation*}
\begin{split}
& \|A_0 - \delta H_0\|_0\leq \frac{1}{2}r_0\delta\quad 
\mbox{ and }\quad \|\nabla^{(m)} A_0\|_0\leq \lambda^m\delta \quad\mbox{ for all } \; m=1\ldots L,
\\ & \|\nabla^{(m)}b_0^{pq}\|_0\leq \lambda^{m+1}\delta^{1/2}
\quad\mbox{ for all } m=0\ldots L, \; p,q=1\ldots d,
\\ & \delta\leq \delta_0, \qquad \lambda\delta^{1/2}\geq 1, \qquad \sigma\geq \sigma_0.
\end{split}    
\end{equation*}
Therefore, (\ref{STA0}) applied on $\bar\omega_\tau$ implies existence of
$\bar v\in \mathcal{C}^2(\bar\omega_\tau,\R^k)$, $\bar
w\in\mathcal{C}^2(\bar\omega_\tau,\R^d)$, with:
\begin{equation}\label{quattro2}
\begin{split}
&\lambda \|\bar v\|_0 + \|\nabla \bar v\|_0\leq C\delta^{1/2},\qquad
\lambda \|\bar w\|_0 + \|\nabla \bar w\|_0\leq C\delta, \\
& \|\nabla^{2}\bar v\|_0\leq C\lambda\sigma^J\delta^{1/2}, \qquad
\|\nabla^{2}\bar w\|_0\leq C\lambda\sigma^J\delta, \\
&\big\|A_0 + [\langle b_0^{pq},\bar v\rangle]_{pq} -
\big(\frac{1}{2}(\nabla \bar v)^T\nabla \bar v +\sym\nabla \bar
w\big)\big\|_0 \leq C \frac{\delta}{\sigma^S},
\end{split}
\end{equation}
where the constant $C$ depends on $d,k,J,S$.

\smallskip

{\bf 2.} Define now $\tilde v\in\mathcal{C}^2(\bar\omega_\tau, \R^{k})$,
$\tilde w\in\mathcal{C}^2(\bar\omega_\tau, \R^{d})$ by:
$$\tilde v = v_0+\bar v, \qquad  \tilde w= w_0-(\nabla v_0)^T\bar v +\bar w - \delta H_0x$$
By a direct calculation, we obtain:
\begin{equation}\label{shakl3}
\begin{split}
&\big(\frac{1}{2}(\nabla \tilde v)^T\nabla \tilde v + \sym\nabla\tilde
w\big) - \big(\frac{1}{2}(\nabla v_0)^T\nabla v_0 + \sym\nabla w_0\big)
\\ & = \frac{1}{2}(\nabla\bar v)^T\nabla\bar v + \sym\nabla\bar w + \sym \big((\nabla
v_0)^T\nabla\bar v\big) -\sym\nabla \big((\nabla v_0)^T\bar v\big) -\delta H_0
\\ & = \big( \frac{1}{2}(\nabla\bar v)^T\nabla\bar v +
\sym\nabla\bar w \big) - \big[\langle\partial_{pq}v_0, \bar v\rangle\big]_{pq} - \delta H_0,
\end{split}
\end{equation}
which implies, recalling $\partial_{pq}v_0 = b_0^{pq}$, that:
$$A- \big( \frac{1}{2}(\nabla\tilde v)^T\nabla\tilde v +
\sym\nabla\tilde w \big) =  (A-\bar A) + \Big(A_0 + \big[\langle b_0^{pq}, \bar v\rangle\big]_{pq}
- \big( \frac{1}{2}(\nabla\bar v)^T\nabla\bar v + \sym\nabla\bar w \big) \Big).$$
Thus, by (\ref{uno2}) and (\ref{quattro2}) we conclude:
$$\|A- \big( \frac{1}{2}(\nabla\tilde v)^T\nabla\tilde v +
\sym\nabla\tilde w \big)\|_0\leq C\frac{\delta}{\sigma^S} + \|A-\bar
A\|_0\leq C\frac{\|\mathcal{D}\|_0}{\sigma^S}+
\frac{\|A\|_{r,\beta}}{{\mathcal{M}}^{r+\beta}} \|\mathcal{D}\|_0^{(r+\beta)/2}$$
where the constant $C$ depends on $\omega, d,k,J,S$, as claimed.

\smallskip

{\bf 3.} Finally, by (\ref{uno2}), (\ref{quattro2}) we get:
\begin{equation*}
\begin{split}
& \|\tilde v - v\|_1 \leq \|v_0-v\|_1+\|\bar v\|_1\leq \|\mathcal{D}\|_0^{1/2}
+ C\delta^{1/2} \leq C \|\mathcal{D}\|_0^{1/2}, \\
& \|\tilde w - w\|_1 \leq \|w_0-w\|_1 + \|\bar w\|_1+ \|(\nabla
v_0)^T\bar v\|_1 +\delta \|H_0x\|_1\\ & \qquad\qquad \, \leq \|\mathcal{D}\|_0^{1/2}
+ C\delta + C(\|\mathcal{D}\|_0^{1/2} + \|\nabla v\|_0)\delta^{1/2} +
\bar C\lambda \|\mathcal{D}\|_0^{1/2}\frac{\delta^{1/2}}{\lambda} + C\delta \\ & \qquad\qquad \, \leq C \|\mathcal{D}\|_0^{1/2} (1+\|\nabla v\|_0),
\end{split}
\end{equation*}    
where $C$ depends on $\omega, d,k,J,S$. It remains to estimate, similarly, the second derivatives:
\begin{equation*}
\begin{split}
& \|\nabla^{(2)}\tilde v\|_0 \leq \|\nabla^{(2)} v_0\|_0 +
\|\nabla^{(2)}\bar v\|_0\leq \bar C \lambda \|\mathcal{D}\|_0^{1/2}+
C\lambda\sigma^J\delta^{1/2} \leq C\lambda\sigma^J
\|\mathcal{D}\|_0^{1/2} = C\mathcal{M}\sigma^J,
\\ & \|\nabla^{(2)}\tilde w\|_0\leq \|\nabla^{(2)}w_0\|_0 +
\|\nabla^{(2)}\bar w\|_0 + \|\nabla^{(2)}\big((\nabla v_0)^T\bar v\big)\|_0
\\ & \qquad \leq \bar C \lambda \|\mathcal{D}\|_0^{1/2}+
C\lambda\sigma^J\delta + C\big(\lambda^2 \|\mathcal{D}\|_0^{1/2}\frac{\delta^{1/2}}{\lambda}
+ \lambda \|\mathcal{D}\|_0^{1/2}\delta^{1/2} +
(\|\mathcal{D}\|_0^{1/2}+\|\nabla v\|_0)\lambda\sigma^J\delta^{1/2}\big)
\\ & \qquad \leq C\lambda\sigma^J \|\mathcal{D}\|_0^{1/2}(1+\|\nabla v\|_0)
= C\mathcal{M}\sigma^J (1+\|\nabla v\|_0),
\end{split}
\end{equation*}    
again with $C$ depending on $\omega, d,k,J,S$.  This ends the proof.
\end{proof}

\medskip

\noindent  We now recall the pre-stage and the Nash-Kuiper iteration scheme for the von K\`arm\`an system (\ref{VK}), which are analogous to the results for the isometric immersion system (\ref{II}) in sections \ref{sec_NK_II} -\ref{sec_NK_II3}. They allow the stage (\ref{STA-VK}) to be iterated, yielding the flexibility statement in Corollary \ref{th_final_VK} below. These results were established in \cite{lew_conv}, with one qualification: \cite[Theorem 5.1]{lew_conv} gives Proposition \ref{th_NK_VK} for $r=0$, while a straightforward inspection of its proof shows that it also applies to $r=1$, as long as $r+\beta\leq 2$. We also note that the following is a version of (\ref{STA-INI0}), written for non-zero referential displacements $v, w$:

\begin{proposition}\label{lem_prestage_VK}
Let $A\in \mathcal{C}(\bar\omega, \R^{d\times d}_{\sym})$
be defined on the closure of an open, bounded domain $\omega\subset\R^d$ 
and let $v\in\mathcal{C}^1(\bar\omega, \R^{k})$,
$w\in\mathcal{C}^1(\bar\omega, \R^{d})$ satisfy:
$$\frac{1}{2}(\nabla v)^T\nabla v + \sym\nabla w<A \quad\mbox{ in } \; \bar\omega.$$ 
Then, for every $\epsilon, \delta >0$ there exists $\bar
v\in\mathcal{C}^\infty(\bar\omega,\R^{k})$, $\bar
w\in\mathcal{C}^\infty(\bar\omega,\R^{d})$, satisfying:
\begin{equation*}
\begin{split}
& \|\bar v-v\|_0 + \|\bar w- w\|_0\leq \epsilon \quad\mbox{ and }\quad
\big\|A - \big(\frac{1}{2}(\nabla \bar v)^T\nabla \bar v +
\sym\nabla\bar w\big)\big\|_0\leq \delta,
\\ & \|\nabla \bar v\|_0\le C(1+\|A - (\frac{1}{2}(\nabla
v)^T\nabla v + \sym\nabla w)\|_0 + \|\nabla v\|_0\big),
\end{split}
\end{equation*}
where $C$ depends only on $\omega, d,k$.
\end{proposition}

\begin{proposition}\label{th_NK_VK}
Assume (\ref{STA-VK})  on an open, bounded 
$\omega\subset\R^d$, for some $J, S\geq
1$. Then, for every $A\in\mathcal{C}^{r,\beta}(\bar\omega, \R^{d\times
  d}_\sym)$ with $r+\beta\in (0, 2]$, 
and $v\in\mathcal{C}^2(\bar\omega,\R^k)$,
$w\in\mathcal{C}^2(\bar\omega,\R^d)$ such that:
$$ 0< \|\mathcal{D}\big\|_0\leq \underline\delta \quad \mbox{ where }
\; \mathcal{D} = A-\big(\frac{1}{2}(\nabla v)^T\nabla v + \sym\nabla w\big),$$
and for every $\alpha$ in the range:
\begin{equation}\label{ran_VK}
0< \alpha <\min\Big\{\frac{r+\beta}{2},\frac{1}{1+2J/S}\Big\},
\end{equation}
there exist $\tilde v\in\mathcal{C}^{1,\alpha}(\bar\omega,\R^k)$,
$\tilde w\in\mathcal{C}^{1,\alpha}(\bar\omega,\R^d)$ with the following properties:
\begin{align*} 
& \|\tilde v - v\|_1\leq C\|\mathcal{D}\|_0^{1/2}, \quad \|\tilde w -
w\|_1\leq C\|\mathcal{D}\|_0^{1/2}(1+\|\nabla v\|_0), \vspace{1mm}\\
& \frac{1}{2}(\nabla \tilde v)^T\nabla \tilde v + \sym\nabla
\tilde w =A \quad\mbox{ in }\; \bar\omega, 
\end{align*}
where the constants $C$ depend only on $\omega, d, k, J, S, r,\beta$.
\end{proposition}

\smallskip

\noindent In conclusion, there follows flexibility of the system (\ref{VK}) in the regularity range (\ref{ran_VK}):

\begin{corollary}\label{th_final_VK}
Assume that (\ref{STA0}) holds. Then, for every
$A\in\mathcal{C}^{r,\beta}(\bar\omega, \R^{d\times d}_\sym)$
defined on an open, bounded $\omega\subset\R^d$, with $r+\beta\in (0,2]$, 
for every $v\in\mathcal{C}^1(\bar\omega,\R^{k})$, $w\in\mathcal{C}^1(\bar\omega,\R^d)$
such that:
\begin{equation}\label{st_sho} 
\frac{1}{2} (\nabla v)^T\nabla v + \sym\nabla w < A \quad\mbox{ in
} \bar\omega,
\end{equation}
for every $\epsilon>0$ and $\alpha$ in (\ref{ran_VK}),
there exist  $\tilde v\in\mathcal{C}^{1,\alpha}(\bar\omega,\R^{k})$,
$\tilde w\in\mathcal{C}^{1,\alpha}(\bar\omega,\R^d)$
satisfying:
\begin{equation*}
\|\tilde v - v\|_0 + \|\tilde w- w\|_0\leq \epsilon
\quad\mbox{ and }\quad \frac{1}{2}(\nabla \tilde v)^T\nabla \tilde v + \sym\nabla
\tilde w = A \quad\mbox{ in }\; \bar\omega.
\end{equation*}
\end{corollary}

\bigskip

\noindent We remark that the strict shortness condition (\ref{st_sho}) in the Monge-Amp\`ere setting above can always be replaced by the corresponding shortness condition, namely:
\begin{equation}\label{sho_only}
\frac{1}{2} (\nabla v)^T\nabla v + \sym\nabla w \leq A \quad\mbox{ in
} \bar\omega,
\end{equation}
Indeed, any displacement pair $(v,w)$ satisfying (\ref{sho_only}) may be approximated by the pair $(v_\epsilon, w_\epsilon) = (v,  w - \epsilon \mathrm{id}_d)$ which satisfies (\ref{st_sho}). Hence, no counterpart of Example \ref{counter_ex} is possible here.

\part{(\ref{STA-VK}) $\Rightarrow$ (\ref{STA0}) and validation of (\ref{STA-VK0})}

\section{Deducing stage (\ref{STA0}) from the base form of (\ref{STA-VK})}\label{seci3}

In the preceding sections, we proved that \eqref{STA0} implies both \eqref{STA-II} and \eqref{STA-VK}. We now establish the converse implication under additional assumptions on the construction underlying \eqref{STA-VK}.
The proofs of \eqref{STA-VK} or \eqref{STA-II} available in the literature (see \cite{CHI, CHI2, CHIL, CS, CDS, DIS1/5, In2026, in_lew, in_lew2, Kallen, Kuiper, lew_conv, lew_improved, lew_improved2, lew_full2d2k, LL_last, lewpak_MA, Nash1, Zhao}) mostly concern the special case of $\bar g=\Id_{d+k}$. They proceed by decomposing a suitably scaled defect into finitely many primitive metrics, followed by the successive addition of perturbations, each designed to correct one primitive metric up to a smaller error. The essential properties of these constructions are recorded in the following version of \eqref{STA-VK}:

\begin{equation}\tag{STA-VK0}\label{STA-VK0}
\left[\quad  \begin{minipage}{12.17cm} 
Let $\omega \subset\R^d$ be open and bounded. There exist integers $J,S\geq 1$, constants $\delta_0\in (0,1), \sigma_0>1$ and $L\in \mathbb{N}$, depending only on $\omega,d,k,J,S$, such that the following holds. Given the fields $A_0\in\mathcal{C}^L(\bar\omega,\R^{d\times d}_\sym)$, $v_0\in \mathcal{C}^{L+2}(\bar\omega,\R^{k})$, $w_0\in \mathcal{C}^{L+2}(\bar\omega,\R^{d})$ and defining $\mathcal{D}_0 = A_0- \big(\frac{1}{2}(\nabla v_0)^T\nabla v_0 +\sym\nabla w_0\big)$, 
and given any parameters $\delta, \lambda, \sigma$,  satisfying the bounds:
\begin{equation}\label{bd}
\begin{split}
& \delta\leq \delta_0,\qquad \lambda\delta^{1/2}\geq 1,
\qquad\sigma\geq \sigma_0, 
\\ & \|\mathcal{D}_0 - \delta H_0\|_0\leq \frac{1}{2}r_0\delta,
\qquad \|\nabla^{(m)} \mathcal{D}_0\big)\|_0\leq \lambda^m\delta \quad\mbox{for }\; m=1\ldots L,
\\ & \max\big\{\|\nabla^{(m+2)} v_0\|_0,
\|\nabla^{(m+2)} w_0\|_0\big\}\leq \lambda^{m+1}\delta^{1/2} \quad\mbox{for }\; m=0\ldots L, 
\end{split}
\end{equation}
the stage construction produces finitely many perturbation displacement pairs $\{(V_i,W_i)\}_{i=1}^P$, such that when defining:
\begin{equation}\label{fo}
v_i = v_{i-1} + V_i, \quad w_i=w_{i-1}-(\nabla v_{i-1})^TV_i + W_i
\quad\mbox{ for  } i=1\ldots P,
\end{equation}
the final fields $v_P,w_P$ have regularity $\mathcal{C}^2 $ and satisfy:
\begin{equation*}
\begin{split}
& \lambda \|(v_P-v_0)\|_0 + \|\nabla (v_P-v_0)\|_0\leq C\delta^{1/2},
\\ & \lambda \|(w_P-w_0)\|_0 + \|\nabla (w_P-w_0)\|_0\leq C(\delta+\delta^{1/2}\|\nabla v_0\|_0), \\
& \|\nabla^{(2)}(v_P-v_0)\|_0\leq C\lambda \sigma^J\delta^{1/2}, \\
& \|\nabla^{(2)}(w_P-w_0)\|_0\leq C\lambda \sigma^J(\delta+\delta^{1/2}\|\nabla v_0\|_0), \\
&\big\|A_0 -\big(\frac{1}{2}(\nabla v_P)^T\nabla v_P +\sym\nabla
w_P\big)\big\|_0\leq C \frac{\delta}{\sigma^S}, 
\end{split}
\end{equation*}
with constants $C$ depending only on $\omega, d, k, J, S$.
Moreover, $\{(V_i,W_i)\}_{i=1}^P$ depend on the initial data only through $\mathcal{D}_0$ and $\nabla^{(2)} v_0$. 
The construction and its estimates use only their prescribed regularity and  estimates \eqref{bd}, but not the Hessian structure of $\nabla^{(2)}v_0$. 
\end{minipage}\quad \right]
\end{equation}

\medskip

\noindent Indeed, for system (\ref{VK}), the construction starts from:
\[\mathcal{D}_0= A_0-\Big(\frac12(\nabla v_0)^T\nabla v_0+\sym \nabla w_0\Big),\] which is the defect associated to the initial displacements $v_0, w_0$, and
it produces a finite number of corrugation pairs $\left( V_i,W_i\right)$ for $i=1\ldots P$. Defining:
\[v_i= v_{i-1}+V_i \qquad w_i = w_{i-1}-(\nabla v_{i-1})^TV_i+W_i\qquad \mbox{ for all }\, i=1\ldots P,\]
one obtains the final displacements $(v_P,w_P)$ with a quantitatively smaller defect and satisfying suitable estimates. All of the stage constructions for the system (\ref{VK}) have the listed additional property: the $P$-tuple $\{(V_i,W_i)\}_{i=1}^P$ depends on the initial data only through $\mathcal{D}_0$ and $\nabla^{(2)} v_0$. Moreover, the justification of the construction of $(V_i, W_i)$ does not require the Hessian structure of $\nabla^{(2)} v_0$, and it is well-defined if $\nabla^{(2)} v_0$ is replaced by any symmetric field $b_0$ satisfying the same estimates. This replacement is performed as follows. In the rules for the construction of the $i$-th perturbation, we write each current Hessian as: 
\[\nabla ^{(2)} v_{i-1} = \nabla^{(2)}v_0 +\sum_{j=1}^{i-1} \nabla^{(2)} V_j,\]
then modify the expression on the right hand side to:
\[[b_0^{pq}]_{pq} +\sum_{j=1}^{i-1} \nabla^{(2)} V_j,\]
serving as replacement input for the construction of the new correction pair $(V_i,W_i)$. We now show that  (\ref{STA-VK0}) implies (\ref{STA0}), with the construction as above:

\begin{theorem} \label{prop_STA1_to_STA0}
Assume that (\ref{STA-VK0}) holds, with some  $J, S\geq 1$. Then (\ref{STA0}) holds
 as well, with the same exponents $J,S$.
\end{theorem}
\begin{proof} 
For the original construction in \eqref{STA-VK0}, the update definition \eqref{fo} yields: 
\begin{equation*}
\begin{split}  
& \big(\frac{1}{2}(\nabla v_i)^T\nabla v_i +\sym\nabla w_i\big) - 
\big(\frac{1}{2}(\nabla v_{i-1})^T\nabla v_{i-1} +\sym\nabla w_{i-1}\big) 
\\ & = \frac{1}{2}(\nabla V_i)^T\nabla V_i+ \sym\nabla W_i + \sym \big((\nabla
v_{i-1})^T\nabla V_i\big) -\sym\nabla \big((\nabla v_{i-1})^TV_i\big)
\\ & = \frac{1}{2}(\nabla V_i)^T\nabla V_i+ \sym\nabla W_i -
\big[\langle \partial_{pq}v_{i-1}, V_i\rangle\big]_{pq},
\end{split}
\end{equation*}
so that the difference of the consecutive defects can be written as:
\begin{equation}\label{sei}
\begin{split}  
& \Big(A_0 - \big(\frac{1}{2}(\nabla v_i)^T\nabla v_i +\sym\nabla
w_i\big)\Big) - \Big(A_0 - \big(\frac{1}{2}(\nabla v_{i-1})^T\nabla v_{i-1} +\sym\nabla w_{i-1}\big)\Big) 
\\ & = - \Big(\frac{1}{2}(\nabla V_i)^T\nabla V_i+ \sym\nabla W_i\Big)
+ \big[\langle \partial_{pq}v_0, V_i\rangle\big]_{pq} +
\sum_{j=1}^{i-1}\big[\langle \partial_{pq}V_j, V_i\rangle\big]_{pq}.
\end{split}
\end{equation}
Let $A_0,b_0$, and $\delta,\lambda,\sigma $ satisfy the assumptions of \eqref{STA0}. We set $v_0=0$, $w_0=0$ and apply \eqref{STA-VK0} with initial defect $A_0$, replacing $\nabla^{(2)} v_0 $ (which equals zero in this case) by $b_0$. This replacement is permitted by the structural assumption,  since the assumptions of \eqref{STA0} give the required bounds on the initial defect $A_0$ and on $b_0$. Define the $i$-th defect for \eqref{STA0} by:
\[\tilde{\mathcal{D}}_i = A_0 + \left[\langle b_0^{pq},v_i\rangle\right]_{pq} -\Big(\frac12(\nabla v_i)^T\nabla v_i + \sym \nabla w_i\Big).\]
Then $\tilde{\mathcal{D}}_0 =A_0 $ and, since $v_0=0$, formula \eqref{fo} implies: 
\[ \tilde{\mathcal{D}}_i -\tilde{\mathcal{D}}_{i-1} = - \Big(\frac{1}{2}(\nabla V_i)^T\nabla V_i+ \sym\nabla W_i\Big)
+ \big[\langle b^{pq}_0, V_i\rangle\big]_{pq} +
\sum_{j=1}^{i-1}\big[\langle \partial_{pq}V_j, V_i\rangle\big]_{pq},\]
which is exactly \eqref{sei} with $\nabla^{(2)} v_0$ replaced by $b_0$. Since the initial defects agree,  the derived $\tilde{\mathcal{D}}_i$ consequently coincide with the defects propagated by the stage construction with $\nabla^{(2)} v_0$ being replaced by $b_0$.
By assumption, the same stage estimates hold. Since the initial fields vanish, the claim follows by taking $\tilde v= v_P$ and $\tilde w =w_P$.
\end{proof}

\section{Two examples}\label{sec_exa}

In this section, we show that (minor modifications of) the stage constructions in \cite{lew_conv} and \cite{CHIL} satisfy the structural properties required in (\ref{STA-VK0}).

\begin{example}\label{ex_first}
In \cite{lew_conv}, stage (\ref{STA-VK}) was proved in the
general case of arbitrary $d, k$, with:
$$ J=\frac{lcm(d_*,k)}{k}, \qquad  S=\frac{lcm(d_*,k)}{d_*}\qquad
\mbox{where } d_* = \frac{(d+1)d}{2},$$  
implying flexibility result in Corollary \ref{th_final_VK} with any exponent:
$$\alpha<\min\Big\{\frac{r+\beta}{2}, \frac{1}{1+2d_*/k}\Big\}.$$ 
In the same paper, (\ref{STA-VK}) was also deduced with:
$$ J=1, \qquad  S=N \qquad \mbox{ for any fixed } N\in\mathbb{N},
\mbox{ when } k\geq 2d_*,$$  
so that Corollary \ref{th_final_VK} holds in that case with exponents:
$$\alpha<\min\Big\{\frac{r+\beta}{2}, 1\Big\}.$$
We claim that both stage constructions in \cite{lew_conv} can be modified to yield (\ref{STA-VK0}), implying therefore (\ref{STA0}) and then (\ref{STA-II}). In conclusion, by Theorems \ref{prop_STA1_to_STA0} and
\ref{prop_STA0_to_STAII}, the assertion of Theorem \ref{th_final_II}
holds with the ranges of $\alpha$ indicated above. The same convex integration result for the isometric
immersion system (\ref{II}) was proved independently for $k=2d_*$ in
\cite{Kallen}, for $k=d$ in \cite{CaoSz2025}, and for arbitrary $d,k$
in \cite{LL_last}, whereas our present analysis allows to deduce the
same results from a more transparent construction for the simpler system (\ref{VK}).

\medskip

\noindent The constructions in \cite{lew_conv} already satisfy the structural properties required in \eqref{STA-VK0}, but need to be slightly modified at the level of the defect decomposition, since the smallness assumption on the defect is not centered around $H_0$, as in \eqref{bd}. To obtain a positive decomposition, a suitable multiple of $\Id_d$ is added to the defect and subsequently compensated by an affine shift in $w$; see \cite[section 3, (3.6), (3.12), (4.11)]{lew_conv}. In the first construction, where each stage consists of $S$ substages, such a shift is introduced at every substage. These large affine shifts prevent the additional factor of $\lambda^{-1}$ in the $\mathcal{C}^0$ bound on $w_P-w_0$, which is essential for transferring the von K\'arm\'an stage to the isometric immersion stage.
Under the centered defect assumption \eqref{bd}, the decomposition can instead be achieved after a simple rescaling, without adding a constant matrix which would require an affine shift to $w$. For a stage consisting of several substages, it remains to ensure that the intermediate defect has a centered bound at each substage. As explained below, this follows from the error estimates of \cite{lew_conv} with a slightly larger choice of frequencies.
In the following, we abbreviate:
\[\mathrm{VK}(v,w) = \frac{1}{2}\nabla v^T\nabla v+ \mathrm{sym}\nabla w\,,\quad \mathcal{D}_0 = A_0 - \mathrm{VK}(v_0,w_0)\,.\]

\smallskip

{\bf 1. (Construction for arbitrary $d,k$)}
For the first construction, the stage comprised of: 
\[P=lcm(d_*,k)=Jk=Sd_\star\]
steps, arranged into $S$ sub-stages, in which the intermediate defect $\mathcal{D}_s=A_0-\mathrm{VK}(v_{sd_*},w_{sd_*})$, $s=0\ldots S-1$, is first decomposed and then removed by the addition of $d_*$ perturbations. Each increment field $(V_i, W_i)$ in (\ref{fo}) is
given (see \cite[section 3, step 4]{lew_conv}) only in function of the
decomposition coefficients $a_{\delta}^s$ of the current intermediate defect, starting with the
initial defect $\mathcal D_0$. The next defect consists of the sum of
terms, of which only one depends on $\nabla^{(2)}v_{i-1}$, and hence on
$\nabla^{(2)}v_0$, in the linear fashion (see \cite[section 3, step 5, (3.13)]{lew_conv}).
Estimating the derivatives of the aforementioned defect depends thus only
on the initial bounds on $\nabla^{(2)}v_0$ in \cite[section 3, step 1,
(3.2)$_3$]{lew_conv}, utilized in \cite[section 3, steps 8, 9]{lew_conv}
to prove the inductive estimates in \cite[section 3, step 6, (3.14)$_2$,
(3.14)$_3$]{lew_conv}.
To verify that the decomposition of the defect without the affine shift is available at every substage, we write, for $s=0\ldots S$:
\[\mathcal D_s = A_0-\mathrm{VK}(v_{sd_*},w_{sd_*})\,,\quad \delta_s=\delta\sigma^{-s}\,. \]
We claim that we can guarantee the centered bound 
\begin{equation}\label{centereddef}
    \|\mathcal D_s - \delta_s H_0\|_0 \leq \frac{r_0}{2}\delta_s\,,\quad \mbox{for all }\, s=0\ldots S\,.
\end{equation} 
For $s=0$, this is exactly the assumption \eqref{bd}. At substage $s$, we apply Lemma \ref{lem_dec_def} to: 
\[ H_s = \frac{\mathcal D_s}{\delta_s} - \frac{\delta_{s+1}}{\delta_s} H_0 =\frac{\mathcal D_s}{\delta_s} - \frac{1}{\sigma} H_0\,. \]
and note that for $\sigma_0 $ sufficiently large, the centered bound implies $\|H_s-H_0\|_0\leq 3r_0/4$.
Decomposing $H_s$ and rescaling the coefficients therefore gives: 
\[ \mathcal D_s-\delta_{s+1}H_0 = \delta_s H_s = \sum_{i=1}^{d_*} (\delta_s^{1/2}a^s_i)^2\eta_i\otimes\eta_i\quad\mbox{ for }\quad 
a_i^s = \bar a_i\left(H_s\right)^{1/2}\,.\]
We then use the amplitudes $\{\delta_s^{1/2}a_i^s\}_{i=1}^{d_*}$ in the perturbation formulas of \cite[section 3, step 4]{lew_conv} and recover the new defect which satisfies:
\[\mathcal D_{s+1}-\delta_{s+1}H_0 = -E_{s+1}\,, \]
 where $E_{s+1}$ is the sum of error terms in \cite[(3.13)]{lew_conv}. The required smallness of $E_{s+1}$ follows from the same inductive estimates as in \cite[(3.14)]{lew_conv}, after the following enlargement of frequencies. Parameters $\lambda, \sigma$ of \eqref{STA-VK0} correspond, respectively, to $l^{-1}$ and the ratio $\lambda l$ of \cite{lew_conv}. Writing:
\[ i = jk +\gamma = sd_*+l,\quad 1\leq \gamma\leq k,\quad 1\leq l\leq d_*\,, \] 
we define the frequencies:
\begin{equation}\label{freq}
    \lambda_0= \lambda\,,\quad \lambda_{i } = K^i \lambda \sigma^{1+j+s/2}\,,
\end{equation}
for some large enough factor $K\geq 1$,  depending only on $\omega, d,k,J,S$. The estimates of \cite[section 3, step 9]{lew_conv} then yield: 
\[\|E_{s+1}\|_0\leq C\left( K^{-1}+K^{-2}\right) \delta_{s+1}\,.\]
Here, the error terms involving the current Hessian (the first term on the right hand side of \cite[(3.13)]{lew_conv}) contribute the factor $K^{-1}$,  the quadratic terms the factor $K^{-2}$. The corresponding derivative estimates propagate as in \cite{lew_conv}. Since the number of substages is fixed, $K$ can be chosen independently of $\delta,\lambda,\sigma$ so that \eqref{centereddef} holds. At the end of the construction, this yields:
 \[\|A_0-\mathrm{VK}(v_P,w_P)\|_0\leq \delta_S |H_0|+\|\mathcal{D}_{S}-\delta_SH_0\|_0\leq \left(|H_0|+\frac{r_0}{2}\right)\frac{\delta}{\sigma^S},\]
 as desired. 
Finally, the perturbation formulas, together with $\|\nabla v_i\|_0\leq \|\nabla v_0\| + C\delta^{1/2}$, give: 
\[\|w_P-w_0\|_0\leq \frac{C}{\lambda}\left( \delta+\delta^{1/2}\|\nabla v_0\|_0\right)\,.\] The remaining estimates follow as in \cite[section 3, steps 7,8]{lew_conv}. The additional fixed factor $K$ affects only the constants, since: 
 \[\lambda_i\delta_s^{1/2}= K^i\lambda \sigma^{1+j}\delta^{1/2}\leq K^P\lambda\sigma^J\delta^{1/2}\,.\]
Hence, the exponents $J,S$ are unchanged. 

\medskip

{\bf 2. (Construction for  $k\geq 2d_*$)} The second construction in \cite{lew_conv} also has the structural properties required in \eqref{STA-VK0}.  The related stage comprises of utilizing $d_*$ simultaneous Nash spirals, equivalent to $P=2$ increment pairs in (\ref{fo}). The increments $\{(V_i, W_i)\}_{i=1}^2$ are then given
in function of the decomposition coefficients of the $N$-th
K\"allen iteration, in which $\nabla^{(2)}v_0$ appears linearly, see 
\cite[section 4, steps 2, 4]{lew_conv}, and the inductive estimates
\cite[section 3, steps 3, (4.7)$_1$ - (4.7)$_3$]{lew_conv} depend only
on the initial bounds on $\nabla^{(2)}v_0$ in \cite[section 3, step 1,
(3.2)$_3$]{lew_conv}. Thus the construction and its estimates remain valid when $\nabla^{(2)} v_0$ is replaced by an arbitrary symmetric tensor $b_0$.
As in the first construction, the centered defect assumption can be used to omit the affine shift in $w$. In this case, the decomposition of the defect is obtained by a K\"allen iteration, see \cite[section 4, step 2]{lew_conv}. The defect assumption allows the omission of the term $\tilde C\,\mathrm{Id}_d$ from the amplitude definition in \cite[section 4, step 2]{lew_conv}, and, in consequence, in the formulas \cite[(4.3)]{lew_conv} and \cite[(4.11)]{lew_conv}. As above, this leads to the sharper increment estimates, required in \eqref{STA-VK0}.
\endproof
\end{example}  

\medskip

\begin{example}\label{ex_last}
In \cite{CHIL}, another construction led to a proof of (\ref{STA-VK}) with:
$$ J=dK+N, \qquad  S=KN \qquad
\mbox{ for any fixed } K,N\in\mathbb{N} \mbox{ when } k\geq d_*-d+1,$$  
so that Corollary \ref{th_final_VK} holds in that dimensionality case with the exponents:
$$\alpha<\min\Big\{\frac{r+\beta}{2}, 1\Big\}.$$
The same paper states the corresponding flexibility result
 for Euclidean isometric immersions in \cite[Theorem 5.1]{CHIL} and outlines its proof. We now
argue that the stage construction in \cite{CHIL} can be easily modified to yield
(\ref{STA-VK0}), so that the claim results from Theorems \ref{prop_STA1_to_STA0},
\ref{prop_STA0_to_STAII}.

\medskip

\noindent Indeed, the stage in \cite{CHIL} comprises of $P=Kd_*$ steps,
arranged into $K$ sub-stages, in each of which the $d_*$ steps correspond to the addition of the increment fields $(V_i, W_i)$ as in (\ref{fo}). The  structural properties of the stage are already those required by \eqref{STA-VK0}. In the first $d$ steps of each substage,
the increments $V_i$ are collinear with
the first codimension vector $e_1$ and depend only on the 
decomposition coefficients $a_{1j}$ of the appropriate
generation of defect, obtained after $N$ applications of K\"allen iterations,  see formula \cite[(3.7)]{CHIL}.  The increments $W_i$ additionally contain the integrations by parts corrections of \cite[(3.8)]{CHIL}, whose dependence on the current Hessian $\nabla^{(2)}v^1_{i-1}$ is linear. Moreover, \cite[Lemma 2.4]{CHIL} applies to arbitrary symmetric matrix fields. Hence these corrections and the error fields in \cite[(3.10), (3.11)]{CHIL} do not require the initial tensor to be a Hessian. 
In the remaining $d_*-d$ steps of the substage, the increments depend on the current coefficients $b_{ij}$ (defined in \cite[section 3, step 8]{CHIL}) and those from preceding substages through \cite[(3.28), (3.29)]{CHIL}. The expansion in \cite[(3.35), (3.36)]{CHIL} separates out the initial Hessian, and its sizes and derivatives are estimated using only the initial bounds in \cite[(3.1)]{CHIL}. Consequently, the construction and its estimates remain valid under the replacement of $\nabla^{(2)} v_0$ by an arbitrary symmetric table $[b_0^{pq}]_{pq}$ satisfying the same estimates.  


\medskip

\noindent The adjustment to the defect decomposition follows Example \ref{ex_first}. In the construction of \cite{CHIL}, the constant matrix added in \cite[(3.5)]{CHIL} is compensated by the affine term in $\bar\Theta_{11}$ in \cite[(3.8)]{CHIL}. To omit this term in every substage, we again set, for $s=0\ldots K$:
\[\mathcal D_s = A_0-\mathrm{VK}(v_{sd_*},w_{sd_*})\,,\quad \delta_s=\delta\sigma^{-Ns}, \]
and guarantee the centered bound: 
\begin{equation}\label{centereddef2}
    \|\mathcal D_s - \delta_s H_0\|_0 \leq \frac{r_0}{2}\delta_s\quad\mbox{ for all } \, s=0\ldots K.
\end{equation} 
 This is done by decomposing $H_s= \delta_s^{-1}\left(\mathcal D_s -\delta_{s+1}H_0\right)$ with \cite[Lemma 2.5]{CHIL}. Notice that the latter lemma is stated with radius $r_0/2$. However, the K\"allen iteration also works with radius $3r_0/4$ (possibly after a frequency increase). Thus, applying this decomposition to $H_s$ and rescaling the resulting coefficients by $\delta_s^{1/2}$ yields the required decomposition of $\mathcal D_s -\delta_{s+1}H_0$ without an added constant matrix. We can therefore replace $\bar\Theta_{11}$ by $\Theta_{11}$ in \cite[(3.8)]{CHIL}.
We then replace the factor $\sigma$ in the definition of  frequencies in \cite[(3.2)]{CHIL} by $\Lambda \sigma$, for a suitably large $\Lambda\geq 1$. The error estimates in
\cite[(3.6), (3.22), (3.34), (3.37)]{CHIL} yield:
\[\|\mathcal D_{s+1}-\delta_{s+1}H_0\|_0 \leq C\frac{\delta_s}{(\Lambda\sigma)^N} = C\Lambda^{-N} \delta_{s+1}\leq \frac{r_0}{2}\delta_{s+1},\]
and the corresponding derivative estimates follow as in \cite{CHIL}, observing that: 
\[\frac{\delta_s}{\delta_{s+1}}=\sigma^N\leq (\Lambda\sigma)^N,\]
so that \cite[(3.16), (3.21)]{CHIL} hold. 
Finally, estimates in \cite[section 3, steps 6 and 11]{CHIL}, with the affine term removed, give the desired bound: 
\[\|w_P-w_0\|_0\leq\frac{C}{\lambda}\big(\delta+\delta^{1/2}\|\nabla v_0\|_0\big).\]
The remaining increment estimates follow from the same argument. As in Example \ref{ex_first}, the enlargement of the frequencies only changes constants, not the exponents $J, S$.
\endproof
\end{example} 

\section*{Appendix: proofs of \ref{stima3.5}, \ref{stima5} and Lemmas \ref{lem_det}, \ref{lem_Tu}, \ref{lem_propa}}\label{sec_app}

In this section, we prove the new mollification and triple commutator estimates. We also show how to propagate a normal frame from a given tangent frame to a nearby tangent frame, without requiring the corresponding ambient metrics $\Id$ and $\bar G$, respectively, to be close.

\bigskip

\noindent {\bf Proof of \ref{stima3.5} in Lemma \ref{lem_stima}}

Write the estimated quantity as the sum of two terms in:
\begin{equation*}
\begin{split}
& \Big|\big((\bar f\circ u)\ast\phi_l - \bar f\circ u \big)(x) \Big|\leq
\int_{\R^d} \Big|(\bar f\circ u)(x-y) - (\bar f\circ u)(x) - \langle (\nabla\bar f)_{\mid u(x)} \nabla u(x), y\rangle \Big|\phi_l(y)\;\mbox{d}y
\\ & \leq \int_{\R^d} \Big|(\bar f\circ u)(x-y) - (\bar f\circ u)(x) - \langle (\nabla\bar f)_{\mid u(x)}, u(x-y)-u(x)\rangle \Big|\phi_l(y)\;\mbox{d}y
\\ & \quad + \int_{\R^d} \Big|\langle (\nabla\bar f)_{\mid u(x)}, u(x-y) - u(x)+\nabla u(x)y\rangle \Big|\phi_l(y)\;\mbox{d}y.
\end{split}
\end{equation*}
The integrand in the first term in the right hand side above may be written as:
\begin{equation*}
\begin{split}
\Big|\int_0^1 \langle (\nabla\bar f)_{\mid t u(x-y) + (1-t)u(x)}-(\nabla \bar f)_{\mid u(x)}, u(x-y) - u(x)\rangle \;\mbox{d}t\Big| & \leq \|\nabla \bar f\|_{0,\beta}|u(x-y)-u(x)|^{1+\beta} \\ & \leq \|\nabla \bar f\|_{0,\beta}\|\nabla u\|_0^{1+\beta} |y|^{1+\beta},
\end{split}
\end{equation*}
whereas the integrand in the second term is estimated by:
\begin{equation*}
\begin{split}
\Big|\langle (\nabla\bar f)_{\mid u(x)}, u(x-y) - u(x)+\nabla u(x)y\rangle \Big| & \leq \frac{1}{2} \|\nabla \bar f\|_{0}\|\nabla^2 u\|_0 |y|^2,
\end{split}
\end{equation*}
Since both these terms are integrated against  $\phi_l(y)$ on $\R^d$, there follows \ref{stima3.5}.
\endproof

\bigskip

\noindent {\bf Proof of \ref{stima5} in Lemma \ref{lem_stima}}

Write the triple commutator as the sum of two commutator terms in:
\begin{equation*}
\begin{split}
& (fgh)\ast\phi_l - (f\ast\phi_l)
(g\ast\phi_l)(h\ast\phi_l) \\ & = \Big((fgh)\ast\phi_l - ((fg)\ast\phi_l)
(h\ast\phi_l)\Big) + \Big((fg)\ast\phi_l - 
(f\ast\phi_l)(g\ast\phi_l)\Big)(h\ast\phi_l)
\end{split}
\end{equation*}
By \ref{stima4} and \ref{stima1} we get:
\begin{equation*}
\begin{split}
& \big\|\nabla^{(m)}\big((fgh)\ast\phi_l - (f\ast\phi_l)
(g\ast\phi_l)(h\ast\phi_l)\big)\big\|_0 \\ & \leq  C l^{2-m} 
\|\nabla(fg)\|_0\|\nabla h\|_0 + C\sum_{r+s=m}
l^{2-r}\|\nabla f\|_0\|\nabla g\|_0\|\nabla^{(s)}(h\ast\phi_l)\|_0
\\ & \leq C l^{2-m}\big(\|f\|_0\|\nabla g\|_{0} \|\nabla h\|_{0}
+\|g\|_0\|\nabla f\|_{0} \|\nabla h\|_{0}\big) + 
C \|\nabla f\|_0\|\nabla g\|_0 \sum_{r+s=m}
l^{2-r}l^{-s}\|h\|_0,
\end{split}
\end{equation*}
which yields \ref{stima5}. The proof is done.
\endproof

\bigskip

\noindent {\bf Proof of Lemma \ref{lem_det}}

The first assertion holds, as for every $i=1\ldots d$, conditions (\ref{assu_imm}) imply:
$$\langle Fe_i, GFe_i\rangle \in \Big[\frac{1}{\gamma},\gamma\Big]\cap \Big[\frac{|Fe_i|^2}{\Lambda}, \Lambda|Fe_i|^2\Big].$$
Consequently, $|Fe_i|^2\in \big[{1}/{\gamma\Lambda}, \gamma\Lambda\big]$, so $|F|^2\in \big[{d}/{\gamma\Lambda}, d\gamma\Lambda\big]$. For the second assertion, note that all the eigenvalues of the symmetric matrix $F^TGF$ lie within the interval $[1/\gamma,\gamma]$. This implies the stated determinant bound.
\endproof

\bigskip

\noindent {\bf Proof of Lemma \ref{lem_Tu}}

Writing $T=\frac{1}{\det(F^TGF)} F\,\mbox{cof}\big(F^TGF\big)$, we easily conclude that
Lemma \ref{lem_det} implies the first assertion, since $\|G\|_0\leq C$ depending on $\Lambda, d, k$. For the second assertion, we use Fa\'a di Bruno's formula
to the composition $T =f\circ (F,G)$, where the map
$f(F, G) = F\big(F^TGF\big)^{-1}$
is smooth and has all its derivatives bounded on the compact
domain given by the condition (\ref{assu_imm}).
\endproof

\bigskip

\noindent {\bf Proof of Lemma \ref{lem_propa}}

{\bf 1.} Fix a unit vector $\xi\in\mathbb{S}^{d-1}$. Reasoning as in Lemma \ref{lem_det}, the first assumption in (\ref{assu_imm2}) yields $|F\xi|\in\big[\gamma^{-1/2}, \gamma^{1/2}\big]$. Since $\big||\tilde F\xi|-|F\xi|\big|\leq |(\tilde F-F)\xi|\leq \|\tilde F-F\|_0\leq \rho$, it follows that $|\tilde F\xi|\in \big[ \gamma^{-1/2} - \rho, \gamma^{1/2} + \rho\big]$, and consequently by the second assumption in (\ref{assu_imm2}) we get:
$$\langle \tilde F^T\tilde G\tilde F\xi, \xi\rangle  = 
\langle \tilde G\tilde F\xi, \tilde F\xi\rangle \in \Big[ \frac{(\gamma^{-1/2}-\rho)^2}{\Lambda}, \Lambda (\gamma^{1/2}+\rho)^2\Big].$$
This proves (\ref{out1}) if only $\rho$ is sufficiently small, in function of $\gamma$.

\smallskip

{\bf 2.} To prove the second assertion, observe first that (\ref{out1}) and Lemma \ref{lem_Tu} yield:
$$\|\tilde T\|_0\leq C \quad \mbox{ where }\; \tilde T = \tilde F(\tilde F^T\tilde G\tilde F)^{-1}\in \mathcal{C}^n(\bar\omega, \R^{(d+k)\times d}),$$
with $C$ that depends on $\gamma, \Lambda, d,k$.
Set $E^i=Ee_i$ for $i=1\ldots k$ and define:
$$\nu^i= \tilde G^{-1}E^i-\tilde T (\tilde F-F)^TE^i.$$
We claim that vectors $\{\nu^i\in\R^{d+k}\}_{i=1}^k$ are $\tilde G$-uniformly linearly independent in $\bar\omega$. Then, $\tilde E$ is defined by setting $\tilde Ee_i = \tilde E^i$ and employing an inductive Gram-Schmidt orthonormalization in:
\begin{equation}\label{GSch}
\begin{split}
&\tilde E^i = \frac{\tilde \nu^i}{\langle \tilde \nu^i, \tilde G\tilde\nu^i\rangle^{1/2}} \quad\mbox{ for all }\; i=1\ldots k,
\\ & \mbox{where: }\; \tilde \nu^i = \nu^i-\sum_{j<i}\langle \nu^i, \tilde G\tilde E^j\rangle \tilde E^j.
\end{split}
\end{equation}
To show the aforementioned linear independence, fix  $\xi\in\mathbb{S}^{k-1}$. The second assumption in (\ref{assu_imm2}) yields the same ellipticity bound for $\tilde G^{-1}$, so the same argument as in Lemma \ref{lem_det} implies:
$$|\tilde G^{-1}E\xi|\in \Big[\frac{|E\xi|}{\Lambda^{1/2}}, \Lambda^{1/2}|E\xi|\Big] = \Big[\frac{1}{\Lambda^{1/2}}, \Lambda^{1/2}\Big], $$
because $|E\xi|=1$. Further, denote $\nu\in\mathcal{C}^n(\bar\omega, \R^{(d+k)\times k})$  where $\nu e_i = \nu^i$, so that:
$$\big||\nu\xi|-|\tilde G^{-1}E\xi|\big| \leq |\tilde T (\tilde F-F)^TE|\leq C \|\tilde T\|_0 \|\tilde F-F\|_0\leq C\rho$$
with $C$ depending on $\gamma,\Lambda,d,k$. Consequently, we get:
$|\nu\xi|\in \big[(2\Lambda)^{-1/2}, (2\Lambda)^{1/2}\big]$ for $\rho$ sufficiently small, and further: $\langle \nu^T\tilde G\nu\xi, \xi\rangle \in \big[2\Lambda^{2}, 2\Lambda^{2}\big]$, resulting in the claimed condition:
$$\frac{1}{2\Lambda^2}\Id_k\leq \nu^T\tilde G\nu \leq 2\Lambda^2\Id_k \quad \mbox{ in }\; \bar\omega.$$
The bound (\ref{out2}) follows now directly by the Fa\`a di Bruno formula, since each $\tilde E^i$ is a smooth function of the quantities $\tilde G, \tilde F, F$, $E^1,\ldots, E^i$, on the compact region of their definiteness.
\endproof

\section*{AI declaration} Assistance of ChatGPT was used in discussions leading to  Example \ref{counter_ex}. The authors take full responsibility for the content of the work.

\end{document}